\documentclass[mnsc,nonblindrev]{informs3_arxiv} 

\OneAndAHalfSpacedXI 

\usepackage{natbib}
 \bibpunct[, ]{(}{)}{,}{a}{}{,}%
 \def\bibfont{\small}%
 \def\bibhang{24pt}%

\usepackage{graphicx}
\usepackage{caption}
\usepackage{subcaption}
\usepackage{amsmath}
\usepackage{amsfonts}
\usepackage{amssymb}
\usepackage{amsbsy}
\usepackage{multirow}
\usepackage[ruled,vlined]{algorithm2e}
\usepackage[pdftex,colorlinks=true,urlcolor=blue,citecolor=black,anchorcolor=black,linkcolor=black]{hyperref}
\usepackage{xcolor}

\TheoremsNumberedThrough     
\ECRepeatTheorems

\EquationsNumberedThrough    

\MANUSCRIPTNO{}

\newcommand{\R}{\mathbb{R}}
\newcommand{\hp}{\hat{p}}
\newcommand{\tP}{\tilde{P}}
\newcommand{\tp}{\tilde{p}}
\newcommand{\tE}{\tilde{E}}
\newcommand{\tVar}{\widetilde{Var}}
\newcommand{\tf}{\tilde{f}}
\newcommand{\cA}{\mathcal{A}}
\newcommand{\cE}{\mathcal{E}}
\newcommand{\cH}{\mathcal{H}}
\newcommand{\bZ}{\mathbf{Z}}
\newcommand{\bY}{\mathbf{Y}}

\begin{document}
\RUNAUTHOR{Bai et al.}

\RUNTITLE{Probabilistic Efficiency in Rare-Event Simulation}

\TITLE{Over-Conservativeness of Variance-Based Efficiency Criteria and Probabilistic Efficiency in Rare-Event Simulation}

\ARTICLEAUTHORS{%
\AUTHOR{Yuanlu Bai}
\AFF{Columbia University, \EMAIL{yb2436@columbia.edu}}
\AUTHOR{Zhiyuan Huang}
\AFF{Tongji University, \EMAIL{huangzy@tongji.edu.cn}}
\AUTHOR{Henry Lam}
\AFF{Columbia University, \EMAIL{henry.lam@columbia.edu}}
\AUTHOR{Ding Zhao}
\AFF{Carnegie Mellon University, \EMAIL{dingzhao@cmu.edu}}
} 

\ABSTRACT{%
In rare-event simulation, an importance sampling (IS) estimator is regarded as efficient if its relative error, namely the ratio between its standard deviation and mean, is sufficiently controlled. It is widely known that when a rare-event set contains multiple ``important regions" encoded by the so-called dominating points, IS needs to account for all of them via mixing to achieve efficiency. We argue that in typical experiments, missing less significant dominating points may not necessarily cause inefficiency, and the traditional analysis recipe could suffer from intrinsic looseness by using relative error, or in turn estimation variance, as an efficiency criterion. We propose a new efficiency notion, which we call \emph{probabilistic efficiency}, to tighten this gap. In particular, we show that under the standard Gartner-Ellis large deviations regime, an IS that uses only the most significant dominating points is sufficient to attain this efficiency notion. Our finding is especially relevant in high-dimensional settings where the computational effort to locate all dominating points is enormous.
}%



\maketitle

\section{Introduction}
\label{sec:intro}

We study the problem of estimating the probabilities of rare events with Monte Carlo simulation, which falls in the domain of rare-event simulation \citep{bucklew2004rare,juneja2006rare,rubino2009rare}. Traditionally, rare-event simulation is of wide interest to a variety of areas such as queueing systems \citep{dupuis2007dynamic,dupuis2009importance,blanchet2007editorial,blanchet2009rare,blanchet2014rare,kroese1999efficient,ridder2009importance,sadowsky1991large,szechtman2002rare}, highly dependable computer systems and communication networks \citep{tuffin2004numerical,lewis1984monte,goyal1989unified,carrasco1992failure,shahabuddin1994importance,kesidis1993effective}, financial risk management \citep{glasserman2003monte,glasserman2005importance,glasserman2008fast} and insurance modeling \citep{asmussen1985conjuate,asmussen2010ruin}. More recently, with the rapid development of intelligent physical systems such as autonomous vehicles and personal assistive robots \citep{ding2021multimodal, arief2021deep}, rare-event simulation is also applied to assess their risks before deployments in public, where the risks are often quantified by the probabilities of violations of certain safety metrics such as crash or injury rate \citep{huang2018accelerated,o2018scalable,zhao2016accelerated, zhao2018accelerated}. The latter problems typically involve complex AI-driven underlying algorithms that deem the rare-event structures rough or difficult. The current work is motivated from the importance of handling such type of rare-event problems (e.g., the U.S. National Artificial Intelligence Research and Development Strategic Plan \citep{national2019national} lists ``developing effective evaluation methods for AI" as a top priority) and provides a step towards rigorously grounded procedures in this direction.

The starting challenge in rare-event simulation is that, by its own nature, the target rare events seldom occur in the simulation experiment when using crude Monte Carlo. In other words, to achieve an acceptable estimation accuracy relative to the target probability, the required simulation size could be huge in order to obtain sufficient hits on the target events. Statistically, this issue is manifested as a large ratio between the standard deviation (per run) to the mean, known as the relative error, that determines the order of a required sample size. In the large deviations regime where the target probability can depend exponentially on the rarity parameter, this in particular means the required sample size is exponentially large.

To address the inefficiency of crude Monte Carlo, a range of variance reduction techniques have been developed. Among them, importance sampling (IS) \citep{siegmund1976importance} has been broadly applied to improve the efficiency. IS uses an alternative probability measure to generate the simulation samples, and then reweighs the outputs via likelihood ratio to guarantee unbiasedness. The goal is that by using this alternate estimator than simply counting the frequency of hits in crude Monte Carlo, one can achieve a small relative error with a much smaller sample size.

To this end, it is also widely known that IS is a ``delicate" technique, in the sense that the IS probability measure needs to be carefully chosen in order to achieve a small relative error. In the typical large deviations setting, the suggestion is to tilt the probability measure to the ``important region". The delicacy appears when there are more than one important regions, in which case all of them need to be accounted for. More specifically, in the light-tailed regime, these important regions are guided precisely by the so-called dominating points, which capture the most likely scenario in a local region of the rare event. Despite the tempting approach to simply shift the distribution center to the globally most likely scenario, it is well established that if not all the dominating points are included in the IS mixture distribution, then the resulting estimator may no longer be efficient in terms of the relative error (see, for instance, the seminal work \cite{glasserman1997counterexamples}).

Our main goal in this paper is to argue that, in potentially many light-tailed problems, the inclusion of all the dominating points in an IS could be \emph{unnecessary}. Our study is motivated from high-dimensional settings where finding all dominating points could be computationally expensive or even prohibitive, yet these problems may arise in recent safety-critical applications (e.g., \cite{arief2021deep,webb2018statistical,bai2022rare}).

To intuitively explain the unnecessity, let us first drill into why all the dominating points are arguably needed in the literature in the first place. Imagine that a rare event set $\cE$ comprises two disjoint ``important" regions, say $\cE_1$ and $\cE_2$, and the dominating points are correspondingly $a_1$ and $a_2$, which are sufficiently ``far away'' from each other, and $a_1$ has say a higher density than $a_2$. Roughly speaking, if an IS scheme only focuses on $\cE_1$ and tilts the distribution center towards $a_1$, then there is a small chance that the sample from this IS distribution hits $\cE_2$, so that the contribution from this sample in the ultimate estimator is non-zero and, moreover, may constitute a large likelihood ratio and consequently elicit a large variance. This unfortunate event of falling into a secondary important region is a source of inefficiency according to the relative error criterion.

Now let's take a step back and think about the following: How likely does the ``unfortunate" event above occur? In a typical Monte Carlo experiment, we argue -- and we will see clearly in experiments -- that this could be very unlikely, to the extent that we shouldn't be worried at all with a reasonable simulation size. Yet, according to the relative error criterion, it seems necessary to worry about this, because it contributes to the variance of the estimator in each single run. This points to that using variance to measure efficiency in rare-event simulation could be \emph{too loose} to begin with. This variance measure, in turn, comes from the Markov inequality that converts relative error into a sufficient relative closeness between the estimate and target probability with high confidence. In other words, this Markov inequality itself could be the source of looseness.

This motivates us to propose what we call \emph{probabilistic efficiency}. Different from all the efficiency criteria in the literature, including asymptotic efficiency (also known as asymptotic optimality or logarithmic efficiency) and bounded relative error \citep{juneja2006rare,ecuyer2010asymptotic}, probabilistic efficiency does not use relative error. Instead, it is a criterion on achieving relative closeness directly. A distinctive element in probabilistic efficiency is the control on the simulation size itself, that we only allow it to grow \emph{moderately} with the underlying rarity parameter. This moderate simulation size, which is often the only feasible option in experiments, suppresses the occurrence of the unfortunate event of falling into a secondary important region. This way, while the variance could blow up, the high-confidence closeness between the estimate and target probability could still be retained.

With this new framework, we show that under standard assumptions in the widely used Gartner-Ellis large deviations regime, an IS that uses only the most significant dominating points is sufficient to ensure probabilistic efficiency. The Gartner-Ellis regime has been used across different applications such as queueing \citep{ridder2009importance,szechtman2002rare,blanchet2014rare}, communication systems \citep{smith1997quick,chen1993importance} and finance \citep{zhang2009rare,blanchet2009efficient,glasserman2005importance,glasserman2008fast}. Our results thus stipulate that in all these problems, in order to obtain a good estimate relative to the ground-truth rare-event probability, we only need to exponentially tilt to the most significant dominating point when there is only one such point, without the use of any mixture. This is a sharp contrast to the established IS recipe. Moreover, this makes the construction of IS in closer line with the large deviations theory that governs the rare-event probability asymptotic. More specifically, in large deviations, the most significant dominating point coincides with the minimizer of the so-called rate function that controls the exponential decay rate of the probability. Our theory thus postulates that to attain probabilistic efficiency, it suffices to consider only this rate function minimizer when constructing the IS.


We close this introduction with further discussions of our study in relation to existing works. First, we contrast our probabilistic efficiency with the notion of ``well-estimated" in the dependability literature (Definition 3 in \cite{tuffin2004numerical}). The latter asserts that all the paths which contribute most significantly to the target quantity need to occur sufficiently likely under the IS distribution. This notion is intuitively similar to the scheme of tilting the most significant dominating points which we suggest, but there are substantial differences in the motivations, implications, and also some technicality. First is that \cite{tuffin2004numerical} uses well-estimated as a sufficient condition to guarantee the efficiency in estimating the variance, not only the rare-event probability itself. In contrast, our work proposes to remove variance as an efficiency criterion, a motive that is essentially orthogonal to \cite{tuffin2004numerical}, and we propose the tilting of only the most significant dominating points as being sufficient to ensure our new criterion of probabilistic efficiency. Our key message is that the existing proposal that suggests using all dominating points, \emph{not} only the most significant ones, is an over-conservative approach and hence worth our remedy. Regarding technicality, \cite{tuffin2004numerical} focuses on dependable systems where the target probability decays polynomially instead of exponentially in the rarity parameter as in our setting. Under his framework, there are only countably many points in the sample space, and thus the target probability or variance can be more readily decomposed (i.e., one could explicitly quantify the ``contribution'' of each point), while in our possibly continuous setting the contribution of each point is less easy to determine. Lastly, well-estimated requires that all the significant paths are no longer rare (i.e., of order $O(1)$) under IS, while we require a less stringent subexponential decay. 

Second, we caution that probabilistic efficiency is not meant to replace existing variance-based efficiency criteria, but rather to \emph{complement} them especially in situations where identifying all dominating points is infeasible due to problem complexity. In problems where the latter is not an issue, it remains ``safer" to use existing criteria, as probabilistic efficiency relies on a more subtle sample size requirement, namely that it is not exponentially large. While this condition is reasonable for realistic problems, it would warrant future experimental diagnostics to detect violations of such a condition. 
Third, we note that our variance-free approach can potentially be adapted to rare-event estimation problems to streamline IS construction beyond the considered light-tailed large deviations regime. However, for some of these problems (e.g., heavy-tailed problems;  \citealt{blanchet2008efficient,hult2012importance,blanchet2012efficient,chen2019efficient}), the efficiency gain appears less dramatic than ours which possesses an exponential speed-up. Moreover, the Gartner-Ellis paradigm that we consider in this paper is arguably the most widely used and forms the basis of analysis for many rare-event problems.


In the following, we first introduce in more detail the background of rare-event simulation and the established efficiency criteria in the literature, all of which involve estimation variance or relative error (Section \ref{sec:background}). Then we show several motivating numerical examples to illustrate how excluding some dominating points in IS appear to give similar and sometimes even better performances than including all these points, the latter suggested predominantly in the literature (Section \ref{sec:prelim numerics}). This motivates our new notion of probabilistic efficiency to explain the observed numerical phenomena (Section \ref{sec:concept}) and the analysis of efficiency guarantees using this new notion (Section \ref{sec:guarantees_GE}). We then show further numerical results to validate our theory and performances of our estimators (Section \ref{sec:additional numerics}). Finally, we give some cautionary notes about probabilistic efficiency which involve the risk of under-estimation, and suggest some future directions (Section \ref{sec:conclusion}).

\section{Problem Setting and Existing Framework}\label{sec:background}
We consider an indexed family of rare events $\{\cA_{\gamma}\}_{\gamma}$, where $\gamma$ denotes a ``rarity parameter" such that as $\gamma\to\infty$, the event $\cA_{\gamma}$ becomes rarer so that $P(\cA_\gamma)\to0$. Our goal is to estimate $p=p(\gamma):=P(\cA_{\gamma})$ using Monte Carlo simulation. Here, the index $\gamma$ is introduced for modeling purpose so that we can speak of asymptotic rate, which is customary in the rare-event simulation literature. 
\subsection{Crude Monte Carlo and Relative Error}
To motivate the various notions that we would discuss momentarily, let us consider using crude Monte Carlo to estimate $p$. This means we utilize the unbiased estimator $Z=I_{\cA_{\gamma}}$ for $p$, where $I_{\cA_\gamma}$ denotes the indicator variable on the event $\cA_\gamma$. Suppose we generate the output $Z$ independently $n$ times, and construct $\hp$ as their sample mean. Intuitively, when $p$ is tiny, this estimator $\hp$ is most likely zero unless $n$ is a huge number, since a long trial length is needed to land at the rare event $\cA_\gamma$.

To describe the above challenge mathematically, we consider the following criterion. For a given tolerance level $\varepsilon>0$ (e.g., $5\%$), we would like an estimator $\hat p$ to satisfy
\begin{equation}
P(|\hat p-p|>\delta p)\leq \varepsilon\label{basic}
\end{equation}
for a certain $0<\delta<1$ when using a simulation size $n$. In \eqref{basic}, the closeness between $\hat p$ and $p$, which represents the error of $\hat p$ in estimating $p$, is measured relative to the magnitude of $p$ itself. This is because for tiny $p$, the estimation is only meaningful if the error is small enough relative to this tiny quantity. Linking to crude Monte Carlo, outputting merely $\hat p=0$, as likely to happen thereby, would be viewed as incurring a substantial error, i.e., lying in the event $|\hat p-p|>\delta p$ in \eqref{basic}. 
Put in another way, we argue that a huge sample size $n$ is needed for crude Monte Carlo to attain \eqref{basic}. Note that by Chebyshev's inequality, we get that for any $0<\delta<1$,
\begin{equation}
    P(|\hp-p|>\delta p)\leq \frac{Var(Z)}{n\delta^2p^2}.\label{Chebyshev}
\end{equation}
Hence, $n\geq \frac{Var(Z)}{\varepsilon\delta^2p^2}$ implies that $P(|\hp-p|>\delta p)\leq\varepsilon$ for $\varepsilon>0$. This means that $\frac{Var(Z)}{\varepsilon\delta^2p^2}$ is a sufficient size for $n$ to achieve \eqref{basic}. This quantity depends on the ratio between the standard deviation $\sqrt{Var(Z)}$ and the mean $p$, which is known as the \emph{relative error}. Here, for crude Monte Carlo the relative error is $\frac{\sqrt{Var(Z)}}{p}=\frac{\sqrt{p(1-p)}}{p}=\sqrt{\frac{1-p}{p}}$ which blows up as $p\to0$. Consequently, the sufficient size for $n$ also blows up as $p\to0$. In particular, if $p$ decays exponentially in $\gamma$ -- a typical scaling in large deviations, then the required $n$ to achieve \eqref{basic} also scales exponentially.

\subsection{Importance Sampling and Asymptotic Efficiency}
The above challenge motivates variance reduction techniques to reduce the sample size requirement. Among the most popular is IS. In this approach, we generate samples from an alternate measure $\tP$ where $\tP$ satisfies $ PI_{\cA_\gamma}\ll\tP$ (i.e., $PI_{\cA_\gamma}$ is absolutely continuous with respect to $\tP$), and use $Z=I_{\cA_{\gamma}}\frac{dP}{d\tP}$ as an unbiased output for $p$, where $\frac{dP}{d\tP}$ is the Radon-Nikodym derivative, or the so-called likelihood ratio, between $P$ and $\tP$. Though this output is always unbiased thanks to the likelihood ratio adjustment, the performance of the IS estimator in terms of variability heavily depends on the choice of the IS probability measure $\tP$. In the literature, several efficiency criteria for IS estimators have been developed. A common criterion is asymptotic efficiency \citep{asmussen2007stochastic,heidelberger1995fast,juneja2006rare}:

\begin{definition}[Asymptotic efficiency]
The IS estimator $Z=I_{\cA_{\gamma}}\frac{dP}{d\tP}$ under $\tP$ is said to achieve \textit{asymptotic efficiency} if $\lim_{\gamma\to\infty}\frac{\log(\tE(Z^2))}{\log p}=2$ where $\tE(\cdot)$ denotes the expectation under $\tP$.\label{def:AE}
\end{definition}

For functions $f,g:\R\to\R$, we say $g(\gamma)$ is subexponential in $f(\gamma)$ as $\gamma\to\infty$ if $\lim_{\gamma\to\infty}\frac{\log g(\gamma)}{f(\gamma)}=0$, i.e. $g(\gamma)=\exp(f(\gamma)o(1))$. For functions $f,g,h$, clearly if $\lim_{\gamma\to\infty}\frac{\log g(\gamma)}{f(\gamma)}=\lim_{\gamma\to\infty}\frac{\log h(\gamma)}{f(\gamma)}$, then $\lim_{\gamma\to\infty}\frac{\log(g(\gamma)/h(\gamma))}{f(\gamma)}=0$ and hence $g(\gamma)/h(\gamma)$ is subexponential in $f(\gamma)$ as $\gamma\to\infty$. By taking $g(\gamma)=\tilde E(Z^2)$, $h(\gamma)=p^2$ and $f(\gamma)=\log p$, we see that the condition in Definition~\ref{def:AE} is equivalent to the condition that $\tE(Z^2)/p^2$ (or $\tVar(Z)/p^2$) is subexponential in $-\log p$ as $\gamma\to\infty$.

From \eqref{Chebyshev} and its subsequent discussion, asymptotic efficiency implies that the required simulation size $n$ to attain a prefixed relative error grows only subexponentially in $-\log p$. As a stronger requirement, $Z$ is said to have a bounded relative error if $\limsup_{\gamma\to\infty}\frac{\tVar(Z)}{p^2}<\infty$, which implies that the required simulation size remains bounded no matter how small $p$ is. The criterion of bounded relative error is sometimes too strict to achieve, so we focus on asymptotic efficiency in this paper. More efficiency criteria could be found in \citet{juneja2006rare,ecuyer2010asymptotic,blanchet2012state}.

\subsection{Large Deviations and Dominating Points}
In the large deviations setting, the classical notion of \textit{dominating points} is used to guarantee asymptotic efficiency of IS \citep{sadowsky1990large}. To explain, let us first recall the so-called \emph{rate function} in the large deviations theory which, intuitively speaking, measures the likelihood of hitting each point on an exponential scale.
More specifically, we consider the standard Gartner-Ellis regime as follows \citep{dembo2009large,bucklew2004rare}. Without loss of generality consider $\gamma>0$. Suppose that  $\cA_{\gamma}=\{\frac{1}{\gamma}X_{\gamma}\in \cE\}$ where $\{X_{\gamma}\}_{\gamma}$ are $\R^d$-valued random variables and $\cE$ is a fixed Borel set in $\R^d$. We define $\mu_{\gamma}(x)=\frac{1}{\gamma}\log E(e^{x^\top  X_{\gamma}}),x\in\R^d$ as the scaled logarithmic moment generating function. We denote $\mathcal D(f)=\{x:f(x)<\infty\}$ as the domain of a function $f$.  With these, we assume the following:
\begin{assumption}
$\mu_{\gamma}(x)$ satisfies the following conditions:
\begin{enumerate}
    \item $\mu(x)=\lim_{\gamma\to\infty}\mu_{\gamma}(x)$ exists for any $x\in\R^d$, where we allow $\infty$ both as a limit value and as an element of the sequence $\{\mu_{\gamma}(x)\}$;
    \item $0\in \mathcal{D}(\mu)^{\circ}$;
    \item $\mu$ is essentially smooth, i.e., $\mathcal{D}(\mu)^{\circ}$ is non-empty, $\mu$ is differentiable everywhere in $\mathcal{D}(\mu)^{\circ}$ and $\mu$ is steep.
\end{enumerate}
\label{asm:mu_GE}
\end{assumption}
Then we define the rate function $I(y)=\sup_{x\in\R^d}\{x^\top y-\mu(x)\},y\in\R^d$ as the Legendre transform of $\mu$. For any set $\cE\subset\R^d$, we denote $I(\cE)=\inf_{y\in\cE}I(y)$. We make the following assumptions for the set $\cE$:
\begin{assumption}
$\cE\subset \R^d$ is a Borel set such that $\overline{\cE}=\overline{\cE^{\circ}}$, $\cE^\circ\cap\mathcal{D}(I)^\circ\ne\emptyset$ and $I(\cE)>0$.
\label{asm:E_GE}
\end{assumption}

Under these assumptions, we have the following result (e.g., adapted from \cite{dembo2009large} Theorem~\ref{thm:strong_probabilistic_efficiency}.3.6):
\begin{theorem}[Gartner-Ellis Theorem]
Suppose that Assumption \ref{asm:mu_GE} holds. For any Borel set $\cE\subset\R^d$, we have that 
$$
-I(\cE^\circ)\leq\liminf_{\gamma\to\infty}\frac{1}{\gamma}\log P\left(\frac{1}{\gamma}X_{\gamma}\in\cE\right)\leq\limsup_{\gamma\to\infty}\frac{1}{\gamma}\log P\left(\frac{1}{\gamma}X_{\gamma}\in\cE\right)\leq -I(\overline{\cE}).
$$
If additionally Assumption \ref{asm:E_GE} holds, then
$$
\lim_{\gamma\to\infty}\frac{1}{\gamma}\log P\left(\frac{1}{\gamma}X_{\gamma}\in\cE\right)=-I(\cE).
$$
\label{thm:GE}
\end{theorem}

Assumptions \ref{asm:mu_GE} and \ref{asm:E_GE} are standard light-tailed conditions on $\mathcal A_\gamma$, which guarantee the considered probability $P(\mathcal A_\gamma)$ to decay exponentially in $\gamma$ with decay rate $I(\cE)$. The rate function $I(y)$ can be viewed as a measurement on the likelihood of hitting $y$ in the exponential scale. The most likely point to hit among $\cE$ is hence given by the minimizer of $I(y)$ over $\cE$, resulting in the overall exponential decay rate $I(\cE)$. Note that, by Theorem \ref{thm:GE}, we have $-\log p=\Theta(\gamma)$ as $\gamma\to\infty$, and hence subexponential (or exponential) in $-\log p$ is equivalent to subexponential (or exponential) in $\gamma$. 

Now we present the concept of dominating points and sets:
\begin{definition}[Dominating Set]
Suppose that Assumptions \ref{asm:mu_GE} and \ref{asm:E_GE} hold. We call $A=\{a_1,\dots,a_r\}\subset\partial\cE$ a dominating set for $\cE$ if 
\begin{enumerate}
    \item For each $i$, $a_i\in\mathcal{D}(I)^\circ$ and there exists a unique $s_{a_i}\in\R^d$ such that $\nabla\mu(s_{a_i})=a_i$;
    \item $\cE\subset\bigcup_{i=1}^m\{x\in\R^d:s_{a_i}^\top (x-a_i)\geq 0\}$;
    \item For any $i$, $A\setminus\{a_i\}$ does not satisfy Condition 2.
\end{enumerate}
We call any point in $A$ a dominating point. For two dominating points $a$ and $a'$, we say $a$ is more significant than $a'$ if $I(a)<I(a')$.
\label{def:dominating_set_GE}
\end{definition}

\begin{figure}[htbp]
    \centering
    \includegraphics[width=0.4\textwidth]{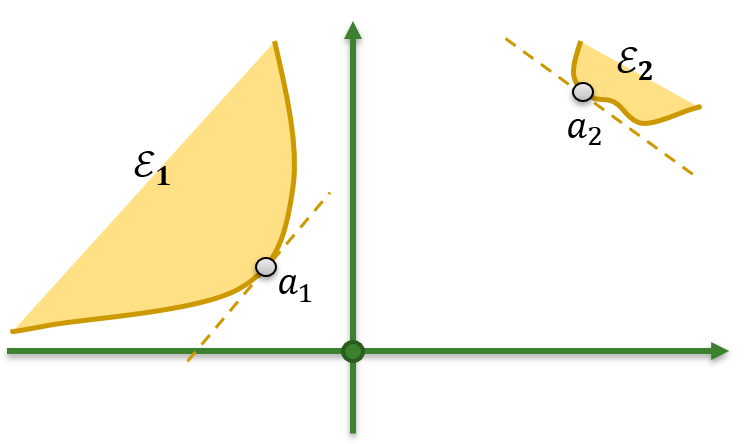}
    \caption{Illustration of rare-event set $\cE=\cE_1\cup\cE_2$ and dominating points.}
    \label{fig:illustration}
\end{figure}

Dominating points can be understood as the ``local minimizers'' of the rate function $I$ in the sense that $a_i$ is the minimizer of $I$ in $\cE\cap\{x\in\R^d:s_{a_i}^\top (x-a_i)\geq 0\}$. To understand this, first, Condition 1 in Definition \ref{def:dominating_set_GE} stipulates that $s_{a_i}$ is the gradient of $I(x)$ at the point $a_i$. Then, in Condition 2, $a_i$ can be seen as the minimizer of $I(x)$ over the set $\{x\in\R^d:s_{a_i}^\top (x-a_i)\geq 0\}$, a fact that follows from the first-order optimality condition in convex function minimization. Thus, Condition 2 stipulates that any point in the rare-event set $\cE$ has an $I$ value not less than one of the $a_i$'s. Geometrically, any points in the set $\cE$ must lie in the half-space tangentially cut by one of the $a_i$'s (i.e., the ``backyard" of the $a_i$). Figure \ref{fig:illustration} is an illustration of a rare-event set $\cE=\cE_1\cup\cE_2$ and the dominating set $\{a_1,a_2\}$. Here, $a_1$ is the global minimum rate point in $\cE$, but $\cE_2$ is not covered by $\{x:s_{a_1}^\top (x-a_1)\geq 0\}$, so $a_2$ is included in the dominating set. Finally, we note that Condition 3 in Definition \ref{def:dominating_set_GE} enforces the dominating set to be the minimal set of points such that the geometric properties in Conditions 1 and 2 are satisfied. 

Note that dominating points may not be local minimizers of the rate function $I$ in $\cE$ (even though they are minimizers in $\cE\cap\{x\in\R^d:s_{a_i}^\top (x-a_i)\geq 0\}$ as discussed above). Nonetheless, the most significant dominating points are indeed global minimizers of $I$ in $\cE$. This is presented in the following theorem:
\begin{theorem}
Suppose that Assumptions \ref{asm:mu_GE} and \ref{asm:E_GE} hold, and $A=\{a_1,\dots,a_r\}$ is a dominating set for $\cE$. Then $I(\cE)=\min_{i=1,\dots,r}I(a_i)$. That is to say, the most significant dominating points are global minimizers of $I$ in $\cE$.
\label{thm:global_minimizer}
\end{theorem}

We should also point out that dominating set defined according to Definition \ref{def:dominating_set_GE} may not be unique. Advantageously, the theory and estimators we present will flexibly apply to any such dominating set.

\subsection{Asymptotically (In)efficient Importance Samplers}
We are now ready to describe the main message of this section, which is the established recipe in constructing efficient IS. The standard proposal is to use a mixture of exponentially tilted distributions, where each exponential tilting is with respect to each dominating point. In particular, suppose that $A=\{a_1,\dots,a_r\}$ is a dominating set. Then the IS distribution is $\tP$ such that
\begin{equation}
    \frac{d\tP}{dP}=\sum_{i=1}^r\alpha_i e^{s_{a_i}^\top X_{\gamma}-\gamma\mu_{\gamma}(s_{a_i})}
    \label{eqn:full_IS_GE}
\end{equation}
with $\sum_{i=1}^r\alpha_i=1,\alpha_i>0,\forall i$. Here, $e^{s_{a_i}^\top X_{\gamma}-\gamma\mu_{\gamma}(s_{a_i})}dP$ is the exponential tilting towards the dominating point $a_i$ and $\alpha_i$'s are the mixing weights. The IS \eqref{eqn:full_IS_GE} is well known to be asymptotically efficient:
\begin{proposition}[Mixture IS is asymptotically efficient]
Suppose Assumptions \ref{asm:mu_GE} and \ref{asm:E_GE} hold, and the dominating set has finite cardinality, i.e., $r<\infty$. Then the IS distribution \eqref{eqn:full_IS_GE} with any fixed $\alpha_i$'s is asymptotically efficient.\label{classical AE}
\end{proposition}

While the proof of Proposition \ref{classical AE} is standard, we include it in the Appendix for self-containedness. Here, we describe the key intuition in justifying the necessity of mixture. First, the likelihood ratio in the considered mixture IS is
\begin{equation*}
    L=\frac{dP}{d\tP}=\frac{1}{\sum_{i=1}^r\alpha_i e^{s_{a_i}^\top X_{\gamma}-\gamma\mu_{\gamma}(s_{a_i})}}
\end{equation*}
and it satisfies that for any $i$,
\begin{equation}
    L\leq \frac{1}{\alpha_i e^{s_{a_i}^\top X_{\gamma}-\gamma\mu_{\gamma}(s_{a_i})}}=\frac{1}{\alpha_i}e^{-s_{a_i}^\top (X_{\gamma}-\gamma a_i)-\gamma (s_{a_i}^\top a_i-\mu_{\gamma}(s_{a_i}))}.\label{interim classical}
\end{equation}
In the exponent in the rightmost expression of \eqref{interim classical}, the second term $\gamma (s_{a_i}^\top a_i-\mu_{\gamma}(s_{a_i}))$ is approximately $\gamma I(a_i)$, and the first term is the ``overshoot" of the sampled $X_\gamma$ compared to the dominating point $a_i$. That is, if $X_\gamma$ is in the ``backyard" of $a_i$, then this term $s_{a_i}^\top (X_{\gamma}-\gamma a_i)\geq0$. The definition of dominating set, especially Condition 2 in Definition \ref{def:dominating_set_GE}, guarantees any $(1/\gamma)X_\gamma$ in $\cE$ must have $s_{a_i}^\top (X_{\gamma}-\gamma a_i)\geq0$ for at least one of the $i$'s. Thus, by decomposing the second moment of $Z=I_{\mathcal A_\gamma}L$ according to the backyards of $a_i$'s, we can ensure that the magnitude of the likelihood ratio, when $X_\gamma$ lies inside the rare event set, is properly controlled. More precisely, write $\cE=\bigcup_{i=1}^r\cE_i$ where each $\cE_i\subset\{x:s_{a_i}^\top(x-a_i)\geq0\}$. Then the second moment of the IS satisfies
$$\tilde E\left[I\left(\frac{1}{\gamma} X_\gamma\in\cE\right)L^2\right]\leq\sum_{i=1}^r\tilde E\left[I\left(\frac{1}{\gamma} X_\gamma\in\cE_i\right)L^2\right]\leq\sum_{i=1}^r\tilde E\left[I\left(\frac{1}{\gamma} X_\gamma\in\cE_i\right)\frac{1}{\alpha_i^2}e^{-2\gamma (s_{a_i}^\top a_i-\mu_{\gamma}(s_{a_i}))}\right]$$
which is approximately bounded by $\sum_{i=1}^re^{-2\gamma I(a_i)}/\alpha_i^2$ and hence $e^{-2\gamma I(\cE)}$ in the exponential scale, thus verifying asymptotic efficiency.


On the other hand, if we miss some dominating points in the construction of the mixture IS, then asymptotic efficiency may fail to be attained. Below we give a simple example to demonstrate this. 

\begin{proposition}[Missed dominating point leads to violation of asymptotic efficiency]
Suppose that we want to estimate $p=P(\frac{1}{\gamma}X_{\gamma}\in(-\infty,-2]\cup[1,\infty))$ where $X_{\gamma}\sim N(0,\gamma)$ under $P$. If the IS distribution is chosen as $X_{\gamma}\sim N(\gamma,\gamma)$, then $\tE(Z^2)/p^2=\Theta(\sqrt{\gamma} e^{3\gamma/2})$ grows exponentially in $-\log p=\Theta(\gamma)$, and hence $Z$ is not asymptotically efficient by definition.
\label{lem:example}
\end{proposition}
In this example, the dominating points are $1$ and $-2$, and $N(\gamma,\gamma)$ is the exponential tilt towards the first dominating point $1$ (for Gaussian distribution, exponential tilting amounts to a mean shift). Here, by considering only this point, it is possible that a generated $X_\gamma$ satisfies $(1/\gamma)X_\gamma\in(-\infty,-2]$ while the overshoot $s_{1}^\top (X_{\gamma}-\gamma 1)$, as explained for \eqref{interim classical}, takes a very negative value. This scenario contributes significantly to the overall variance and ultimately violates asymptotic efficiency.



Our main insight in this paper is a rebuke of the above viewpoint. More specifically, we argue that missing inferior dominating point, such as the example in Proposition \ref{lem:example}, can still result in a good IS according to our beginning criterion \eqref{basic}. A core ingredient of this assertion is to question the use of asymptotic efficiency, or more generally variance-based efficiency criteria. Before delving into the theory, let us first present some numerical results to shed light on how much difference it makes to use different numbers of dominating points in the IS mixture. This is the focus of our next section. 

\section{Motivating Experimental Results}\label{sec:prelim numerics}
We run three numerical examples to demonstrate that missing dominating points in IS construction, while provably leads to asymptotic inefficiency, could perform well empirically. This thus suggests an inadequacy in using asymptotic efficiency, or more generally variance-based criteria, to measure the performances of rare-event estimators. Besides, the computationally demanding example in Section \ref{sec:numerical_mnist} justifies the motivation why we seek to reduce the number of used dominating points in the IS mixture.
\label{sec:numerical}


\subsection{Large Deviations of an I.I.D. Sum}
\label{sec:numerical_glasserman}

We consider the problem of estimating the tail probability involving a sum of random variables, where $Y_1,Y_2,...$ are i.i.d and we are interested in $$P(|S_m|\geq a m ), $$
where $S_m=\sum_{i=1}^m Y_i$. We consider $m$ as the rarity parameter $\gamma$ presented in Section \ref{sec:background}. Using the notation in the Gartner-Ellis regime, we have $X_{\gamma}=S_{m}$ and $\cE=(-\infty,-a]\cup[a,\infty)$. Then $\mu(x)=\log E(e^{x Y_1})$ and we suppose Assumption~\ref{asm:mu_GE} is satisfied. By Theorem~\ref{thm:GE}, when $|EY_1|<a$, if $s_a$ and $s_{-a}$ satisfy $\nabla\mu(s_a)=a$ and  $\nabla\mu(s_{-a})=-a$, then  we have 
\begin{equation} \label{eq:counter_rate_1}
    -\lim_{m\to \infty} \frac{1}{m}  \log P(S_m \geq am)=s_a a -\mu (s_a) = I(a)
\end{equation}
and \begin{equation} \label{eq:counter_rate_2}
     -\lim_{m\to \infty} \frac{1}{m}  \log P(S_m \leq -am)=-s_{-a} a -\mu (s_{-a}) = I(-a) . 
\end{equation}

For this problem, \citet{glasserman1997counterexamples} Section 3 provides two estimators, $\hat \alpha(m)$ and $\hat \beta(m)$. Specifically, we have $\hat \alpha(m)=\exp (-s_a S_m +m \mu(s_a)) I_{\{ |S_m| \geq am \}}$ with samples of $Y_i$ generated from exponentially tilted distribution using $s_a$. The estimator $\hat \beta(m)=\exp (-s_a S_m +m \mu(s_a))I_{\{ S_m \geq am \}} +  \exp (-s_{-a} S_m' +m \mu(s_{-a})) I_{\{ S_m' \leq -am \}}$ with $S_m$ and $S_m'$ constructed from independent sequences of i.i.d. $Y_i$'s and $Y_i'$'s generated from exponentially tilted distributions using $s_a$ and and $s_{-a}$ respectively. That is, $\hat\beta(m)$ attempts to estimate $P( S_m \geq am )$ and $P( S_m \leq -am )$ separately using different IS samples and sum up these estimates. Here, $\hat{\alpha}(m)$ only uses one dominating point whereas $\hat{\beta}(m)$ uses both points (note that even though $\hat\beta(m)$ does not use the mixture IS scheme in \eqref{eqn:full_IS_GE}, the idea is similar in that it accounts for both dominating points). In our experiment, we follow \citet{glasserman1997counterexamples} to set $Y_1=A-B$ with $A \sim N(1.5,1)$, $B\sim Exp(1)$ and $A,B$ independent, and $a=1.5$ (in this case, $s_a =(\sqrt{5}-1)/2 $, $s_{-a}=-2+\sqrt{2}$, and $I(a) \approx 0.2902$, $I(-a) \approx 0.7044$).


We run numerical experiments with $m=10$, $30$, $50$ and $100$. The results using $10^4$ samples are shown in Table \ref{table:glasserman}. By comparing the numbers in the second and third rows, we observe that $\hat \alpha(m)$ and $\hat \beta(m)$ have very similar empirical performances. However, note that:
\begin{proposition}
Under the problem specification above, $\hat \beta(m)$ is asymptotically efficient while $\hat\alpha(m)$ is not. In fact, $\tilde E(\hat \alpha^2(m))\to \infty$ as $m\to\infty$ where $\tilde E$ denotes the expectation under the exponential tilting towards $a$.\label{prop:example1}
\end{proposition}

In view of Proposition \ref{prop:example1}, $\hat \alpha(m)$ is arguably a very poor estimator as it bears an exploding variance. We therefore see an apparent discrepancy between empirical performances and theoretical guidance -- The theoretically bad variance does not result in poor empirical performances. Proposition \ref{prop:example1} is proved in \citet{glasserman1997counterexamples}, where the asymptotic efficiency of $\hat \beta(m)$ follows from their Proposition 1, while the variance behavior of $\hat \alpha(m)$ appears in their Theorem 1. 

\begin{table}[htbp]
\caption{Point estimates (and 95\% CI) using IS estimators for the tail probability with different $m$.}
\centering
\resizebox{\columnwidth}{!}{%
\begin{tabular}{|l|l|l|l|l|}
\hline
$m$  & 10                & 30                & 50                 & 100                \\ \hline
$\hat \alpha(m)$ & 8.22($\pm$0.26) $\times10^{-3}$  & 1.60($\pm$0.07) $\times10^{-5}$ & 3.77($\pm$0.18) $\times10^{-8}$ & 1.34($\pm$0.08) $\times10^{-14}$\\ \hline
$\hat \beta(m)$ & 8.29($\pm$0.26) $\times10^{-3}$& 1.60($\pm$0.07) $\times10^{-5}$ & 3.77($\pm$0.18) $\times10^{-8}$ & 1.34($\pm$0.08) $\times10^{-14}$ \\ \hline
\end{tabular}} \label{table:glasserman}
\end{table}

\subsection{Overshoot Probability of Random Walk}
\label{sec:numerical_sum}

We consider the problem of estimating the overshoot probability of the finite-horizon maximum of a random walk. We define the probability of interest as 
\begin{equation*}
p=P\left(\max_{m=1,...,d} S_m \geq a\right),
\end{equation*}
where $S_m=\sum_{i=1}^m Y_i$ and $Y_i$'s are Gaussian distributed with mean 0, standard deviation $\sigma$, and pairwise correlation $-0.02$, i.e., $corr(Y_i,Y_j)=-0.02$ for any $i,j\in \{1,...,d\}$ with $i\neq j$.  Suppose that the rarity parameter is $\gamma=1/\sigma^2\to\infty$. We note that we can reformulate this target rare event as $\left\{\frac{1}{\gamma}X_{\gamma}\in\bigcup_{m=1}^d \cH_m\right\}$ where $X_{\gamma}=\gamma(Y_1,\dots,Y_d)^\top $ and $\cH_m = \{x \in \mathbb{R}^d : \sum_{i=1}^m x_i \geq a\}  $ with $x_i$ denoting the $i$th element in $x$. This decomposition allows us to construct an IS estimator using the $d$ dominating points corresponding to each half-space $\cH_m$. More specifically, in this example, $I(y)=\frac{1}{2}y^\top  \Sigma^{-1} y$, and the dominating points, ranking from the most to the least significant (i.e., increasing rate function value), are $a_1=\frac{a\Sigma e_d}{e_d^\top \Sigma e_d},a_2=\frac{a\Sigma e_{d-1}}{e_{d-1}^\top \Sigma e_{d-1}},\dots,a_d=\frac{a\Sigma e_1}{e_1^\top \Sigma e_1}$, where $e_i$ denotes the vector with 1 in the first $i$ elements  and 0 for the rest. 

\begin{figure}[htbp]
      \centering
     \subfloat[$\sigma=0.3$. \label{fig:overshoot_03_2}
     ]{\includegraphics[width=0.45\textwidth]{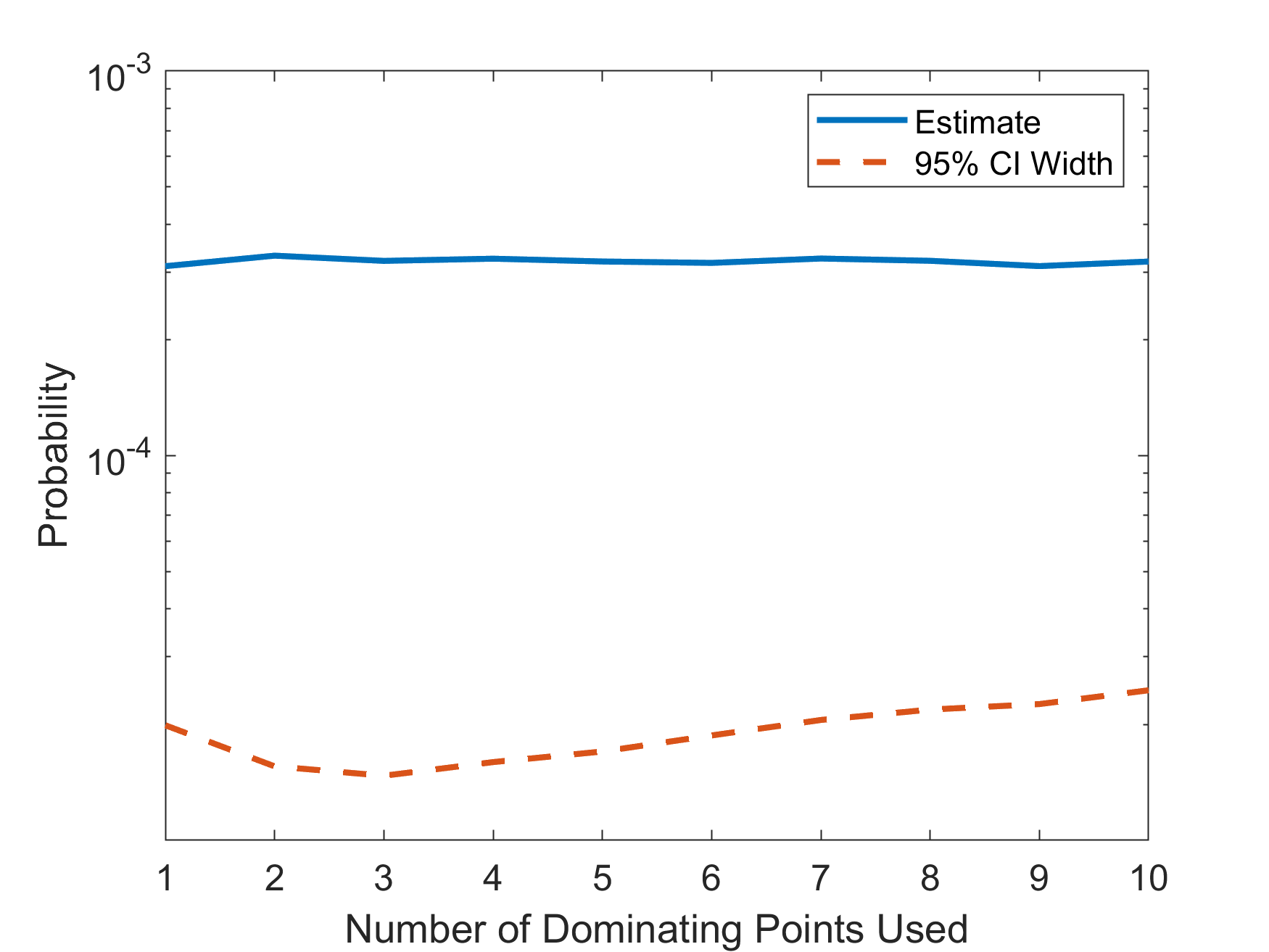}} \quad 
     \subfloat[ $\sigma=0.2$.\label{fig:overshoot_02_2}]{\includegraphics[width=0.45\textwidth]{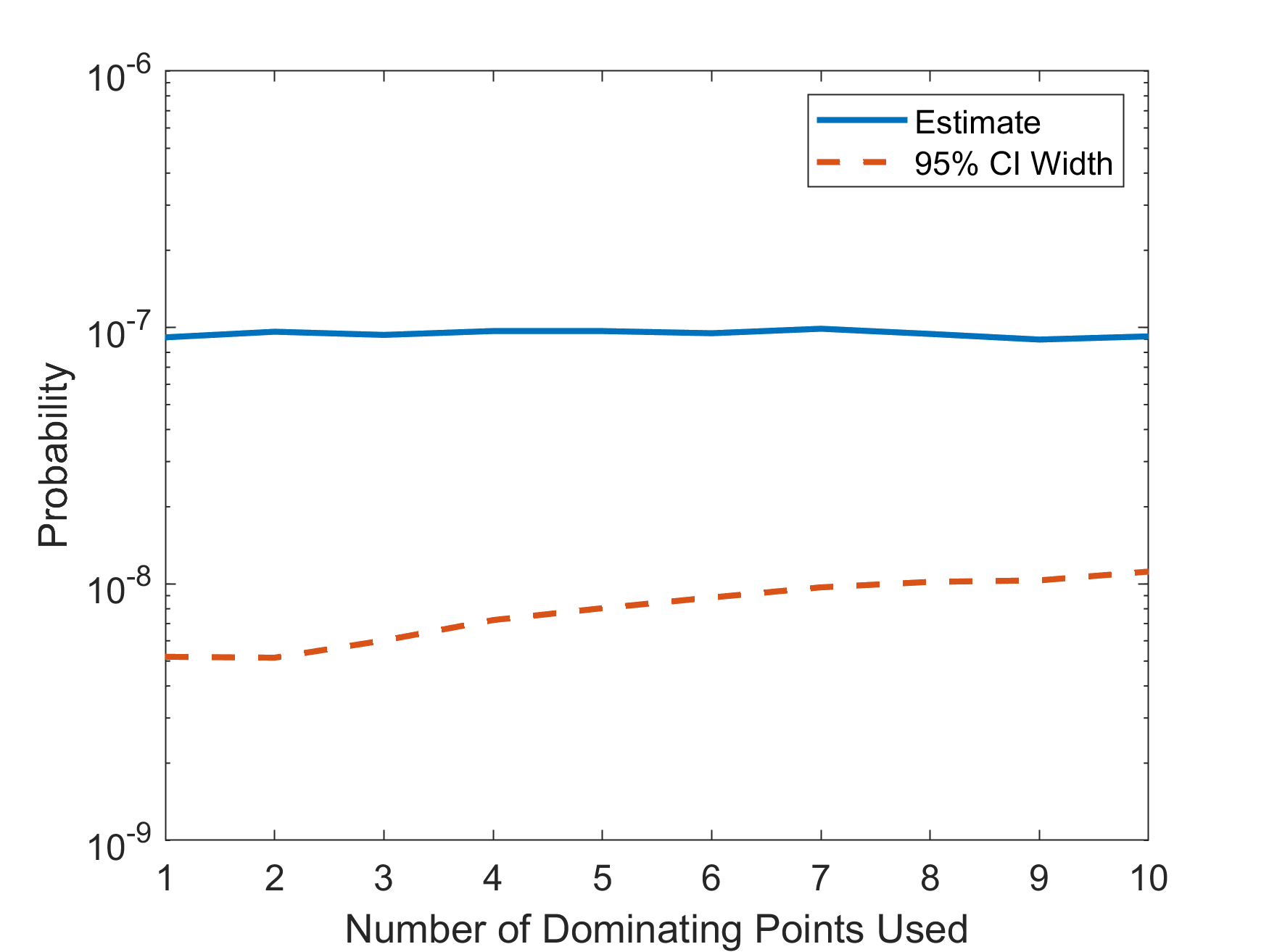}}
     \caption{Simulation results for the random walk experiment. Point estimates and CI widths for IS estimators using different numbers of dominating points.}
     \label{fig:overshoot_IS2}
 \end{figure}

In our experiments, we fix $a$ and vary $\sigma$ for different rarity levels. In addition, we set $d=10$. We generate $10^4$ samples from IS distributions using a varying, partial list of dominating points. That is, we choose the IS distributions for  $X_{\gamma}$ as $\frac{1}{k}\sum_{i=1}^k\phi(x;\gamma a_i,\gamma \Sigma)$, where $\phi(x;a,b)$ denotes the Gaussian density with mean $a$ and variance $b$, for $k=1,\dots,d$. The performances of these IS estimators are shown in Table~\ref{table:overshoot2}. We observe that these estimators with different numbers of dominating points all perform similarly.  In particular, we present the cases with $\sigma=0.2$ and $\sigma=0.3$ respectively in Figure~\ref{fig:overshoot_IS2}. We observe that in both cases the performances of the IS estimators are almost independent of the number of used dominating points, with the probability estimates all comparable while using more dominating points slightly increases the CI width.

On the other hand, Proposition~\ref{classical AE} implies that the IS using all dominating points is asymptotically efficient while we have:
\begin{proposition}
Under the problem specification above, the IS estimator that exponentially tilts towards the most significant dominating point $a_1$, i.e., $X_\gamma$ distributed as $N(\gamma a_1,\gamma \Sigma)$, is not asymptotically efficient.
\label{prop:example2}
\end{proposition}

Thus, like in Section \ref{sec:numerical_glasserman}, there appears a mismatch between theoretical guidance and empirical observation. The asymptotic inefficiency of simple exponential tilting towards only the most significant dominating point does not result in a poor experimental performance.

\begin{table}[htbp] 
\caption{Point estimates (and  95\% CIs) from IS estimators using different numbers of dominating points for the overshoot probability. \# denotes the number of dominating points used in the IS estimator. }
\centering
\resizebox{\columnwidth}{!}{%
\begin{tabular}{|l|l|l|l|l|l|l|}
\hline
 $\sigma$  & 0.2                & 0.22                & 0.24                 & 0.26                & 0.28                & 0.3                 \\ \hline
 \#  & prob (with CI)        & prob    (with   CI)        & prob    (with   CI)        & prob  (with    CI)        & prob   (with   CI)        & prob  (with    CI)        \\ \hline
1  & 9.15($\pm$0.52) $\times 10^{-8}$ & 1.24($\pm$0.19) $\times 10^{-6}$ & 7.96($\pm$0.85) $\times 10^{-6}$ & 3.52($\pm$0.29) $\times 10^{-5}$ & 1.15($\pm$0.08) $\times 10^{-4}$ & 3.11($\pm$0.20) $\times 10^{-4}$ \\ \hline
2  & 9.63($\pm$0.52) $\times 10^{-8}$ & 1.24($\pm$0.08) $\times 10^{-6}$ & 8.38($\pm$0.52) $\times 10^{-6}$ & 3.71($\pm$0.20) $\times 10^{-5}$ & 1.23($\pm$0.06) $\times 10^{-4}$ & 3.31($\pm$0.16) $\times 10^{-4}$ \\ \hline
3  & 9.36($\pm$0.60) $\times 10^{-8}$ & 1.15($\pm$0.07) $\times 10^{-6}$ & 7.87($\pm$0.43) $\times 10^{-6}$ & 3.61($\pm$0.19) $\times 10^{-5}$ & 1.21($\pm$0.06) $\times 10^{-4}$ & 3.21($\pm$0.15) $\times 10^{-4}$ \\ \hline
4  & 9.69($\pm$0.72) $\times 10^{-8}$ & 1.18($\pm$0.08) $\times 10^{-6}$ & 7.98($\pm$0.49) $\times 10^{-6}$ & 3.64($\pm$0.20) $\times 10^{-5}$ & 1.21($\pm$0.06) $\times 10^{-4}$ & 3.25($\pm$0.16) $\times 10^{-4}$ \\ \hline
5  & 9.68($\pm$0.80) $\times 10^{-8}$ & 1.17($\pm$0.09) $\times 10^{-6}$ & 7.92($\pm$0.54) $\times 10^{-6}$ & 3.59($\pm$0.22) $\times 10^{-5}$ & 1.20($\pm$0.07) $\times 10^{-4}$ & 3.20($\pm$0.17) $\times 10^{-4}$ \\ \hline
6  & 9.50($\pm$0.89) $\times 10^{-8}$ & 1.15($\pm$0.10) $\times 10^{-6}$ & 7.79($\pm$0.60) $\times 10^{-6}$ & 3.55($\pm$0.25) $\times 10^{-5}$ & 1.19($\pm$0.08) $\times 10^{-4}$ & 3.17($\pm$0.19) $\times 10^{-4}$ \\ \hline
7  & 9.89($\pm$0.97) $\times 10^{-8}$ & 1.20($\pm$0.11) $\times 10^{-6}$ & 8.13($\pm$0.66) $\times 10^{-6}$ & 3.69($\pm$0.27) $\times 10^{-5}$ & 1.22($\pm$0.08) $\times 10^{-4}$ & 3.26($\pm$0.21) $\times 10^{-4}$ \\ \hline
8  & 9.44($\pm$1.02) $\times 10^{-8}$ & 1.16($\pm$0.11) $\times 10^{-6}$ & 7.93($\pm$0.70) $\times 10^{-6}$ & 3.62($\pm$0.29) $\times 10^{-5}$ & 1.20($\pm$0.09) $\times 10^{-4}$ & 3.21($\pm$0.22) $\times 10^{-4}$ \\ \hline
9  & 8.97($\pm$1.03) $\times 10^{-8}$ & 1.11($\pm$0.12) $\times 10^{-6}$ & 7.63($\pm$0.72) $\times 10^{-6}$ & 3.48($\pm$0.30) $\times 10^{-5}$ & 1.16($\pm$0.09) $\times 10^{-4}$ & 3.11($\pm$0.23) $\times 10^{-4}$ \\ \hline
10 & 9.23($\pm$1.12) $\times 10^{-8}$ & 1.15($\pm$0.13) $\times 10^{-6}$ & 7.87($\pm$0.77) $\times 10^{-6}$ & 3.55($\pm$0.32) $\times 10^{-5}$ & 1.19($\pm$0.10) $\times 10^{-4}$ & 3.20($\pm$0.25) $\times 10^{-4}$ \\ \hline
\end{tabular}}\label{table:overshoot2}
\end{table}

\subsection{Robustness Assessment for an MNIST Classification Model}
\label{sec:numerical_mnist}
We consider a rare-event probability estimation problem from an image classification task. Our goal is to estimate the probability of misclassification when the input of a prediction model is perturbed by tiny noise. This probability estimate is of interest as a robustness measure of the prediction model~\citep{webb2018statistical}. 
More specifically, suppose that the prediction model $g$ is able to predict the label of input $x_0$, i.e. $g(x_0)=c$ where $c$ is the true label of $x_0$. Then $P(g(x_0+\varepsilon)\ne c)$ where $\varepsilon$ is a random perturbation can be used to measure the robustness of $g$.

In particular, we consider the classification problem on MNIST dataset which contains 70,000 images of handwritten digits and each image consists of $28\times 28$ pixels. We train a 2-ReLU-layer neural network with 20 neurons in each layer using 60,000 training data, which achieves approximately 95\% of testing data accuracy in predicting the digits. We perturb a fixed input (that is correctly predicted) with a Gaussian noise with mean 0 and standard deviation $\sigma$ on each of the 784 dimensions to assess the robustness of the prediction. Note that the rarity of this problem is determined by the value of $\sigma$, and we let the rarity parameter $\gamma=1/\sigma^2\to\infty$. The target rare event can be reformulated as $\left\{\frac{1}{\gamma}X_{\gamma}\in\{x:g(x_0+x)\neq c)\}\right\}$ where $X_{\gamma}=\gamma\varepsilon$.

We apply mixtures of exponential tiltings as IS estimators for this problem, namely by considering the IS distribution $\frac{1}{k}\sum_{i=1}^k\phi(x;\gamma a_i,\gamma I)$, where $\phi(x;a,b)$ denotes the Gaussian density with mean $a$ and variance $b$ as in Section \ref{sec:numerical_sum}, for $k=1,2,\ldots$. Here $a_i,i=1,2,\ldots$ denote the dominating points. In order to compute these points, we apply the scheme introduced in~\citet{huang2018designing} and~\citet{bai2022rare}, which sequentially searches for dominating points by minimizing the rate function on the rare-event set that excludes the half-spaces cut from previous more significant dominating points. In the Gaussian case with piecewise linear rare-event set boundary as in our current example, each iteration amounts to finding the highest-density point on a piecewise-linear-boundary set, which can be conducted using mixed integer programming (see Algorithm \ref{algo:dominating_points} in the Appendix).
Due to the high dimensionality of the input space and the complexity of the neural network predictor, the number of dominating points in this problem is huge. We implemented this sequential searching algorithm and it took a week to find the first 100 dominating points. Since we stopped the algorithm prematurely, the actual number of dominating points can be much larger. We run IS distributions with different numbers of dominating points (ranging from 1 to 41) and magnitudes of $\sigma$ (ranging from $0.1$ to $0.2$) and report the estimated probabilities and CIs. We use $10^5$ samples for IS estimators and $10^7$ samples for crude Monte Carlo estimators. 
\begin{figure}[htbp]
      \centering
     \subfloat[ \label{fig:MNIST_one}
     ]{\includegraphics[width=0.45\textwidth]{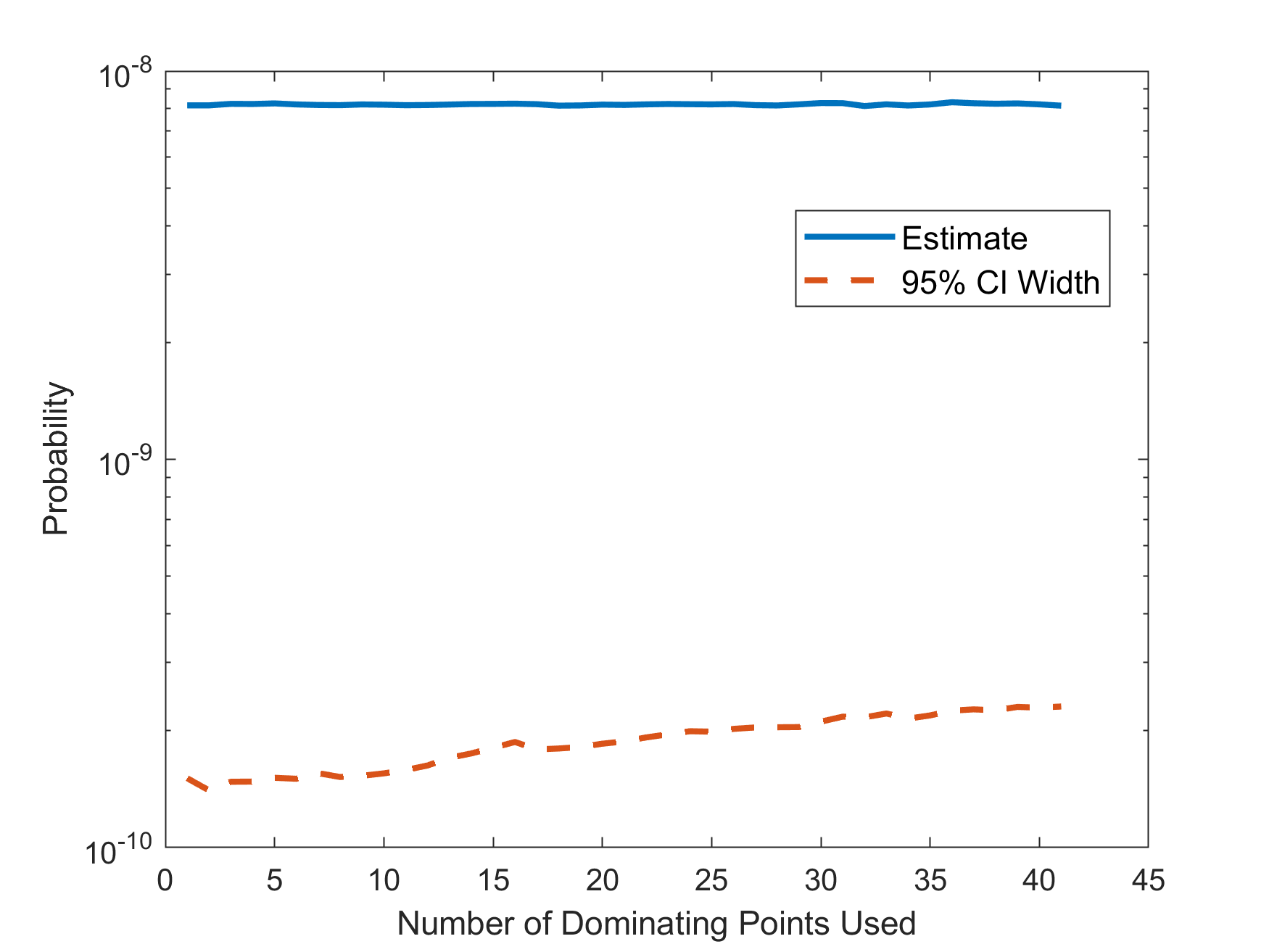}} \quad 
     \subfloat[ \label{fig:MNIST_compare}]{\includegraphics[width=0.45\textwidth]{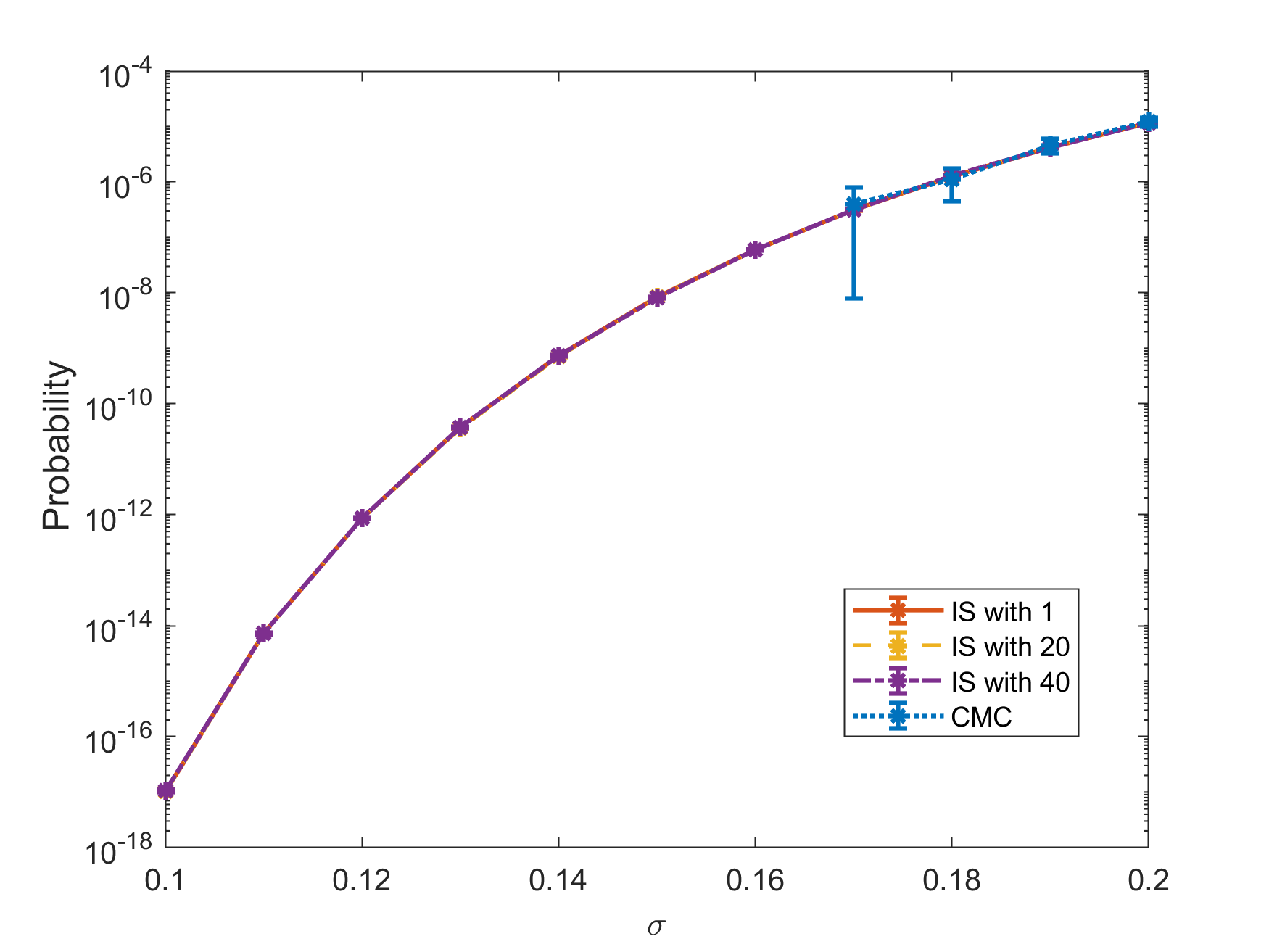}}
     \caption{Simulation results for the MNIST experiment. (a) Point estimates and CI widths from IS estimators using different numbers of dominating points. (b) Point estimates from IS estimators using different numbers of dominating points (IS with 1, 20, and 40) and crude Monte Carlo (CMC), with vertical error bars representing their 95\% CIs (the CIs for the IS estimates are extremely narrow).}
     \label{fig:MNIST}
 \end{figure}


Figures~\ref{fig:MNIST_one} and \ref{fig:MNIST_compare} show the results. Missing less significant dominating points does not seem to make noticeable differences in this problem. As shown in Figure~\ref{fig:MNIST_one}, when we fix the rarity of the problem, the estimate is not sensitive to the number of dominating points. The CI width has an increasing trend as the number of dominating points gets larger, indicating that additional dominating points can in fact even hurt performances. 

In Figure~\ref{fig:MNIST_compare}, we vary the rarity of the problem and compare the performances of different IS estimators and crude Monte Carlo. Note that estimates using crude Monte Carlo are unavailable for rarer configurations due to its inefficiency. We observe that the estimates from different IS estimators overlap visually in all considered cases, which indicates that the differences among these estimates are negligible. We also note that these estimates are consistent with the crude Monte Carlo estimates (when available), which shows their correctness. 

Nonetheless, once again we have an apparent mismatch between theoretical inefficiency and good empirical performances:
\begin{proposition}
Under the problem specification above, the IS estimator that exponentially tilts towards the most significant dominating point, i.e., $X_\gamma$ distributed as $N(\gamma a_1,\gamma I)$, is not asymptotically efficient. 
\label{prop:example3}
\end{proposition}

\section{Probabilistic Efficiency}
\label{sec:concept}
Section \ref{sec:prelim numerics} shows that IS estimators that miss some dominating points could perform competitively compared to estimators that consider all of them, thus suggesting a gap between the notion of asymptotic efficiency and empirical performances. In light of this, we propose the concept of \emph{probabilistic efficiency} as a relaxation of asymptotic efficiency. The key of probabilistic efficiency is to consider the high-probability relative discrepancy of the estimator from the ground truth directly, instead of using the relative error or equivalently the estimation variance. The latter, as can be seen in the arguments in Section \ref{sec:background}, provides a sufficient, but not necessary, condition on the required sample size. In particular, there is an intrinsic looseness brought by the Markov or Chebyshev inequality \eqref{Chebyshev} that converts relative error into the required sample size. 

To proceed, we first define the following:
\begin{definition}[Minimal relative discrepancy]
For any estimator $\hp$ of $p$ and any $\varepsilon>0$, the \emph{minimal relative discrepancy} of $\hp$, at tolerance level $\varepsilon$, is given by
\begin{equation}
    \delta_{\varepsilon}(\hp,p):=\inf\{\delta>0:\tP(|\hp-p|>\delta p)\leq \varepsilon\}.
    \label{eqn:delta}
\end{equation}\end{definition}
The minimal relative discrepancy $\delta_{\varepsilon}(\hp,p)$ measures the relative accuracy of the estimator $\hp$, in that it gives the smallest relative discrepancy of $\hp$ from $p$ that can be achieved with probability $1-\varepsilon$. Thus the smaller is $\delta_{\varepsilon}(\hp,p)$, the more accurate is $\hp$. Note that in \eqref{eqn:delta}, the probability $\tilde P$ is the one generating the estimator $\hat p$.

We say that $\hp$ is probabilistically efficient if $\delta_{\varepsilon}(\hp,p)$ can be made small in some sense, without needing to use a gigantic amount of computation. More precisely, we propose the following notions:
\begin{definition}[Probabilistic Efficiency]
Suppose that $\{\cA_{\gamma}\}_{\gamma}$ is an indexed family of rare events and $p=P(\cA_{\gamma})\to 0$ as $\gamma\to\infty$. Consider an estimator $\hp$ obtained from $n=n(\gamma)$ independent replications of $Z$. For any $\varepsilon>0$, we define $\delta_{\varepsilon}(\hp,p)$ as in \eqref{eqn:delta}. Then 
\begin{enumerate}
    \item We call $Z$ \emph{strongly probabilistically efficient} if we can choose $n$ subexponential in $-\log p$ such that, for any $\varepsilon>0$,  $\lim_{\gamma\to\infty}\delta_{\varepsilon}(\hp,p)=0$;    
    \item We call $Z$ \emph{weakly probabilistically efficient} if we can choose $n$ subexponential in $-\log p$ such that, for any $\varepsilon>0$, $\limsup_{\gamma\to\infty}\delta_{\varepsilon}(\hp,p)<1$.
    
\end{enumerate}
\end{definition}

Note that strong probabilistic efficiency matches the usual notion in statistical estimation. That is, the estimator approaches the target parameter as $\gamma\to\infty$. In contrast, weak probabilistic efficiency only cares about a correct magnitude. While this may appear less desirable, in rare-event estimation a correct magnitude can be viewed as sufficient as the target quantity is very small, and this weaker notion allows more flexibility in constructing estimators. We also contrast our proposed probabilistic efficiency with a notion named probabilistic bounded relative error proposed in \citet{tuffin2012probabilistic}, where the IS measure is randomly chosen and efficiency is achieved if the resulting random relative error of the IS estimator is bounded by some constant with high probability, which is conceptually different from our notion.

The following shows that probabilistic efficiency is a relaxation of asymptotic efficiency:
\begin{proposition}
If $Z$ is asymptotically efficient, then $Z$ is strongly probabilistically efficient.
\label{prop:ae_imply_pe}
\end{proposition}
\proof{Proof of Proposition \ref{prop:ae_imply_pe}.}
For any unbiased estimator $Z$, we have that for any $\varepsilon>0$,
\begin{equation*}
    \tP\left(|\hp-p|>\sqrt{\frac{\tVar(Z)}{\varepsilon np^2}}p\right)\leq \varepsilon
\end{equation*}
and hence by definition,
\begin{equation*}
    \delta_{\varepsilon}(\hp,p)\leq \sqrt{\frac{\tVar(Z)}{\varepsilon np^2}}.
\end{equation*}
If $Z$ is asymptotically efficient, $\frac{\tVar(Z)}{p^2}$ grows at most subexponentially in $-\log p$, so we could choose $n$ subexponentially growing in $-\log p$ such that $\lim_{\gamma\to\infty}\delta_{\varepsilon}(\hp,p)=0$ for any $\varepsilon>0$. By definition, $Z$ is strongly probabilistically efficient.\hfill\Halmos
\endproof

While asymptotic efficiency implies strong probabilistic efficiency, we note that these two notions are not equivalent. In Section \ref{sec:prelim numerics}, Propositions \ref{prop:example1}--\ref{prop:example3} show that asymptotic efficiency does not hold for the considered IS estimators in all the presented examples, but Theorem \ref{PE examples} in Section \ref{sec:verification} will show that strong probabilistic efficiency actually holds for all of them.

Now we explain how probabilistic efficiency helps us understand the influence of missing some dominating points. Recall the example where $\cA_{\gamma}=\{\frac{1}{\gamma}X_{\gamma}\in \cE\}$ and $\cE$ comprises two disjoint and faraway pieces $\cE_1$ and $\cE_2$. The dominating points are respectively $a_1$ and $a_2$ (recall Figure \ref{fig:illustration}). Denote $p_j=P(\frac{1}{\gamma}X_{\gamma}\in\cE_j),j=1,2$, and we assume that $p_2$ is exponentially smaller than $p_1$. If we focus on $\cE_1$ and simply use the exponential tilting towards $a_1$ as the IS distribution, then we face the risk of having a sample falling into $\cE_2$ while the associated likelihood ratio is very high, which leads to asymptotic inefficiency. However, experimentally if we run the simulation with a moderate sample size, then most likely none of the samples fall into $\cE_2$. Conditional on not hitting $\cE_2$, we actually get an estimate close to $p_1$, which is in turn close to $p$. In other words, even if the resulting IS estimator is not asymptotically efficient, it could still give a good estimate in terms of its distance to $p$, as long as the sample size is not overly big. The latter is precisely the paradigm of probabilistic efficiency.

More concretely, we have the following theorem:
\begin{theorem}[Achieving strong probabilistic efficiency]
Suppose that $\{\cA_{\gamma}\}_{\gamma}$ is an indexed family of rare events and $p=P(\cA_{\gamma})\to 0$ as $\gamma\to\infty$. We write $\cA_{\gamma}=\cA_{\gamma}^1\cup \cA_{\gamma}^2$ where $\cA_{\gamma}^1$ and $\cA_{\gamma}^2$ are two disjoint events. Denote $p_j=P(\cA_{\gamma}^j),j=1,2$. Assume that 
\begin{enumerate}
    \item $\frac{p_1}{p}\to1$ as $\gamma\to\infty$;
    \item We have an asymptotically efficient IS estimator for $p_1$ obtained from $Z_1=I_{\cA_{\gamma}^1}\frac{dP}{d\tP}$ under $\tP$. This implies that there exists $n=n(\gamma)$ growing subexponentially in $-\log p$ such that $\frac{\tVar(Z_1)}{np_1^2}\to 0$ as $\gamma\to\infty$;
    \item $\tp_2:=\tP(\cA_{\gamma}^2)$ satisfies that $n\tp_2\to0$ as $\gamma\to\infty$.
\end{enumerate}
Let $\hp$ be the sample mean of $n$ independent replications of $Z=I_{\cA_{\gamma}}\frac{dP}{d\tP}$ under $\tP$. For any $\varepsilon>0$, define $\delta_{\varepsilon}(\hp,p)$ as in \eqref{eqn:delta}. Then we have 
\begin{equation*}
    \delta_{\varepsilon}(\hp,p)\leq\sqrt{\frac{\tVar(Z_1)}{np_1^2(\varepsilon-n\tp_2)}}+\frac{p_2}{p}\to 0\text{ as }\gamma\to\infty.
\end{equation*}
Hence $Z$ is strongly probabilistically efficient. 
\label{thm:strong_probabilistic_efficiency}
\end{theorem}
\proof{Proof of Theorem \ref{thm:strong_probabilistic_efficiency}.}
Suppose that we sample $\omega_1,\dots,\omega_n$ under $\tP$. Let
\begin{equation*}
    \hp=\frac{1}{n}\sum_{i=1}^nI_{\cA_{\gamma}}(\omega_i)\frac{dP}{d\tP}(\omega_i)
\end{equation*} 
and 
\begin{equation*}
    \hp_j=\frac{1}{n}\sum_{i=1}^nI_{\cA_{\gamma}^j}(\omega_i)\frac{dP}{d\tP}(\omega_i),j=1,2.
\end{equation*}
Clearly $\hp=\hp_1+\hp_2$. For simplicity, denote 
\begin{align}
    \delta&:=\sqrt{\frac{\tVar(Z_1)}{np_1^2(\varepsilon-n\tp_2)}}+\frac{p_2}{p},\label{delta def}\\
    \tilde{\delta}&:=\frac{p_2}{p}.\label{tilde delta def}
\end{align}
Then we have that 
\begin{align*}
    \tP(|\hp-p|>\delta p) &\leq \tP(|\hp_1-p_1|>(\delta-\tilde{\delta})p \text{ or }|\hp_2-p_2|>\tilde{\delta}p)\\
    & \leq \tP(|\hp_1-p_1|>(\delta-\tilde{\delta})p)+\tP(|\hp_2-p_2|>\tilde{\delta}p)\\
    &\leq \tP(|\hp_1-p_1|>(\delta-\tilde{\delta})p_1)+\tP(\hp_2>0)\\
    &\leq \frac{\tVar(Z_1)}{n(\delta-\tilde{\delta})^2p_1^2}+n\tilde{p}_2\\
    &=\varepsilon.
\end{align*}
where the second inequality follows from a union bound, the third inequality follows from $p_1\leq p$ and that $|\hp_2-p_2|>\tilde{\delta}p$ implies $\hp_2>0$, the fourth inequality follows from Chebyshev's inequality in the first term and a union bound in the second term, and the last equality follows from the definitions in \eqref{delta def} and \eqref{tilde delta def}. Thus $\delta_{\varepsilon}(\hp,p)\leq\delta$. Finally, $\delta\to0$ by a direct use of the assumptions.\hfill\Halmos
\endproof

Similarly, if we relax the assumption that $p_1/p\to 1$ as $\gamma\to\infty$, we get sufficient conditions for weak probabilistic efficiency:

\begin{theorem}[Achieving weak probabilistic efficiency]
Suppose that $\{\cA_{\gamma}\}_{\gamma}$ is an indexed family of rare events and $p=P(\cA_{\gamma})\to 0$ as $\gamma\to\infty$. We write $\cA_{\gamma}=\cA_{\gamma}^1\cup \cA_{\gamma}^2$ where $\cA_{\gamma}^1$ and $\cA_{\gamma}^2$ are two disjoint events. Denote $p_j=P(\cA_{\gamma}^j),j=1,2$. Assume that 
\begin{enumerate}
    \item $\liminf_{\gamma\to\infty}\frac{p_1}{p}= c$ where $0<c\leq 1$;
    \item We have an asymptotically efficient IS estimator for $p_1$ obtained from $Z_1=I_{\cA_{\gamma}^1}\frac{dP}{d\tP}$ under $\tP$. This implies that there exists $n=n(\gamma)$ growing subexponentially in $-\log p$ such that  $\frac{\tVar(Z_1)}{np_1^2}\to 0$ as $\gamma\to\infty$;
    \item $\tp_2:=\tP(\cA_{\gamma}^2)$ satisfies that $n\tp_2\to0$ as $\gamma\to\infty$.
\end{enumerate}
Let $\hp$ be the sample mean of $n$ independent replications of $Z=I_{\cA_{\gamma}}\frac{dP}{d\tP}$ under $\tP$. For any $\varepsilon>0$, define $\delta_{\varepsilon}(\hp,p)$ as in \eqref{eqn:delta}. Then we have 
\begin{equation*}
    \limsup_{\gamma\to\infty}\delta_{\varepsilon}(\hp,p)\leq\limsup_{\gamma\to\infty}\left(\sqrt{\frac{\tVar(Z_1)}{np_1^2(\varepsilon-n\tp_2)}}+\frac{p_2}{p}\right)=1-c<1.
\end{equation*}
Hence $Z$ is weakly probabilistically efficient. 
\label{thm:weak_probabilistic_efficiency}
\end{theorem}
\proof{Proof of Theorem \ref{thm:weak_probabilistic_efficiency}.}
Following the proof of Theorem \ref{thm:strong_probabilistic_efficiency}, we still get
\begin{equation*}
    \delta_{\varepsilon}(\hp,p)\leq\delta:=\sqrt{\frac{\tVar(Z_1)}{np_1^2(\varepsilon-n\tp_2)}}+\frac{p_2}{p}.
\end{equation*}
Under the conditions of Theorem \ref{thm:weak_probabilistic_efficiency}, now $\limsup_{\gamma\to\infty}\delta= 1-c<1$. By the definition, $Z$ is weakly probabilistically efficient.\hfill\Halmos
\endproof

We note that in Theorems \ref{thm:strong_probabilistic_efficiency} and \ref{thm:weak_probabilistic_efficiency}, $\cA_{\gamma}^1$ and $\cA_{\gamma}^2$ could be very general events. In particular, when we use dominating points to decompose the rare-event set, $\cA_{\gamma}^1$ and $\cA_{\gamma}^2$ are not necessarily each governed by only one dominating point but could be more as long as the assumptions hold. Moreover, there can be multiple ways to split $\cA_{\gamma}$ into $\cA_{\gamma}^1$ and $\cA_{\gamma}^2$, and as long as one of these ways validates the assumptions in Theorem \ref{thm:strong_probabilistic_efficiency} or \ref{thm:weak_probabilistic_efficiency} then probabilistic efficiency is guaranteed. This provides flexibility in using Theorems \ref{thm:strong_probabilistic_efficiency} and \ref{thm:weak_probabilistic_efficiency}; Sections \ref{sec:experiment_validate_random_walk} and \ref{sec:experiment_validate_mnist} will demonstrate this in some specific examples.

According to the theorems, supposing that we have found some dominating points while the remaining ones are known to be less significant and ``far from'' the current ones, we could simply use the current mixture IS distribution instead of keep searching. The remaining question is how we could detect that the remaining dominating points are negligible, i.e. the assumptions of the theorem are satisfied. Besides, probabilistic efficiency only implies that the point estimate is reliable in some sense. This raises questions on inference such as the construction of valid CIs. In the next section, we will make these discussions precise and show our answers under minimal assumptions in the standard Gartner-Ellis regime. 

\section{Probabilistically Efficient Estimation in the Gartner-Ellis Regime}
\label{sec:guarantees_GE}


We study probabilistically efficient IS in the widely used Gartner-Ellis regime introduced in Section~\ref{sec:background}. Our key result is that probabilistic efficiency can be readily achieved by using only the most significant dominating points, under essentially no more assumptions than what is needed to derive the Gartner-Ellis large deviations asymptotic. 



We first consider the case where there is only one most significant dominating point, which is a common scenario (e.g., in all the examples in Section \ref{sec:prelim numerics}):

\begin{theorem}[Using the most significant dominating point is probabilistically efficient]
Consider the problem of estimating $p=P(\frac{1}{\gamma}X_{\gamma}\in\cE)$. Suppose that Assumptions \ref{asm:mu_GE} and \ref{asm:E_GE} hold. Suppose also that the dominating set $A$ has finite cardinality with a unique most significant dominating point $a$, i.e., $I(a)<I(\tilde a)$ for all other $\tilde a\in A$. Then the IS distribution $\tilde P$ given by the exponential tilting towards $a$, i.e.,
\begin{equation}
\frac{d\tilde P}{dP}(\omega)=e^{s_{a}^\top X_\gamma-\gamma\mu_\gamma(s_{a})}\label{IS single}
\end{equation}
is strongly probabilistically efficient.\label{thm:unique}
\end{theorem}

Theorem \ref{thm:unique} stipulates that we only need the most significant dominating point in constructing an efficient IS. This result, which is in sharp contrast to the established IS recipe that suggests using all dominating points, explains the good empirical performance of the ``poor" estimators in Section \ref{sec:prelim numerics}. Moreover, the proposal in Theorem \ref{thm:unique} is in closer line with the Gartner-Ellis asymptotic theory, in that the use of the most significant dominating point, which is also the minimizer of the rate function (recall Theorem \ref{thm:global_minimizer}), governs both the large deviations asymptotic and the construction of efficient IS. Lastly, regarding the assumptions needed, the only additional condition beyond the standard Gartner-Ellis assumptions (i.e., Assumptions \ref{asm:mu_GE} and \ref{asm:E_GE}) is the finite cardinality of the dominating set. In fact, if we have multiple most significant points, we have a natural generalization:


\begin{theorem}
[Mixing most significant dominating points]
Consider the problem of estimating $p=P(\frac{1}{\gamma}X_{\gamma}\in\cE)$. Suppose that Assumptions \ref{asm:mu_GE} and \ref{asm:E_GE} hold. Suppose also that the dominating set $A$ has finite cardinality with $k$ most significant dominating points $a_1,\ldots,a_k$, i.e., $I(a_1)=\cdots=I(a_k)<I(a')$ for all $a'\in A\setminus\{a_1,\ldots,a_k\}$. Then the IS distribution $\tilde P$ given by the mixture of exponential tiltings towards $a_1,\ldots,a_k$, i.e.,
\begin{equation}
    \frac{d\tP}{dP}(\omega)=\sum_{i=1}^k\alpha_i e^{s_{a_i}^\top X_{\gamma}-\gamma\mu_{\gamma}(s_{a_i})}
    \label{eqn:IS_distribution_GE}
\end{equation}
where $\sum_{i=1}^k\alpha_i=1,\alpha_i>0,\forall i$, is strongly probabilistically efficient.\label{thm:general_point_estimate_GE}
\end{theorem}

That is, we use mixture to account for all the most significant dominating points when there are multiple of them. The proofs of Theorems \ref{thm:unique} and \ref{thm:general_point_estimate_GE} amount to verifying the assumptions in Theorem \ref{thm:strong_probabilistic_efficiency} using the Gartner-Ellis conditions. In particular, Conditions 1 and 2 in Theorem \ref{thm:strong_probabilistic_efficiency} can be routinely verified, while Condition 3 is checked by showing that $\tilde p_2$ is in fact exponentially decaying in $\gamma$, which requires an application of the Gartner-Ellis theorem under the IS distribution. The verification of Condition 3 especially reveals a key phenomenon that, under the exponential tilting to the most significant dominating point(s), the probability of an IS sample hitting onto the ``backyards" of other dominating points is exponentially small, which in turn fulfills the notion of probabilistic efficiency.



\proof{Proofs of Theorems \ref{thm:unique} and \ref{thm:general_point_estimate_GE}.}
We focus on Theorem \ref{thm:general_point_estimate_GE} since Theorem \ref{thm:unique} is a special case therein. It suffices to verify all the assumptions in Theorem \ref{thm:strong_probabilistic_efficiency}. For $\cA_{\gamma}=\{\frac{1}{\gamma}X_{\gamma}\in\cE\}$, by Theorem \ref{thm:GE}, $p=P(\cA_{\gamma})$ satisfies that $\lim_{\gamma\to\infty}\frac{1}{\gamma}\log p=-I(\cE)<0$, so $p\to 0$ as $\gamma\to\infty$. If $k=|A|$, then the IS estimator from \eqref{eqn:IS_distribution_GE} already uses all the dominating points and thus is asymptotically efficient by Proposition \ref{classical AE}. Hence, from now on, we assume that $k<|A|$. For convenience, we denote $a_{k+1}$ as a next most significant point other than $a_1,\ldots,a_k$, i.e., $I(a_{k+1})>I(a_1)=\cdots=I(a_k)$ and $I(a_{k+1})\leq I(a')$ for all $a'\in A\setminus\{a_1,\ldots,a_k\}$. We split $\cE$ into $\cE_1=\cE\cap\bigcup_{i=1}^k\{x\in\R^d:s_{a_i}^\top(x-a_i)\geq 0\}$ and $\cE_2=\cE\setminus\cE_1$, and define $\cA_j=\{\frac{1}{\gamma}X_{\gamma}\in\cE_j\}$, $p_j=P(\cA_j)$ for $j=1,2$. 

First, by Theorem \ref{thm:GE}, we have that $\limsup_{\gamma\to\infty}\frac{1}{\gamma}\log p_2\leq -I(\overline{\cE_2})$. Since $\overline{\cE_2}\subset \bigcup_{j=k+1}^{|A|}\{x:s_{a_j}^\top(x-a_j)\geq 0\}$, we know that $I(\overline{\cE_2})\geq I(\bigcup_{j=k+1}^{|A|}\{x:s_{a_j}^\top(x-a_j)\geq 0\})=I(a_{k+1})$, and hence $\limsup_{\gamma\to\infty}\frac{1}{\gamma}\log p_2\leq-I(a_{k+1})<-I(\cE)$. As a result, $p_2/p\to 0$ as $\gamma\to\infty$. This verifies Assumption 1 in Theorem \ref{thm:strong_probabilistic_efficiency}.

Second, by the definition, $\{a_1,\dots,a_k\}$ is a dominating set for $\cE_1$, so the IS estimator $Z_1=I(\frac{1}{\gamma}X_{\gamma}\in\cE_1)\frac{dP}{d\tP}(\omega)$ is asymptotically efficient by Proposition \ref{classical AE}. This verifies Assumption 2 in Theorem \ref{thm:strong_probabilistic_efficiency}.

Third, we would prove that $\tp_2$ decays exponentially in $\gamma$ (hence also exponentially in $-\log p$), and hence $n\tp_2\to 0$ for subexponentially growing $n$ which verifies Assumption 3 in Theorem \ref{thm:strong_probabilistic_efficiency}. Indeed, we have 
\begin{equation*}
    \tp_2=\sum_{i=1}^k\alpha_i\tP_i\left(\frac{1}{\gamma}X_{\gamma}\in\cE_2\right)
\end{equation*}
where $\frac{d\tP_i}{dP}(\omega)=e^{s_{a_i}^\top X_{\gamma}-\gamma\mu_{\gamma}(s_{a_i})}$. Denote $\tilde E_i$, $\tilde\mu_{\gamma,i}$ and $\tilde I_i$ as the corresponding expectation, scaled logarithmic moment generating function and rate function under $\tP_i$. Then, under $\tP_i$, we have 
\begin{align*}
   \tilde{\mu}_{\gamma,i}(x)&=\frac{1}{\gamma}\log\tE_i\left(e^{x^\top X_{\gamma}}\right)=\frac{1}{\gamma}\log E\left(e^{x^\top X_{\gamma}+s_{a_i}^\top X_{\gamma}-\gamma\mu_{\gamma}(s_{a_i})}\right)=\mu_{\gamma}(x+s_{a_i})-\mu_{\gamma}(s_{a_i})
\end{align*}
and thus $\tilde{\mu}_i(x)=\lim_{\gamma\to\infty}\tilde{\mu}_{\gamma,i}(x)=\mu(x+s_{a_i})-\mu(s_{a_i})$. Then the rate function is 
\begin{align*}
    \tilde{I}_i(y)&=\sup_{x\in\R^d}\{x^\top y-\tilde{\mu}_i(x)\}\\
    &=\sup_{x\in\R^d}\{x^\top y-\mu(x+s_{a_i})+\mu(s_{a_i})\}\\
    &=\sup_{x\in\R^d}\{(x+s_{a_i})^\top y-\mu(x+s_{a_i})\}-s_{a_i}^\top y+\mu(s_{a_i})\\
    &=I(y)-s_{a_i}^\top y+\mu(s_{a_i}).
\end{align*}
For any $y\in\overline{\cE_2}$, we have that $I(y)\geq I(a_{k+1})$ and that $s_{a_i}^\top (y-a_i)\leq 0$ for $i=1,\dots,k$, and thus $\tilde{I}_i(y)=I(y)-s_{a_i}^\top y+\mu(s_{a_i})\geq I(a_{k+1})-s_{a_i}^\top a_i+\mu(s_{a_i})=I(a_{k+1})-I(a_i)$. Therefore, $\tilde{I}_i(\overline{\cE_2})\geq I(a_{k+1})-I(a_i)=I(a_{k+1})-I(a_1)>0$. From the above derivations, Assumption \ref{asm:mu_GE} still holds for $\tP_i$. By Theorem \ref{thm:GE}, 
\begin{equation*}
    \limsup_{\gamma\to\infty}\frac{1}{\gamma}\log\tP_i\left(\frac{1}{\gamma}X_{\gamma}\in\cE_2\right)\leq -\tilde{I}_i(\overline{\cE_2})\leq -(I(a_{k+1})-I(a_1))<0.
\end{equation*}
Overall, we have $\tp_2$ decays exponentially in $\gamma$.

Now we have verified all the assumptions in Theorem \ref{thm:strong_probabilistic_efficiency}.
\hfill\halmos
\endproof

We comment that if we know any one of the most significant points, say $a$, among several such points, satisfies $p_1/p\to c$ for some $0<c\leq1$, where $p_1=P(\mathcal A_\gamma\cap\{s_a^\top(\frac{1}{\gamma}X_{\gamma}-a)\geq0\})$ is the rare-event probability ``contributed" from $a$, then using the IS that exponentially tilts only to $a$, i.e., \eqref{IS single}, is weakly probabilistically efficient. This can be shown by a similar argument to the proofs of Theorems \ref{thm:unique} and \ref{thm:general_point_estimate_GE} above. Such an approach is in contrast to using the IS mixture in \eqref{eqn:IS_distribution_GE} suggested by Theorem \ref{thm:general_point_estimate_GE} that achieves strong, instead of only weak, probabilistic efficiency. Nonetheless, knowing $p_1/p\to c$ typically requires information on the multiplicative factor in front of the exponential decay dictated by the large deviations rate function, which in turn requires derivation of \emph{exact asymptotic} that is only known for a relatively small number of problems. 

Next, besides point estimates, we investigate inference using probabilistically efficient IS estimators, in particular how to construct (asymptotically) valid CIs. First, we consider the interval 
\begin{equation}\label{eq:loose_CI_PE}
\mathcal I_1=\left[\hat p- \left(\sqrt{\frac{2\hat{V}\log(4/\alpha)}{n}}+\frac{7\log(4/\alpha)M_{\gamma}}{3(n-1)}\right),\ \hat p+ \left(\sqrt{\frac{2\hat{V}\log(4/\alpha)}{n}}+\frac{7\log(4/\alpha)M_{\gamma}}{3(n-1)}\right)\right]
\end{equation}
where $\hat{V}$ is the sample variance and $M_{\gamma}=\max_{i=1,\dots,k}\left\{\frac{1}{\alpha_i}e^{-\gamma(s_{a_i}^\top a_i-\mu_{\gamma}(s_{a_i}))}\right\}$ is deterministic. The following theorem provides an asymptotic coverage guarantee for this CI:



\begin{theorem}[Constructing confidence intervals with probabilistically efficient estimators]
Under the same setting as Theorem \ref{thm:general_point_estimate_GE}, suppose we sample $X^{(1)},\dots,X^{(n)}$ i.i.d. from $\tP$ and let $Z^{(i)}=I(\frac{1}{\gamma}X^{(i)}\in\cE)\frac{dP}{d\tP},i=1,\dots,n$. Use $\hat p$ and $\hat{V}$ to respectively denote the sample mean and sample variance of $Z^{(i)}$'s. If $n$ is subexponentially growing in $-\log p$ (or $\gamma$) as $\gamma\to\infty$, then, for any $0<\alpha<1$,
\begin{equation*}
\liminf_{\gamma\to\infty}\tP\left(p\in\mathcal I_1\right)\geq1-\alpha
\end{equation*}
where $M_{\gamma}=\max_{i=1,\dots,k}\left\{\frac{1}{\alpha_i}e^{-\gamma(s_{a_i}^\top a_i-\mu_{\gamma}(s_{a_i}))}\right\}$. 
\label{thm:general_interval1_GE}
\end{theorem}

In Theorem \ref{thm:general_interval1_GE}, note that even if we neglect the higher-order term (in terms of $n$) $\frac{7\log(4/\alpha)M_{\gamma}}{3(n-1)}$, the CI half-width is $\sqrt{2\log(4/\alpha)}$ times $\sqrt{\frac{\hat{V}}{n}}$, which is more conservative than the Central Limit Theorem (CLT) based interval 
\begin{equation}\label{eq:tight_CI_PE}
    \mathcal I_2=\left[\hp- z_{1-\alpha/2}\sqrt{\frac{\hat{V}}{n}},\ \hp+ z_{1-\alpha/2}\sqrt{\frac{\hat{V}}{n}}\right]
\end{equation}
where $z_{1-\alpha/2}$ is the $(1-\alpha/2)$-quantile of the standard normal distribution. For instance, when $\alpha=0.05$, we have $\sqrt{2\log(4/\alpha)}\approx 2.96$, while $z_{1-\alpha/2}\approx 1.96$. Our next theorem shows that, under stronger conditions, the CLT-based CI \eqref{eq:tight_CI_PE} is also asymptotically valid.
\begin{theorem}[Constructing tight confidence intervals with probabilistically efficient estimators]
Under the same setting as Theorem \ref{thm:general_point_estimate_GE}, suppose we sample $X^{(1)},\dots,X^{(n)}$ i.i.d. from $\tP$ and let $Z^{(i)}=I(\frac{1}{\gamma}X^{(i)}\in\cE)\frac{dP}{d\tP},i=1,\dots,n$. Use $\hat p$ and $\hat{V}$ to respectively denote the sample mean and sample variance of $Z^{(i)}$'s. In this case, we could choose $n$ subexponentially growing in $-\log p$ (or $\gamma$) such that $\frac{M_{\gamma}^2}{n\tVar(Z_1^{(1)})}\to 0$ and $\frac{\tE^2|Z_1^{(1)}-p_1|^3}{n\tVar^3(Z_1^{(1)})}\to 0$ as $\gamma\to\infty$ where $M_{\gamma}$ is as defined in Theorem~\ref{thm:general_interval1_GE}. Then, for any $0<\alpha<1$,
\begin{equation*}
\liminf_{\gamma\to\infty}\tP\left(p\in\mathcal I_2\right)\geq1-\alpha.
\end{equation*}
\label{thm:general_interval2_GE}
\end{theorem}



Theorem \ref{thm:general_interval2_GE} tightens the interval in Theorem \ref{thm:general_interval1_GE} to using the CLT-based critical value $z_{1-\alpha/2}$ with a more careful choice of sample size $n$.

Finally, we prove that if we use all the dominating points in the mixture, so that the estimator satisfies the classical notion of asymptotic efficiency, then, under conditions similar to Theorem \ref{thm:general_interval2_GE}, the CLT-based interval possesses an even stronger guarantee that the asymptotic coverage probability is exactly $1-\alpha$.

\begin{theorem}[Asymptotically exact confidence intervals with asymptotically efficient estimators]
Consider the problem of estimating $p=P(\frac{1}{\gamma}X_{\gamma}\in\cE)$. Suppose that Assumptions \ref{asm:mu_GE} and \ref{asm:E_GE} hold, and the dominating set is finite. The IS estimator is  $Z=I(\frac{1}{\gamma}X_{\gamma}\in\cE)\frac{dP}{d\tP}(\omega)$ under $\tP$ given by \eqref{eqn:full_IS_GE}. We sample $X^{(1)},\dots,X^{(n)}$ i.i.d. from $\tP$ and let $Z^{(i)}=I(\frac{1}{\gamma}X^{(i)}\in\cE)\frac{dP}{d\tP},i=1,\dots,n$. Use $\hat p$ and $\hat{V}$ to respectively denote the sample mean and sample variance of $Z^{(i)}$'s. In this case, we could choose $n$ at least subexponentially growing in $-\log p$ (or $\gamma$) such that $\frac{M_{\gamma}^2}{n\tVar(Z^{(1)})}\to 0$ and $\frac{\tE^2|Z^{(1)}-p|^3}{n\tVar^3(Z^{(1)})}\to 0$ as $\gamma\to\infty$ where $M_{\gamma}=\max_{i=1,\dots,r}\left\{\frac{1}{\alpha_i}e^{-\gamma(s_{a_i}^\top a_i-\mu_{\gamma}(s_{a_i}))}\right\}$. Then, for any $0<\alpha<1$,
\begin{equation*}
\lim_{\gamma\to\infty}\tP\left(p\in\mathcal I_2\right)=1-\alpha.
\end{equation*}
\label{thm:tight_interval_full_IS}
\end{theorem}

We make several remarks regarding the properties of the CLT-based CI $\mathcal I_2$ in Theorems \ref{thm:general_interval2_GE} and \ref{thm:tight_interval_full_IS}. First, it appears that probabilistically efficient samples sacrifice some looseness in terms of CI coverage compared to asymptotically efficient samples, as the guarantee is  valid in Theorem \ref{thm:general_interval2_GE} but exact in Theorem \ref{thm:tight_interval_full_IS}. Second, in Theorem \ref{thm:general_interval2_GE}, like Theorem \ref{thm:general_interval1_GE}, the sample size $n$ is required to be not overly big, manifested by the subexponential growth requirement. This is in contrast to Theorem \ref{thm:tight_interval_full_IS} that does not impose any upper bound on $n$. This ties to the key idea of probabilistic efficiency that, when the sample size is not overly big, there is a negligible chance of any sample hitting the rare-event region not corresponding to the most significant points. Thus the CI constructed from a probabilistically efficient estimator, much like the point estimate, is valid only when the sample size is not overly big, while asymptotically efficient estimators do not impose such a restriction. Lastly, we see the requirement on $n$ given by $\frac{M_{\gamma}^2}{n\tVar(Z^{(1)})}\to 0$ and $\frac{\tE^2|Z^{(1)}-p|^3}{n\tVar^3(Z^{(1)})}\to 0$ in Theorems \ref{thm:general_interval2_GE} and \ref{thm:tight_interval_full_IS}. While these conditions can be difficult to verify in practice, we should note that they are lower bound requirements, and imposed not only for CIs constructed from probabilistically efficient estimators, but also for classical asymptotically efficient estimators as well (to our best knowledge, conditions on the adequacy of sample size to attain CI coverage guarantees for these classical estimators is not known in the literature). In the next section, we will investigate the performances of all these CIs with reasonable sample sizes.

Lastly, to close this section, we briefly note that Algorithm \ref{algo:dominating_points} in Appendix \ref{app:algorithm} shows generally how to identify and compute dominating points, sequentially starting from the most significant one. Moreover, Appendix \ref{sec:guarantees} studies parallel results to this section for an alternative asymptotic regime to Gartner-Ellis that could be suitable for some situations involving highly complex systems.

\section{Further Numerical Experiments and Discussions}\label{sec:additional numerics}

We have shown several examples in Section~\ref{sec:prelim numerics} where IS estimators using only one or a small number of dominating points perform competitively compared with asymptotically efficient IS estimators that use all dominating points. In fact, we have shown in each example in Sections \ref{sec:numerical_glasserman}, \ref{sec:numerical_sum} and \ref{sec:numerical_mnist} that the simple estimator using the most significant dominating point is not asymptotically efficient. In this section, we argue that they are all probabilistically efficient, which is a direct consequence of Theorem \ref{thm:unique}. We then numerically assess the validity of the conditions in Theorem \ref{thm:strong_probabilistic_efficiency}, which forms the underlying basis in justifying probabilistic efficiency. Finally, we test the confidence intervals constructed using our probabilistically efficient estimators discussed in Section \ref{sec:guarantees_GE} and compare with intervals constructed from asymptotically efficient estimators.




\subsection{Verifying Conditions for Probabilistic Efficiency}\label{sec:verification}
We first state the strong probabilistic efficiency of all the proposed estimators that use only the most significant dominating points in Section \ref{sec:prelim numerics}:
\begin{theorem}
Under the problem specifications in Sections \ref{sec:numerical_glasserman}, \ref{sec:numerical_sum} and \ref{sec:numerical_mnist}, the IS estimators that use only the most significant dominating points, namely $\hat\alpha(m)$ in Section \ref{sec:numerical_glasserman}, $X_\gamma$ distributed as $N(\gamma a_1,\gamma \Sigma)$ in Section \ref{sec:numerical_sum} and $N(\gamma a_1,\gamma I)$ in Section \ref{sec:numerical_mnist}, are all strongly probabilistically efficient.\label{PE examples}
\end{theorem}

Next, we validate the underpinning mechanism of how probabilistic efficiency arises in these examples. Note that the main basis of the strong probabilistic efficiency of these estimators, which follows from Theorem \ref{thm:unique}, is Theorem \ref{thm:strong_probabilistic_efficiency}. In particular, Theorem \ref{thm:strong_probabilistic_efficiency} states three conditions that allow one to conclude strong probabilistic efficiency. Among them, the second condition is a property about asymptotic efficiency for an estimator that applies to a more restrictive rare-event set, which has been well-established in the asymptotic efficiency literature (basically, by mixing the exponential tiltings towards all the dominating points associated with the more restrictive rare-event set). Conditions 1 and 3 are more delicate. In the setting with a unique most significant dominating point, say $a$, the former requires a small proportion of the  ``contribution" from the less significant dominating points other than $a$ over the total rare-event probability, i.e., $p_2/p\to0$ where $p_2=P(\mathcal A_\gamma\setminus\{s_a^\top(x-a)\geq0\})$. The latter requires a small probability of sampling any points in the rare-event set that does not belong to the backyard of $a$, i.e., $\tilde P(\text{some of the $n$ samples hits\ }\mathcal A_\gamma\setminus\{s_a^\top(x-a)\geq0\})\to0$ or, as a sufficient condition, $n\tilde p_2\to0$ where $\tilde p_2=\tilde P(\mathcal A_\gamma\setminus\{s_a^\top(x-a)\geq0\})$. Our next goal is to assess the smallness and decreasing trends (as rarity grows) of $p_2/p$ and $\tilde P(\text{some of the $n$ samples hits\ }\mathcal A_\gamma\setminus\{s_a^\top(x-a)\geq0\})$ that drive Theorem \ref{thm:strong_probabilistic_efficiency}. 

\subsubsection{Large Deviations of an I.I.D. Sum.}
For the experiment in Section~\ref{sec:numerical_glasserman}, we use the probabilistically efficient estimator $\hat{\alpha} (m)$. Correspondingly, we have $p_1= P(S_m \geq am)$ and $p_2=P(S_m \leq - am)$. Table~\ref{table:p1_p2_iid_sum} shows these values as $m$ varies, which we approximate respectively by using estimators $\hat \beta_1(m)=\exp (-s_a S_m +m \mu(s_a))I_{\{ S_m \geq am \}}$ and $ \hat \beta_2(m)=  \exp (-s_{-a} S_m' +m \mu(s_{-a})) I_{\{ S_m' \leq -am \}}$, with $s_a, s_{-a}$ defined in Section~\ref{sec:numerical_glasserman}, generated by the same IS samples used in $\hat\beta(m)$. From Table~\ref{table:p1_p2_iid_sum}, we observe that the estimate of $p_2 /p$ is $0.008$ with $m=10$ and decreases to $5.00 \times 10^{-19}$ as $m=100$. This shows that $p_2 /p$ is small and approaches 0 as $m$ increases, which matches Condition 1 in Theorem \ref{thm:strong_probabilistic_efficiency} ($m$ is the rarity parameter here). 

Next, we examine $\tilde{p}_2= \tilde P (S_m \leq - am) $. We generate $10^7$ samples from the strongly probabilistically efficient IS distribution.
We observe that none of the samples fall into $\{S_m \leq - am\}$, which indicates that $\tilde{p}_2$ is extremely small so that $\tilde P(\text{some of the $n$ samples hits\ }\mathcal A_\gamma\setminus\{s_a^\top(x-a)\geq0\})$, with $n=10^4$ in our experiment here, is close to zero. This matches Condition 3 of Theorem \ref{thm:strong_probabilistic_efficiency}.

\begin{table}[]
\caption{Estimates of $p_1$ and $p_2$ for the i.i.d sum example in Section \ref{sec:numerical_glasserman} with $10^4$ samples.}\label{table:p1_p2_iid_sum}
\centering
\begin{tabular}{|l|l|l|l|l|}
\hline
$m$        & 10          & 30          & 50          & 100      \\ \hline
$p_1$   &  $8.33 \times 10^{-3}$ & $1.59 \times 10^{-5}$    & $3.59 \times 10^{-8} $    & $1.33 \times 10^{-14}$  \\ \hline
$p_1/p$ & 0.992 & $\approx 1$ & $\approx 1$ & $\approx 1$        \\ \hline
$p_2$   & $6.56 \times 10^{-5}$    & $3.15 \times 10^{-11}$    & $1.85\times 10^{-17} $    & $6.68 \times 10^{-33}$ \\ \hline
$p_2/p$ & 0.008 & $1.98 \times 10^{-6}$    & $5.16\times 10^{-10}$  & $5.00 \times 10^{-19}$  \\ \hline
$p$ & $ 8.40 \times 10^{-3}$ & $ 1.59 \times 10^{-5}$ & $3.59 \times 10^{-8}$  & $ 1.33 \times 10^{-14}$  \\ \hline
\end{tabular}
\end{table}


\subsubsection{Overshoot Probability of a Random Walk.}\label{sec:experiment_validate_random_walk}
For the experiment in Section \ref{sec:numerical_sum}, we consider the most significant dominating point $a=a_1$ and our probabilistically efficient estimator is the exponential tilting towards the most significant dominating point only, i.e., $X_\gamma$ distributed as $N(\gamma a_1,\gamma \Sigma)$. We define $\cA_\gamma^1=\mathcal A_\gamma\cap\{s_a^\top(x-a)\geq0\}$ and $\cA_\gamma^2=\mathcal A_\gamma\setminus\{s_a^\top(x-a)\geq0\}$, and the corresponding probabilities $p_1= P(\mathcal A_\gamma\cap\{s_a^\top(x-a)\geq0\})$ and $p_2=P(\mathcal A_\gamma\setminus\{s_a^\top(x-a)\geq0\})$. We compute $p_1$, $p_2$, and also the contribution of each of the nine less significant dominating points in $p_2$. More precisely, we define $\mathcal{B}_\gamma^2=\mathcal A_\gamma\cap\{s_{a_2}^\top(x-{a_2})\geq0\} \setminus \{s_{a_1}^\top(x-{a_1})\geq0\}$, $\mathcal{B}_\gamma^3=\mathcal A_\gamma\cap\{s_{a_3}^\top(x-{a_3})\geq0\} \setminus( \{s_{a_1}^\top(x-{a_1})\geq0\}\cup\{s_{a_2}^\top(x-{a_2})\geq0\})$,..., $\mathcal{B}_\gamma^{10}=\mathcal A_\gamma\cap\{s_{a_{10}}^\top(x-{a_{10}})\geq0\} \setminus \left( \cup_{j=1}^9 \{s_{a_j}^\top(x-{a_j})\geq0\}  \right) $, and use $p_{a_2}=P(\mathcal B_\gamma^2) ,..,p_{a_{10}}=P(\mathcal B_\gamma^{10})$ to denote the contribution of dominating points $a_2,...,a_{10}$ (with decreasing significance). For each probability $p_1,p_{a_2},...,p_{a_{10}}$, we construct an IS estimator using the ``corresponding'' dominating points $a_1,...,a_{10}$, e.g., for $p_{a_2}$ the IS distribution is mean shifted to $a_2$. Then we estimate $p_2$ through ${p}_2= {p}_{a_2}+...+{p}_{a_{10}}$ and $p$ through ${p}={p}_1+{p}_2$. Table~\ref{table:p1_p2_random_walk} presents the results estimated using independently generated $10^4$ samples from the corresponding IS distributions. 
We observe that $p_2/p$ has larger values than those in the previous experiment, in that $p_2/p \approx 0.25$ at $\sigma=0.3$ and $p_2/p \approx 0.15$ at $\sigma=0.2$. Nonetheless, $p_2/p$'s value decreases rapidly as $\sigma$ decreases, i.e., the problem becomes rarer, which suggests the trend $p_2/p \to 0$ in Condition 1 of Theorem \ref{thm:strong_probabilistic_efficiency}. Additionally, we observe from the values of $p_{a_2}/p,...,p_{a_{10}}/p$ that the contribution of each less significant dominating point vanishes rapidly with decreasing $\sigma$.

\begin{table}[]
\caption{Estimates of $p_1$, $p_2$ and the contributions of the less significant dominating points for the random walk example in Section \ref{sec:numerical_sum} with $10^4$ samples, where we use $\cA_\gamma^1=\mathcal A_\gamma\cap\{s_a^\top(x-a)\geq0\}$ and $\cA_\gamma^2=\mathcal A_\gamma\setminus\{s_a^\top(x-a)\geq0\}$.}\label{table:p1_p2_random_walk}
\centering
\begin{tabular}{|l|l|l|l|l|l|l|}
\hline
$\sigma$  & 0.3         & 0.28        & 0.26        & 0.24        & 0.22        & 0.2         \\ \hline
$p_1$   & $2.42 \times 10^{-4}$ & $9.07 \times 10^{-5}$     & $2.82 \times 10^{-5}$    & $6.12 \times 10^{-6}$     & $ 9.79 \times 10^{-7}$    & $ 7.92 \times 10^{-8}$     \\ \hline
$p_1/p$ & 0.7434      & 0.7512      & 0.7776      & 0.7960      & 0.8311      & 0.8549      \\ \hline
$p_2$   & $8.36 \times 10^{-5}$     & $ 3.00 \times 10^{-5}$     & $8.05 \times 10^{-6} $    & $1.57 \times 10^{-6}$     & $1.99 \times 10^{-7} $     & $ 1.34 \times 10^{-8}$      \\ \hline
$p_2/p$ & 0.2566      & 0.2488      & 0.2224      & 0.2040      & 0.1689      & 0.1451      \\ \hline
$p$      & $3.26 \times 10^{-4}$ & $1.21 \times 10^{-4}$  & $ 3.62\times 10^{-5}$ & $  7.66 \times 10^{-6} $ & $ 1.18 \times 10^{-6}$  & $ 9.26 \times 10^{-8}$  \\ \hline
$p_{a_2}$   & $  5.29 \times 10^{-5}$    & $1.95 \times 10^{-5}$    & $5.50 \times 10^{-6}$     & $ 1.14 \times 10^{-6}$    & $ 1.56 \times 10^{-7}$    & $ 1.11 \times 10^{-8}$    \\ \hline
$p_{a_2}/p$ & 0.1625   & 0.1611   & 0.1519   & 0.1489   & 0.1323   & 0.1194   \\ \hline
$p_{a_3}$   & $ 2.06 \times 10^{-5}$     & $ 7.46 \times 10^{-6}$     & $ 1.88 \times 10^{-6}$     & $ 3.32 \times 10^{-7}$     & $ 3.57 \times 10^{-8}$     & $ 2.06 \times 10^{-9}$     \\ \hline
$p_{a_3}/p$ & 0.0632    & 0.0618    & 0.0519    & 0.0432    & 0.0303    & 0.0223    \\ \hline
$p_{a_4}$   & $ 7.35 \times 10^{-6}$     & $ 2.48 \times 10^{-6}$     & $5.59 \times 10^{-7}$      & $ 7.85 \times 10^{-8}$     & $ 6.69 \times 10^{-9}$     & $2.96 \times 10^{-10}$     \\ \hline
$p_{a_4}/p$ & 0.0226    & 0.0205    & 0.0154    & 0.0102    & 0.0057    & 0.0032    \\ \hline
$p_{a_5}$   & $2.23 \times 10^{-6}$     & $5.67 \times 10^{-7}$     & $1.02 \times 10^{-7}$     & $ 1.19 \times 10^{-8}$    & $ 6.86 \times 10^{-10} $    & $ 1.95 \times 10^{-11} $    \\ \hline
$p_{a_5}/p$ & 0.0069    & 0.0047    & 0.0028    & 0.0015    & $5.82 \times 10^{-4}$ & $2.11 \times 10^{-4}$  \\ \hline
$p_{a_6}$    & $ 4.18 \times 10^{-7}$     & $ 8.07 \times 10^{-8}$    & $ 1.13 \times 10^{-8}$      & $ 1.00 \times 10^{-9}$     & $ 3.53 \times 10^{-11}$    & $ 4.58 \times 10^{-13}$     \\ \hline
$p_{a_6}/p$ & 0.0013    & $6.68 \times 10^{-4}$  & $3.11 \times 10^{-4}$  & $1.30 \times 10^{-4}$  & $3.00 \times 10^{-5}$  & $4.94 \times 10^{-6}$  \\ \hline
$p_{a_7}$   & $ 3.51 \times 10^{-8}$     & $ 5.00 \times 10^{-9}$     & $4.07 \times 10^{-10} $ & $1.84 \times 10^{-11}$     & $ 3.32 \times 10^{-13}$     & $ 1.75 \times 10^{-15}$     \\ \hline
$p_{a_7}/p$ & $1.08 \times 10^{-4}$ & $4.14\times 10^{-5}$ & $1.12 \times 10^{-5}$ & $2.39 \times 10^{-6}$& $2.82 \times 10^{-7}$ & $1.89\times 10^{-8}$ \\ \hline
$p_{a_8}$  & $ 5.43 \times 10^{-10}$     & $3.95 \times 10^{-11}$    & $ 1.43 \times 10^{-12}$      & $ 2.81 \times 10^{-14}$     & $ 1.58 \times 10^{-16}$     & $ 1.98 \times 10^{-19}$     \\ \hline
$p_{a_8}/p$ & $1.67 \times 10^{-6}$ & $3.28 \times 10^{-7}$ & $3.94 \times 10^{-8}$ & $3.66 \times 10^{-9}$ & $1.34 \times 10^{-10}$ & $2.14\times 10^{-12}$ \\ \hline
$p_{a_9}$   & $1.32 \times 10^{-13}$     & $ 2.93 \times 10^{-15}$     & $2.50 \times 10^{-17}$     & $ 7.30 \times 10^{-20}$     & $ 3.43 \times 10^{-23}$     & $1.46 \times 10^{-27}$     \\ \hline
$p_{a_9}/p$ & $4.05 \times 10^{-10}$ & $  2.43 \times 10^{-11}$ & $6.91\times 10^{-13}$ & $9.50\times 10^{-15}$ & $2.91\times 10^{-17}$ & $1.57\times 10^{-20}$ \\ \hline
$p_{a_{10}}$   & $2.35 \times 10^{-24}$     & $ 1.37 \times 10^{-27} $     & $1.36 \times 10^{-31}$     & $ 1.25 \times 10^{-36}$     & $3.95 \times 10^{-43}$     & $ 1.37 \times 10^{-51}$     \\ \hline
$p_{a_{10}}/p$ & $7.21 \times 10^{-21}$ & $1.13\times 10^{-23}$ & $3.75 \times 10^{-27}$ & $1.63 \times 10^{-31}$ & $3.36\times 10^{-37}$ & $1.48 \times 10^{-44}$ \\ \hline
\end{tabular}
\end{table}

Next, we present the probabilities $\tilde p_1=\tilde P(\mathcal A_\gamma\cap\{s_a^\top(x-a)\geq0\})$ and $\tilde p_2= \tilde P(\mathcal A_\gamma\setminus\{s_a^\top(x-a)\geq0\})$ under the probabilistically efficient IS distribution. For $\tilde p_2$, we also present the contributions from the dominating points $a_2,...,a_{10}$, denoted as $\tilde{p}_{a_2},...,\tilde{p}_{a_{10}}$ with $\tilde{p}_{a_i}=\tilde{P}\left( \mathcal A_\gamma  \cap\{s_{a_i}^\top(x-a_i)\geq0\}  \setminus \left( \cup_{j=1}^{i-1} \{s_{a_j}^\top(x-a_j)\geq0\} \right) \right)$ for $i=2,...,10$. The probabilities are estimated using the proportion of samples falling into the corresponding sets based on $10^4$ samples drawn from the probabilistically efficient IS distribution. The results are presented in Table~\ref{table:tilde_p2_random_walk}.
We observe that $\tilde{p}_2$ generally decreases from 0.0129 at $\sigma=0.3$ to 0.0039 with $\sigma=0.2$. Furthermore, the decreasing trends also appear in each individual contribution from the less significant dominating points, where most of the probabilities (e.g. $\tilde p_{a_4},...,\tilde p_{a_{10}}$) already vanish when $\sigma=0.2$. Based on the value of $\tilde{p}_2$, we estimate the probability $\tilde P(\text{some of the $n$ samples hits\ }\mathcal A_\gamma\setminus\{s_a^\top(x-a)\geq0\})$ through $1 - (1- \tilde{p}_2)^n$. We denote this probability as $\tilde{p}_{hit}$ and present the results with $n=10^4$ (the sample size we use in Section~\ref{sec:numerical}) in the last row of Table~\ref{table:tilde_p2_random_walk}. We observe that there are samples falling into $\mathcal A_\gamma\setminus\{s_a^\top(x-a)\geq0\}$ with approximately probability 1 when we use $n=10^4$ samples. This close-to-1 probability, unfortunately, is quite different from what our Condition 3 in Theorem \ref{thm:strong_probabilistic_efficiency} would entail and cannot explain the good performance of our probabilistically efficient estimator.

\begin{table}[]
\caption{Estimates of $\tilde p_1$, $\tilde p_2$ and the contributions of the less significant dominating points under probabilistically efficient IS for the random walk example in Section \ref{sec:numerical_sum} with $10^4$ samples, where we use $\cA_\gamma^1=\mathcal A_\gamma\cap\{s_a^\top(x-a)\geq0\}$ and $\cA_\gamma^2=\mathcal A_\gamma\setminus\{s_a^\top(x-a)\geq0\}$.}\label{table:tilde_p2_random_walk}
\centering
\begin{tabular}{|l|l|l|l|l|l|l|}
\hline
$\sigma$      & 0.3    & 0.28   & 0.26   & 0.24   & 0.22   & 0.2    \\ \hline
$\tilde p_1$       & 0.5057 & 0.4916 & 0.496  & 0.5053 & 0.4969 & 0.4985 \\ \hline
$\tilde p_2$ & 0.0129 & 0.0141 & 0.0086 & 0.0076 & 0.0057 & 0.0039 \\ \hline
$\tilde p_{a_2}$      & 0.0111 & 0.0111 & 0.007  & 0.0062 & 0.0048 & 0.0034 \\ \hline
$\tilde p_{a_3}$       & 0.0011 & 0.0024 & 0.0013 & 0.0013 & 0.0008 & 0.0005 \\ \hline
$\tilde p_{a_4}$      & 0.0006 & 0.0006 & 0.0002 & 0.0001 & 0.0001 & 0      \\ \hline
$\tilde p_{a_5}$      & 0.0001 & 0      & 0.0001 & 0      & 0      & 0      \\ \hline
$\tilde p_{a_6}$      & 0      & 0      & 0      & 0      & 0      & 0      \\ \hline
$\tilde p_{a_7}$       & 0      & 0      & 0      & 0      & 0      & 0      \\ \hline
$\tilde p_{a_8}$       & 0      & 0      & 0      & 0      & 0      & 0      \\ \hline
$\tilde p_{a_9}$       & 0      & 0      & 0      & 0      & 0      & 0      \\ \hline
$\tilde p_{a_{10}}$      & 0      & 0      & 0      & 0      & 0      & 0      \\ \hline
$\tilde p_{hit}$      & $\approx 1$      & $\approx 1$      & $\approx 1$      & $\approx 1$      & $\approx 1$      & $\approx 1$      \\ \hline
\end{tabular}
\end{table}


To this end, we verify the conditions in Theorem~\ref{thm:strong_probabilistic_efficiency} using an alternative construction of $\cA_\gamma^1$ and $\cA_\gamma^2=\cA_\gamma \setminus \cA_\gamma^1$. Here, in our previous construction, we have chosen the $\cA_\gamma^1$ to be the half-space cut by a dominating point and it turns out that the corresponding $\tilde{p}_2$ is not small and thus the condition of Theorem~\ref{thm:strong_probabilistic_efficiency} appears to fail. However, as discussed right after Theorem \ref{thm:weak_probabilistic_efficiency}, our main theorems allow more flexibility in choosing our $\cA_\gamma^1$, and as long as we find a suitable way to construct $\cA_\gamma^1$ to satisfy the needed conditions, Theorem~\ref{thm:strong_probabilistic_efficiency} can be used to explain our estimator's good performance. 

Here is how we can construct a suitable alternative $\cA_\gamma^1$ for Theorem~\ref{thm:strong_probabilistic_efficiency}. From the proofs of Propositions \ref{prop:example2} and \ref{prop:example3}, we know that our probabilistically efficient estimator is not asymptotically efficient if and only if $\min_{x\in\cA_\gamma}(x+a_1)^\top\Sigma^{-1}(x+a_1)<4 a_1^\top\Sigma^{-1} a_1$. In other words, if we split the rare-event set $\cA_\gamma$ into two parts, say $\cA^1_\gamma=\{\cA_\gamma \cap \{(x+a_1)^\top\Sigma^{-1}(x+a_1) \geq 4 a_1 ^\top\Sigma^{-1} a_1 \} \}$ and $\cA^2_\gamma=\{\cA_\gamma \setminus \{(x+a_1 )^\top\Sigma^{-1}(x+a_1) \geq 4 a_1 ^\top\Sigma^{-1} a_1 \} \}$, then our probabilistically efficient estimator is asymptotically efficient for estimating $P(\cA^1_\gamma)$ because $(x+a_1)^\top\Sigma^{-1}(x+a_1) \geq 4 a_1 ^\top\Sigma^{-1} a_1 $ for all $x\in \cA^1_\gamma$. This implies that, with these choices of $\cA^1_\gamma$ and $\cA^2_\gamma$, we satisfy Condition 2 in Theorem~\ref{thm:strong_probabilistic_efficiency} (since the IS estimator using the most significant dominating point is asymptotically efficient for estimating $P(\cA^1_\gamma)$). 

Next we check Conditions 1 and 3 in Theorem~\ref{thm:strong_probabilistic_efficiency}. We define $p_1= P(\cA^1_\gamma)$, $p_2= P(\cA^2_\gamma)$,  $\tilde  p_1= \tilde P(\cA^1_\gamma)$, and  $\tilde{p}_2= \tilde P(\cA^2_\gamma)$ for our newly constructed $\cA_\gamma^1$ and $\cA_\gamma^2$. We first show $p_2/p \to 0$ and $n \tilde p_2 \to 0$. We note that $\mathcal A_\gamma \cap \{s_a^\top(x-a)\geq0\} \subseteq \cA^1_\gamma$ because $\{s_a^\top(x-a)\geq0\} \subseteq  \{(x+a_1 )^\top\Sigma^{-1}(x+a_1) \geq 4 a_1 ^\top\Sigma^{-1} a_1 \}$. Hence we also have $\cA^2_\gamma \subseteq A_\gamma\setminus\{s_a^\top(x-a)\geq0\}$, which leads to $p_2 \leq p^{old}_2$ and $\tilde  p_2 \leq \tilde  p^{old}_2$ where $p^{old}_2$ and $\tilde  p^{old}_2$ refer to the $p_2$ and $\tilde p_2$ evaluated using our old constructions $\cA_\gamma^1=\mathcal A_\gamma\cap\{s_a^\top(x-a)\geq0\}$ and $\cA_\gamma^2=\mathcal A_\gamma\setminus\{s_a^\top(x-a)\geq0\}$. By $p^{old}_2/p \to 0$ and $n \tilde p^{old}_2 \to 0$ as $\gamma \to \infty$ from Theorem~\ref{thm:unique} we have $p_2/p \to 0$ and $n \tilde p_2 \to 0$. That is, our current new construction $\cA_\gamma^1$ for Theorem \ref{thm:strong_probabilistic_efficiency} would satisfy the conditions therein, and we would like to numerically verify especially Conditions 1 and 3.
Indeed, to estimate $p_2$, we construct an IS estimator that mixes the exponential tiltings towards the dominating points for $\cH_i \cap \{(x+a_1)^\top\Sigma^{-1}(x+a_1) \leq 4 a_1 ^\top\Sigma^{-1} a_1 \} $ with $i=1,...,10$. To estimate $\tilde p_1$ and $\tilde p_2$, we directly generate samples from the probabilistically efficient IS distribution. We define $\tilde p_{hit} = \tilde P(\text{some of the $n$ samples hits\ } \cA^2_\gamma)$ and estimate $\tilde p_{hit}$ through $1 - (1-\tilde p_2)^n$ with $n=10^4$. 
The results are presented in Table~\ref{table:tilde_p_prime_2_randomwalk}. We observe that the values of $p_2/p$ are now extremely small (smaller than $10^{-8}$ in all cases). Furthermore, the values of $\tilde p_2$ are also small, which lead to $\tilde p_{hit} < 0.01$ in all cases when $\sigma$ varies from 0.2 to 0.3. These results now justify Conditions 1 and 3 of Theorem \ref{thm:strong_probabilistic_efficiency} and explain the good performance of our probabilistically efficient estimator in the experiment.

\begin{table}[]
\caption{Estimates of $p_2$ with $10^4$ samples, $\tilde p_1$ and $\tilde p_2$ with $10^7$ samples, and $\tilde p_{hit}$ with $n=10^4$ for random walk example in Section \ref{sec:numerical_sum}, where we use $\cA^1_\gamma=\{\cA_\gamma \cap \{(x+a_1)^\top\Sigma^{-1}(x+a_1) \geq 4 a_1 ^\top\Sigma^{-1} a_1 \} \}$ and $\cA^2_\gamma=\{\cA_\gamma \setminus \{(x+a_1 )^\top\Sigma^{-1}(x+a_1) \geq 4 a_1 ^\top\Sigma^{-1} a_1 \} \}$.}\label{table:tilde_p_prime_2_randomwalk}
\centering
\begin{tabular}{|l|l|l|l|l|l|l|}
\hline
$\sigma$      & 0.3    & 0.28   & 0.26   & 0.24   & 0.22   & 0.2    \\ \hline
$p_2$        & $3.85 \times 10^{-13}$ & $2.58 \times 10^{-12}$ & $9.77 \times 10^{-15}$  & $8.58 \times 10^{-18}$ & $2.32 \times 10^{-22}$ & $1.90 \times 10^{-30}$ \\ \hline
$p_2/p$        & $1.19 \times 10^{-9}$ & $2.15 \times 10^{-8}$ & $2.72 \times 10^{-10}$ & $ 1.08 \times 10^{-12}$ & $2.00 \times 10^{-16}$ & $1.89  \times 10^{-23}$\\ \hline
$\tilde p_1$        & 0.5149 & 0.5121 & 0.5098  & 0.5078 & 0.5059 & 0.5038 \\ \hline
$\tilde p_2$  & $3 \times 10^{-7}$ & $5 \times 10^{-7}$ & $7 \times 10^{-7}$ & $1 \times 10^{-7}$ & $1 \times 10^{-7}$ & $4 \times 10^{-7}$ \\ \hline
$\tilde p_{hit}$    & 0.003    & 0.005   & 0.007    & 0.001  & 0.001 & 0.004   \\ \hline
\end{tabular}
\end{table}

\subsubsection{MNIST Example.}\label{sec:experiment_validate_mnist}
Like Section \ref{sec:numerical_sum}, the experiment in Section \ref{sec:numerical_mnist} also uses a probabilistically efficient estimator based on the exponential tilting towards the most significant dominating point only, i.e., $X_\gamma$ distributed as $N(\gamma a_1,\gamma I)$. We define $\cA_\gamma^1=\mathcal A_\gamma\cap\{s_a^\top(x-a)\geq0\}$ and $\cA_\gamma^2=\mathcal A_\gamma\setminus\{s_a^\top(x-a)\geq0\}$, and the corresponding probabilities $p_1= P(\mathcal A_\gamma\cap\{s_a^\top(x-a)\geq0\})$ and $p_2=P(\mathcal A_\gamma\setminus\{s_a^\top(x-a)\geq0\})$ which are shown in Table~\ref{table:p1_p2_mnist}. Note that in this MNIST example, the total number of dominating points is large and unknown. Thus we only present the contribution of the first 10 dominating points, i.e. $a_1,a_2,...,a_{10}$ in Table~\ref{table:p1_p2_mnist}, denoted by $p_{a_2}=P(\mathcal A_\gamma\cap\{s_{a_2}^\top(x-{a_2})\geq0 \setminus \{s_{a_1}^\top(x-a_1)\geq0\}\})$,..., $p_{a_{10}}=P\left( \mathcal A_\gamma\cap\{s_{a_{10}}^\top(x-{a_{10}})\geq0\} \setminus \left( \cup_{j=1}^9 \{s_{a_j}^\top(x-{a_j})\geq0\}  \right)\right)$. 
Again, we estimate each of the probabilities $p_1, p_{a_2},...,p_{a_{10}}$ using the IS estimator with the corresponding dominating point, i.e. the IS distribution is exponentially tilted using dominating points $a_1,..a_{10}$ respectively. 
We borrow the values of $p$ from Table~\ref{table:accurate_mnist} where each estimate is computed using crude Monte Carlo, and we estimate $p_2$ through ${p}_2=p-{p}_1$. 
We observe that the ratio $p_2/p$ decreases from $0.2154$ to $0.1931$ as we decrease the value of $\sigma$ from $0.2$ to $0.17$, i.e., the problem becomes rarer. We also observe that some individual relative contribution slightly increases in this experiment. However, these increases do not affect the decreasing trend of the total relative contribution of the less significant dominating points. For example, $p_{a_2}/p$ and $p_{a_3}/p$ both increase slightly as $\sigma$ decreases (from 0.1074 and 0.0035 with $\sigma=0.2$ to 0.1313 and 0.0042 with $\sigma=0.17$ respectively), but the relative contribution of the rest of the less significant dominating points (excluding the first 10) is 0.0812, 0.0658, 0.0633, and 0.0362 for $\sigma=0.2,0.19,0.18,0.17$ respectively, which vanishes fast as $\sigma$ decreases.

\begin{table}[]
\caption{Estimates of $p_1$, $p_2$ and the contributions of the less significant dominating points for the MNIST example in Section \ref{sec:numerical_mnist} with $10^4$ samples, where we use $\cA_\gamma^1=\mathcal A_\gamma\cap\{s_a^\top(x-a)\geq0\}$ and $\cA_\gamma^2=\mathcal A_\gamma\setminus\{s_a^\top(x-a)\geq0\}$.}\label{table:p1_p2_mnist}
\centering
\begin{tabular}{|l|l|l|l|l|}
\hline
$\sigma$ & 0.2      & 0.19     & 0.18     & 0.17     \\ \hline
$p_1$  & $9.18 \times 10^{-6}$ & $3.34 \times 10^{-6}$ & $ 1.02 \times 10^{-6}$ & $2.54 \times 10^{-7}$ \\ \hline
$p_1/p$  & 0.7846 & 0.7920 & 0.7892 & 0.8069 \\ \hline
$p_2$    & $ 2.52 \times 10^{-6}$  & $ 8.77 \times 10^{-7}$ & $2.72 \times 10^{-7}$  & $6.09 \times 10^{-8}$  \\ \hline
$p_2/p$  & 0.2154   & 0.2080   & 0.2108   & 0.1931   \\ \hline
$p$     & $1.17 \times 10^{-5}$ & $4.22 \times 10^{-6}$  & $ 1.29 \times 10^{-6}$ & $3.15 \times 10^{-7}$ \\ \hline
$p_{a_2}$  & $1.26 \times 10^{-6}$  & $4.95 \times 10^{-7}$  & $1.61 \times 10^{-7}$  & $ 4.14 \times 10^{-8}$  \\ \hline
$p_{a_2}/p$ & 0.1074 & 0.1173 & 0.1245 & 0.1313 \\ \hline
$p_{a_3}$  & $4.12 \times 10^{-8}$  & $1.22 \times 10^{-8}$  & $ 2.93 \times 10^{-9}$  & $ 1.32 \times 10^{-9}$  \\ \hline
$p_{a_3}/p$  & 0.0035 & 0.0029 & 0.0023 & 0.0042 \\ \hline
$p_{a_4}$  & $ 1.24 \times 10^{-7}$ & $4.45 \times 10^{-8}$  & $1.26 \times 10^{-8}$  & $ 3.22 \times 10^{-9}$  \\ \hline
$p_{a_4}/p$  & 0.0106 & 0.0105 & 0.0098 & 0.0102 \\ \hline
$p_{a_5}$  & $ 2.48 \times 10^{-8}$ & $ 9.09 \times 10^{-9}$  & $ 2.81 \times 10^{-9}$  & $ 7.01 \times 10^{-10}$  \\ \hline
$p_{a_5}/p$  & 0.0021 & 0.0022 & 0.0022 & 0.0022 \\ \hline
$p_{a_6}$  & $3.60 \times 10^{-8}$  & $1.31 \times 10^{-8}$  & $ 3.28 \times 10^{-9}$  & $ 8.11 \times 10^{-10}$ \\ \hline
$p_{a_6}/p$ & 0.0031 & 0.0031 & 0.0025 & 0.0026 \\ \hline
$p_{a_7}$  & 0        & 0        & 0        & 0        \\ \hline
$p_{a_8}$  & $8.80 \times 10^{-8}$  & $2.62 \times 10^{-8}$  & $8.03 \times 10^{-9}$  & $2.02 \times 10^{-9}$  \\ \hline
$p_{a_8}/p$ & 0.0075 & 0.0062 & 0.0062 & 0.0064 \\ \hline
$p_{a_9}$  & 0        & 0        & 0        & 0        \\ \hline
$p_{a_{10}}$ & 0        & 0        & 0        & 0        \\ \hline
\end{tabular}
\end{table}

Table~\ref{table:tilde_p2_mnist} presents the estimates of probabilities $\tilde p_1=\tilde P(\mathcal A_\gamma\cap\{s_a^\top(x-a)\geq0\})$ and $\tilde p_2= \tilde P(\mathcal A_\gamma\setminus\{s_a^\top(x-a)\geq0\})$ under the probabilistically efficient IS distribution. The probabilities $\tilde{p}_{a_2},...,\tilde{p}_{a_{10}}$, defined by  $\tilde{p}_{a_i}=\tilde{P}\left( \mathcal A_\gamma  \cap\{s_{a_i}^\top(x-a_i)\geq0\}  \setminus \left( \cup_{j=1}^{i-1} \{s_{a_j}^\top(x-a_j)\geq0\} \right) \right)$, for $i=2,...,10$,  are also shown to illustrate the contributions of the dominating points $a_2,...,a_{10}$ for $\tilde p_2$.
Again, we find that $\tilde p_2$ decreases from $0.0090$ to $0.0085$ as $\sigma$ decreases from 0.2 to 0.17, i.e., the problem becomes rarer. From the individual contribution, we observe that all the probabilities $\tilde p_{a_3},...,\tilde p_{a_{10}}$ decrease rapidly, except $\tilde p_{a_2}$ that slightly increases as $\sigma$ decreases. We use the value of $\tilde{p}_2$ to estimate probability $\tilde{p}_{hit}= \tilde P(\text{some of the $n$ samples hits\ }\mathcal A_\gamma\setminus\{s_a^\top(x-a)\geq0\})$ through $1 - (1- \tilde{p}_2)^n$. The last row in Table~\ref{table:tilde_p2_mnist} presents the results with $n=10^4$ (the sample size we use in Section~\ref{sec:numerical}). Like in Section~\ref{sec:experiment_validate_random_walk}, we observe that there are samples falling into $\mathcal A_\gamma\setminus\{s_a^\top(x-a)\geq0\}$ with approximately probability 1 and hence this result cannot explain the good performance of our probabilistically efficient estimator.

\begin{table}[]
\caption{Estimates of $\tilde p_1$, $\tilde p_2$ and the contributions of the less significant dominating points for the MNIST example under the probabilistically efficient IS in Section \ref{sec:numerical_mnist} with $10^4$ samples.}\label{table:tilde_p2_mnist}
\centering
\begin{tabular}{|l|l|l|l|l|}
\hline
$\sigma$       & 0.2    & 0.19   & 0.18   & 0.17   \\ \hline
$\tilde p_1$        & 0.4728 & 0.4745 & 0.4760  & 0.4775 \\ \hline
$\tilde p_2$ & 0.0090  & 0.0087 & 0.0086 & 0.0085 \\ \hline
$\tilde p_{a_2}$        & 0.0069 & 0.007  & 0.0072 & 0.0072 \\ \hline
$\tilde p_{a_3}$        & 0.0003 & 0.0003 & 0.0003 & 0.0002 \\ \hline
$\tilde p_{a_4}$        & 0.0006 & 0.0006 & 0.0006 & 0.0006 \\ \hline
$\tilde p_{a_5}$        & 0      & 0      & 0      & 0      \\ \hline
$\tilde p_{a_6}$        & 0.0003 & 0.0003 & 0.0001 & 0.0001 \\ \hline
$\tilde p_{a_7}$        & 0      & 0      & 0      & 0      \\ \hline
$\tilde p_{a_8}$        & 0.0006 & 0.0003 & 0.0003 & 0.0003 \\ \hline
$\tilde p_{a_9}$        & 0      & 0      & 0      & 0      \\ \hline
$\tilde p_{a_{10}}$       & 0      & 0      & 0      & 0      \\ \hline
$\tilde p_{hit}$    & $\approx 1$     & $\approx 1$      & $\approx 1$      & $\approx 1$      \\ \hline
\end{tabular}
\end{table}

Similar to Section~\ref{sec:experiment_validate_random_walk}, we consider an alternative construction of $\cA_\gamma^1$ and $\cA_\gamma^2=\cA_\gamma \setminus \cA_\gamma^1$ to explain our performance. From the proofs of Propositions \ref{prop:example2} and \ref{prop:example3}, we know that our probabilistically efficient estimator is not asymptotically efficient if and only if $\min_{x\in\cA_\gamma}(x+a_1)^\top(x+a_1)<4 a_1^\top a_1$. We split the rare-event set $\cA_\gamma$ into two parts, namely $\cA^1_\gamma=\{\cA_\gamma \cap \{(x+a_1)^\top(x+a_1) \geq 4 a_1 ^\top a_1 \} \}$ and $\cA^2_\gamma=\{\cA_\gamma \setminus \{(x+a_1 )^\top(x+a_1) \geq 4 a_1 ^\top a_1 \} \}$, and our probabilistically efficient estimator is asymptotically efficient for estimating the probabilities of $\cA^1_\gamma$ and $\mathcal A_\gamma \cap \{s_a^\top(x-a)\geq0\} \subseteq \cA^1_\gamma$. We define $p_1= P(\cA^1_\gamma)$, $p_2= P(\cA^2_\gamma)$,  $\tilde  p_1= \tilde P(\cA^1_\gamma)$, and  $\tilde{p}_2= \tilde P(\cA^2_\gamma)$. The use of our newly constructed $\cA_\gamma^1$ for Theorem \ref{thm:strong_probabilistic_efficiency} can be theoretically shown to satisfy Conditions 1 and 3 therein like in Section~\ref{sec:experiment_validate_random_walk}. We now check the numerical values of $p_2$ and $\tilde{p}_2$ to verify these conditions. We use the mixture of all 100 dominating points as the IS distribution for estimating $p_2$. We generate $10^6$ samples for $\sigma$ varying from 0.17 to 0.2 and find no samples falling into $\cA^2_\gamma$, which indicates that $p_2$ (and hence $p_2/p$) is extremely small in all cases. 
We generate $10^6$ samples from the probabilistically efficient IS distribution to estimate $\tilde p_2$ and observe no samples hitting  $\cA^2_\gamma$ for the same range of $\sigma$. In this case, $\tilde p_{hit} = \tilde P(\text{some of the $n$ samples hits\ } \cA^2_\gamma)$ with $n=10^4$ would be close to zero due to the extremely small values of $\tilde p_2$.
These results match Conditions 1 and 3 of Theorem \ref{thm:strong_probabilistic_efficiency} and hence explain the good performance of our probabilistically efficient estimator in the experiment.


\subsubsection{Two-sided Overshoot Probability of a Random Walk.}\label{sec:two-sided random walk}
So far we have considered examples on strongly probabilistically efficient estimators. Here, we consider an additional example where we use a weakly probabilistically efficient estimator. We follow the problem setting in Section~\ref{sec:numerical_sum}, where we consider the overshoot probability of the finite-horizon maximum of a random walk. However, we modify the probability of interest as 
\begin{equation}\label{eq:two-sides-random-walk}
p=P\left(\max_{m=1,...,d} | S_m| \geq a\right),
\end{equation}
where we replace $S_m=\sum_{i=1}^m Y_i$ by its absolute value. The rest of the settings are the same as in Section~\ref{sec:numerical_sum}, i.e., we have $Y_i$'s are Gaussian distributed with mean 0, standard deviation $\sigma$, pairwise correlation $-0.02$, and rarity parameter $\gamma=1/\sigma^2\to\infty$. The target rare event is $\left\{\frac{1}{\gamma}X_{\gamma}\in \left( \bigcup_{m=1}^d \cH^+_m \right) \bigcup \left( \bigcup_{m=1}^d \cH^-_m  \right) \right\}$ where $X_{\gamma}=\gamma(Y_1,\dots,Y_d)^\top $, $\cH^+_m = \{x \in \mathbb{R}^d : \sum_{i=1}^m x_i \geq a\}  $, and $\cH^-_m = \{x \in \mathbb{R}^d : \sum_{i=1}^m x_i \leq -a\}  $ with $x_i$ denoting the $i$th element in $x$. In this case, the rate function is still $I(y)=\frac{1}{2}y^\top  \Sigma^{-1} y$ and there are two most significant dominating points $a_{1}=\frac{a\Sigma e_d}{e_d^\top \Sigma e_d}$ and $-a_{1}=-\frac{a\Sigma e_d}{e_d^\top \Sigma e_d}$  where $e_d \in \R^d$ denotes the vector with 1 in all $d$ elements.

For this experiment, we first introduce an asymptotically efficient estimator. From Proposition~\ref{classical AE}, we know that the IS estimator using dominating points $a_1,...,a_d, -a_1,...,-a_d$ with $a_1,...,a_d$ defined in Section~\ref{sec:numerical_sum} is asymptotically efficient. Next, we show the IS estimator using the dominating point $a_1$ is a weakly probabilistically efficient estimator and is not asymptotically efficient:
\begin{theorem}
Under the problem specifications in Section~\ref{sec:numerical_sum} and rare-event probability defined in \eqref{eq:two-sides-random-walk}, the IS estimator that use only one of the most significant dominating points, namely $a_1$ in Section~\ref{sec:numerical_sum}, and $X_\gamma$ distributed as $N(\gamma a_1,\gamma \Sigma)$ in Section \ref{sec:numerical_sum}, is weakly probabilistically efficient but not  asymptotically efficient.\label{thm:Weakly-PE-examples}
\end{theorem}

Compared to the one-sided overshoot example in Section~\ref{sec:numerical_sum}, here the rare-event set has two most significant dominating points $a_1$ and $-a_1$. Because we only use the first one instead of mixing both of the most significant dominating points in our IS, we only have $p_1/p\to1/2$ instead of 0 and thus weak probabilistic efficiency instead of strong probabilistic efficiency holds as guided by Theorem~\ref{thm:weak_probabilistic_efficiency}. 

To empirically verify Theorem~\ref{thm:weak_probabilistic_efficiency}, let us consider a partition of the rare event set $\cA_\gamma = \left\{  \left( \bigcup_{m=1}^d \cH^+_m \right) \bigcup \left( \bigcup_{m=1}^d \cH^-_m  \right)   \right\}$, where we have $\cA^1_\gamma=\{\cA \cap \{(x+a_1)^\top\Sigma^{-1}(x+a_1) \geq 4 a_1 ^\top\Sigma^{-1} a_1 \} \}$ and $\cA^2_\gamma=\{\cA \setminus \{(x+a_1 )^\top\Sigma^{-1}(x+a_1) \geq 4 a_1 ^\top\Sigma^{-1} a_1 \} \}$. We define $p_1= P(\cA^1_\gamma)$, $p_2= P(\cA^2_\gamma)$,  $\tilde  p_1= \tilde P(\cA^1_\gamma)$, and  $\tilde{p}_2= \tilde P(\cA^2_\gamma)$, where $\tilde P$ is the IS distribution exponentially tilted using the dominating point $a_1$. In our experiments, we set $d=10$, fix $a_1,-a_1$ and vary $\sigma$ for different rarity levels. For each case, we generate $10^4$ samples from IS distributions using the above asymptotically efficient estimator and weakly probabilistically efficient estimator. The results are presented in Table~\ref{table:2-side-prob-ci}. We observe that although our weakly probabilistically efficient estimator underestimates the rare-event probability in all considered cases, the estimates have relatively tight CIs and provide a good estimation on the magnitude of the rare-event probability, i.e., the estimates are around 0.5 of the estimates given by the asymptotically efficient estimator.

\begin{table}[]
\caption{Point estimates (and  95\% CIs) from the asymptotically efficient estimator and the weakly probabilistically efficient estimator for the two-sided overshoot probability. ``AE'' denotes the asymptotically efficient estimator 
 and ``PE'' denotes the weakly probabilistically efficient estimator.}\label{table:2-side-prob-ci}
 \centering
\resizebox{\columnwidth}{!}{%
\begin{tabular}{|l|l|l|l|l|l|l|}
\hline
$\sigma$   & 0.3      & 0.28     & 0.26     & 0.24     & 0.22     & 0.2      \\ \hline
AE & $6.15 (\pm 0.47) \times 10^{-4}$ &  $2.52 (\pm 0.20) \times 10^{-4}$ & $7.53 (\pm 0.68)\times 10^{-5}$  & $ 1.56(\pm 0.15)\times 10^{-5}$  & $2.26 (\pm 0.24)\times 10^{-6}$ & $1.93 (\pm 0.24 )\times 10^{-7}$  \\ \hline
PE &  $3.18 (\pm 0.18)\times 10^{-4}$ & $1.19 (\pm 0.08)\times 10^{-4}$  & $ 3.57(\pm 0.29)\times 10^{-5}$  & $ 8.00(\pm 0.64)\times 10^{-6}$ &  $ 1.17(\pm 0.11)\times 10^{-6}$  & $ 9.43(\pm 0.63)\times 10^{-8}$  \\ \hline
AE/PE & 0.516 & 0.470 & 0.474 & 0.512 & 0.517 & 0.488 \\ \hline
\end{tabular}}
\end{table}

In Table~\ref{table:2-sided-validate}, we investigate the numerical values of $p_1$, $p_2$, $\tilde p_1$ and $\tilde p_2$ and check the values of $p_1/p$ and $\tilde p_{hit}= \tilde P(\text{some of the $n$ samples hits\ } \cA^2_\gamma)$ with $n=10^4$. We estimate $p_1$ and $p_2$ using the asymptotically efficient estimator with $10^7$ independently generated samples. For the estimation of $\tilde p_1$ and $\tilde p_2$, we generate $10^7$ samples using the weakly probabilistically efficient IS distribution. We observe that in all cases the values of $p_1/p $ are very close to 1/2. On the other hand, the probabilities of falling into $\cA^2_\gamma$ are all below $10^{-6}$, which lead to $\tilde p_{hit}$ valued smaller than $0.01$. These results verify the conditions in Theorem~\ref{thm:weak_probabilistic_efficiency} that explain the weak probabilistic efficiency of the IS estimator.

\begin{table}[]
\caption{Estimates of $p_1$, $p_2$, $\tilde p_1$ and $\tilde p_2$ with $10^7$ samples, and $\tilde p_{hit}$ with $n=10^4$ for two-sided overshoot probability.}\label{table:2-sided-validate}
\centering
\begin{tabular}{|l|l|l|l|l|l|l|}
\hline
$\sigma$         & 0.3      & 0.28     & 0.26     & 0.24     & 0.22     & 0.2      \\ \hline
$p_1$ &  $3.23 \times 10^{-4}$& $1.20 \times 10^{-4}$ & $3.58 \times 10^{-5}$ &$7.91 \times 10^{-6}$ & $ 1.17\times 10^{-6}$& $9.48 \times 10^{-8}$ \\ \hline
$p_1/p$          & 0.501 & 0.500 & 0.500 & 0.499 & 0.501 & 0.499 \\ \hline
$p_2$            & $ 3.22\times 10^{-4}$ &  $1.20 \times 10^{-4}$& $ 3.59\times 10^{-5}$&  $ 7.94\times 10^{-6}$&  $ 1.16\times 10^{-6}$&  $9.51 \times 10^{-8}$\\ \hline
$p$             & $6.45 \times 10^{-4}$ & $2.40 \times 10^{-4}$ & $ 7.17\times 10^{-5}$ & $1.59 \times 10^{-5}$ & $ 2.33\times 10^{-6}$ &  $ 1.90\times 10^{-7}$\\ \hline
$\tilde p_1$      & 0.515& 0.512 & 0.510 & 0.507 & 0.506 & 0.504 \\ \hline
$\tilde p_2$      & $2 \times 10^{-7}$  & $6 \times 10^{-7}$ & $ 5\times 10^{-7}$ & $6 \times 10^{-7}$ & $4 \times 10^{-7}$ &  $6 \times 10^{-7}$ \\ \hline
$\tilde p_{hit}$ & 0.002 & 0.006 & 0.005 & 0.006 & 0.004 & 0.006 \\ \hline
\end{tabular}
\end{table}

\subsection{Illustration of Confidence Intervals}
We investigate the performances of the CIs proposed in Section~\ref{sec:guarantees_GE}. In particular, we construct CIs \eqref{eq:loose_CI_PE} and \eqref{eq:tight_CI_PE} from probabilistically efficient estimators, namely the IS schemes using only the most significant point in Sections \ref{sec:numerical_glasserman}, \ref{sec:numerical_sum} and \ref{sec:numerical_mnist}. For convenience, we call interval \eqref{eq:loose_CI_PE} the ``loose CI" and interval \eqref{eq:tight_CI_PE} the ``tight CI", since the latter has a shorter length and matches the CLT-based interval. For comparison, we also construct CI \eqref{eq:tight_CI_PE} from asymptotically efficient estimators. In particular, for the settings in Sections \ref{sec:numerical_glasserman} and \ref{sec:numerical_sum}, these estimators are built from mixtures of exponential tiltings towards all the dominating points. For the setting in Section \ref{sec:numerical_mnist}, computing all dominating points requires insurmountable resources (as discussed therein), and so we use the mixture of 100 dominating points as a proxy of an asymptotically efficient estimator (100 is the total number of dominating points we discover using one-week's computation). 




In the experiments, we compare the coverage rates of all three intervals described above. These coverage rates are obtained from a large number of experimental repetitions. Since the ground truths of these problems are unknown, we run a gigantic amount of simulation runs using either asymptotically efficient estimators or crude Monte Carlo to obtain highly accurate estimates, which serve as the ``truths'' when estimating the coverage of the CIs. The exact number of simulation runs used in our ISs, number of experimental repetitions, and number of runs to approximate the ground truths are specified in the discussion of each example below.

\subsubsection{Large Deviations of an I.I.D. Sum.}

For the experiment in Section~\ref{sec:numerical_glasserman}, we use the asymptotically efficient estimator $\hat{\beta} (m)$ to obtain highly accurate estimates for all values of $m$ as the ground truths. These estimates are presented in Table \ref{table:accurate_sum_iid}. Our probabilistically efficient estimator $\hat{\alpha} (m)$ is computed using $10^4$ independently generated samples. From this, we apply CIs \eqref{eq:loose_CI_PE} and \eqref{eq:tight_CI_PE}. We also construct CIs \eqref{eq:tight_CI_PE} using asymptotically efficient estimator $\hat{\beta} (m)$ (which is used to approximate the ground truth) with $10^4$ independently generated samples. We approximate the coverage rates using $10^5$ experimental iterations. Moreover, we compute the average CI width for each type of CIs. The experiment results are shown in Table \ref{table:ci_cover_iid_sum}. 

From Table \ref{table:ci_cover_iid_sum}, we observe that the coverage rates of tight CIs by our probabilistically efficient estimator are close to $95\%$ in three out of the four cases, but is $3\%$ below $95\%$ in one case ($m=10$). On the other hand, the loose CIs are valid but perform conservatively with more than $99\%$ coverage rates and wider average widths in all cases. The tight CIs by asymptotically efficient estimators provide valid coverage in all four cases. In the problems with rarer probabilities (i.e. $m=30,50,100$), the tight CIs by probabilistically efficient estimators perform similarly as the CIs by asymptotically efficient estimators in terms of both CI width and coverage. This shows the competitiveness of CIs using probabilistically efficient estimators for rarer problems.

\begin{table}[]
\caption{Highly accurate point estimates (and 95\% CI) using asymptotically efficient estimators for the problem in Section~\ref{sec:numerical_glasserman}. The estimates are computed with $10^7$ samples for $m=10,30,50$ and $5\times 10^6$ samples for $m=100$.} \label{table:accurate_sum_iid}
\centering
\begin{tabular}{|l|l|l|l|l|}
\hline
$m$                 & 10       & 30       & 50       & 100      \\ \hline
$\hat{p}$ & $8.85 (\pm 0.0084) \times 10^{-3}$ &  $1.58 (\pm 0.0021) \times 10^{-5}$  &  $3.76 (\pm 0.0056) \times 10^{-8}$ &  $1.34 (\pm 0.0034) \times 10^{-14}$ \\ \hline
\end{tabular}
\end{table}

\begin{table}[]
\caption{Coverage rates and average CI widths of the loose confidence intervals (``Loose CI by PE''), the tight confidence intervals (``Tight CI by PE'') for probabilistically efficient estimators, and the tight confidence intervals for asymptotically efficient estimators (''Tight CI by AE'') in the experiments of Section~\ref{sec:numerical_glasserman}.}\label{table:ci_cover_iid_sum}
\centering
\begin{tabular}{|l|l|l|l|l|l|}
\hline
&$m$   & 10     & 30     & 50     & 100    \\ \hline
\multirow{2}{*}{Loose CI by PE} &Coverage Rate  & 0.998 & 0.999 & 0.9994 & 0.999 \\ \cline{2-6}
&  Average Width  & $9.12 \times 10^{-4}$ & $ 2.31 \times 10^{-6}$ & $ 6.37 \times 10^{-9}$ & $ 2.81 \times 10^{-15}$ \\ \hline \hline
\multirow{2}{*}{Tight CI by PE} & Coverage Rate & 0.921 & 0.951 & 0.950 & 0.950 \\ \cline{2-6}
& Average Width  & $5.30 \times 10^{-4}$ & $1.31 \times 10^{-6}$ & $3.54 \times 10^{-9}$ & $1.52 \times 10^{-15}$ \\ \hline \hline
\multirow{2}{*}{Tight CI by AE} & Coverage Rate  & 0.950 &0.960  & 0.949 &  0.950\\ \cline{2-6}
& Average Width  & $5.30 \times 10^{-4}$ & $1.30 \times 10^{-6}$ & $3.55 \times 10^{-9}$ & $1.52 \times 10^{-15}$ \\ \hline 
\end{tabular}
\end{table}

\subsubsection{Overshoot Probability of a Random Walk.}

For the experiment in Section~\ref{sec:numerical_sum}, we use the asymptotically efficient estimator that mixes all dominating points to approximate the ground truths presented in Table \ref{table:accurate_random_walk}. Our probabilistically efficient estimator is computed using $10^4$ independently generated samples. We construct CIs \eqref{eq:loose_CI_PE} and \eqref{eq:tight_CI_PE} based on this estimator. For comparison we also construct CI \eqref{eq:tight_CI_PE} from asymptotically efficient estimator using $10^4$ samples independently generated from the ones used to approximate the ground truth. We use $10^5$ experimental repetitions to estimate the coverage rates of all CIs. The coverage rates and average widths are presented in Table \ref{table:ci_cover_random_walk}. 

From Table \ref{table:ci_cover_random_walk}, the loose CIs perform conservatively with near to $99\%$ coverage rates in all cases. On the other hand, the tight CIs constructed from our probabilistically efficient estimators have coverage rates below $95\%$ in most of the cases. Moreover, as $\sigma$ decreases (the probability become rarer), the coverage rates first drop from around 0.93 (with $\sigma=0.3$) to around 0.89 (with $\sigma=0.2$), then they improve as $\sigma$ further decreases and reaches around $95\%$ when $\sigma=0.1$. The tight CIs by asymptotically efficient estimators have more stable coverage rates than the CIs by probabilistically efficient estimators, but also suffer under-coverage in several cases (e.g., 0.86 with $\sigma=0.28$). We also observe that the tight CIs by the probabilistically efficient estimators have better average widths than the CIs by the asymptotically efficient estimators with smaller $\sigma$ (e.g. $\sigma=0.1,0.12$). The results show the validity of the CIs with probabilistically efficient estimators as $\sigma \to 0$, but also that the coverage rate may not always monotonically improve as the problem becomes rarer. 

\begin{table}[]
\caption{Higly accurate point estimates (and 95\% CI) using asymptotically efficient estimators for the problem in Section~\ref{sec:numerical_sum}. The estimates are computed with $10^7$ samples.} \label{table:accurate_random_walk}
\centering
\begin{tabular}{|l|l|}
\hline
$\sigma$ & $\hat{p}$\\ \hline
0.1 & $ (5.57\pm 0.039)  \times 10^{-26}$\\ \hline
0.12 & $ (1.30\pm 0.008)  \times 10^{-18}$\\ \hline
0.14 & $ (3.82\pm 0.020)  \times 10^{-14}$\\ \hline
0.16 & $ (3.18\pm 0.015)  \times 10^{-11}$\\ \hline
0.18 & $ (3.32\pm 0.014)  \times 10^{-9}$\\ \hline
0.2 & $9.51 (\pm 0.035)  \times 10^{-8}$\\ \hline
0.22 &  $7.93 (\pm 0.004) \times 10^{-6}$\\ \hline
0.24 & $1.55 (\pm 0.010)\times 10^{-5}$\\ \hline
0.26 & $3.59 (\pm 0.024) \times 10^{-6}$ \\ \hline
0.28 & $1.20(\pm 0.003) \times 10^{-4}$ \\ \hline
0.3 &  $3.22(\pm 0.009) \times 10^{-4}$\\ \hline
\end{tabular}
\end{table}

\begin{table}[]
\caption{Coverage rates and average CI widths of the loose confidence intervals (``Loose CI by PE''), the tight confidence intervals (``Tight CI by PE'') for probabilistically efficient estimators, and the tight confidence intervals for asymptotically efficient estimators (''Tight CI by AE'') in the experiments of Section~\ref{sec:numerical_sum}.}\label{table:ci_cover_random_walk}
\centering
\begin{tabular}{|l|l|l||l|l||l|l|}
\hline
      &    \multicolumn{2}{l||}{Loose CI by PE}     &   \multicolumn{2}{l||}{Tight CI by PE}    &   \multicolumn{2}{l|}{Tight CI by AE}     \\ \hline
$\sigma$ & Coverage & Width &  Coverage & Width  & Coverage  & Width \\ \hline
0.1 & 0.9997 & $ 1.49 \times 10^{-26}$ & 0.949 & $ 7.88 \times 10^{-27}$  & 0.926 & $ 2.41 \times 10^{-26}$ \\ \hline
0.12 & 0.999 & $ 3.58 \times 10^{-19}$ & 0.938 & $ 1.99 \times 10^{-19}$  & 0.917 & $ 5.08 \times 10^{-19}$ \\ \hline
0.14 & 0.997 & $ 1.07 \times 10^{-14}$ & 0.914 & $ 6.15 \times 10^{-15}$  &0.951  & $ 1.32\times 10^{-14}$ \\ \hline
0.16 & 0.994 & $ 9.37 \times 10^{-12}$ & 0.897 & $5.54 \times 10^{-12}$  & 0.908 & $ 9.67 \times 10^{-12}$ \\ \hline
0.18 & 0.990 & $ 9.89 \times 10^{-10}$ & 0.892 & $5.95 \times 10^{-10}$  & 0.951 & $ 8.72 \times 10^{-10}$ \\ \hline
0.2 & 0.988 & $ 2.86 \times 10^{-8}$ & 0.895 & $ 1.75 \times 10^{-8}$  & 0.964 & $ 2.24 \times 10^{-8}$ \\ \hline
0.22 & 0.988 & $ 3.15 \times 10^{-7}$ & 0.901 &$ 1.93 \times 10^{-7}$   & 0.936 & $ 2.45 \times 10^{-7}$ \\ \hline
0.24 & 0.990 & $ 2.20 \times 10^{-6}$ & 0.913 &$ 1.35 \times 10^{-6}$   & 0.959 & $ 1.55 \times 10^{-6}$ \\ \hline
0.26 & 0.991& $9.67 \times 10^{-6}$ & 0.918  &$ 6.00 \times 10^{-6}$  & 0.932 & $ 6.49 \times 10^{-6}$ \\ \hline
0.28 & 0.992& $ 2.99\times 10^{-5}$ & 0.924 & $1.86 \times 10^{-5}$  & 0.861 & $ 2.01 \times 10^{-5}$ \\ \hline
0.3 & 0.993 & $7.71 \times 10^{-5}$  & 0.927 &$4.80 \times 10^{-5}$   & 0.956 & $ 4.90 \times 10^{-5}$ \\ \hline
\end{tabular}
\end{table}

\subsubsection{MNIST Example.}
For the experiment in Section~\ref{sec:numerical_mnist}, we use $10^{10}$ runs of crude Monte Carlo to approximate the ground truths, which are shown in Table~\ref{table:accurate_mnist}. Note that the estimate for $\sigma=0.17$ is relatively less accurate than other estimates, revealed by the CI width in the magnitude of around 0.1 of the probability estimate. We obtain our probabilistically efficient estimator by generating $10^4$ independent samples and construct CIs \eqref{eq:loose_CI_PE} and \eqref{eq:tight_CI_PE} based on this estimator. Since locating all dominating points to construct asymptotically efficient estimator is computationally infeasible in this example, we use IS estimators that mix the most significant 100 dominating points (the number of dominating points obtained from our sequential mixed integer programming procedure in Algorithm \ref{algo:gaussian}) as a proxy. We construct CI \eqref{eq:tight_CI_PE} from this estimator using $10^4$ samples. We use $10^5$ experimental repetitions to estimate the coverage rates and average widths of the CIs from probabilistically efficient estimators and $10^3$ repetitions for the CIs from IS estimators using 100 dominating points (we use repetition size $10^3$ instead of $10^5$ because of the long computational time caused by a large number of mixtures in the IS distribution). The results are presented in Table~\ref{table:ci_cover_mnist}. 

From Table~\ref{table:ci_cover_mnist}, we observe that in three out of the four cases, the tight CIs constructed from probabilistically efficient estimators provide coverage rates that are slightly below $95\%$. Similar to the previous random walk overshoot problem, the coverage rates are closer to $95\%$ for rarer problems (e.g., the coverage is $94.9\%$ for $\sigma=0.17$). The under-coverage is alleviated when we use more than one dominating point, as shown in the row of ``Tight CI by AE'' (where we use 100 dominating points). On the other hand, the loose CIs have higher than nominal coverage rates in all cases, but are conservative since the rates are around $97.5\%-99.5\%$. Again, we observe the validity of the CIs with probabilistically efficient estimators as the rare-event probability decreases, which validates our analysis.

\begin{table}[]
\caption{Highly accurate point estimates (and 95\% CI) using asymptotically efficient estimators for the problem in Section~\ref{sec:numerical_mnist}. The estimates are computed with $10^{10}$ samples.} \label{table:accurate_mnist}
\centering
\begin{tabular}{|l|l|l|l|l|}
\hline
$\sigma$           & 0.17     & 0.18     & 0.19     & 0.2      \\ \hline
$\hat{p}$ &  $3.15 (\pm 0.15) \times 10^{-7}$ &  $1.29 (\pm 0.031) \times 10^{-6}$ & $4.22 (\pm 0.056) \times 10^{-6}$ &  $1.17 (\pm 0.0094) \times 10^{-5}$ \\ \hline
\end{tabular}
\end{table}

\begin{table}[]
\caption{Coverage rates and average CI widths of the loose confidence intervals (``Loose CI by PE''), the tight confidence intervals (``Tight CI by PE'') for probabilistically efficient estimators, and the tight confidence intervals for asymptotically efficient estimators (''Tight CI by AE'') in the experiments of Section~\ref{sec:numerical_mnist}.}\label{table:ci_cover_mnist}
\centering
\begin{tabular}{|l|l|l|l|l|l|}
\hline
& $\sigma$             & 0.17     & 0.18     & 0.19     & 0.2      \\ \hline
\multirow{2}{*}{Loose CI by PE  } & Coverage Rate           & 0.996   & 0.978   & 0.980   & 0.977   \\ \cline{2-6}
&Average Width  & $6.32\times 10^{-8}$ & $2.61 \times 10^{-7}$  & $1.13 \times 10^{-6}$ & $3.20 \times 10^{-6}$ \\ \hline \hline
\multirow{2}{*}{Tight CI by PE   } & Coverage Rate        & 0.949   & 0.874   & 0.885   & 0.877   \\ \cline{2-6}
&Average Width  & $4.19\times 10^{-8}$ & $1.73 \times 10^{-7}$  & $7.50 \times 10^{-7}$ & $2.12 \times 10^{-6}$ \\ \hline \hline
\multirow{2}{*}{Tight CI by AE} & Coverage Rate & 0.958 & 0.933 & 0.945 & 0.951  \\ \cline{2-6}
&Average Width  & $4.77\times 10^{-8}$ & $1.96 \times 10^{-7}$  & $6.13 \times 10^{-7}$ & $1.90 \times 10^{-6}$ \\ \hline
\end{tabular}
\end{table}

\subsubsection{Summary of Experimental Observations on Confidence Interval Construction. } From the CI construction for the three examples in Section~\ref{sec:numerical_glasserman} investigated above, we draw several conclusions: 1) Tight CIs by probabilistically efficient estimators appear to have close to the nominal coverage rate when the problem is rare enough; 2) Loose CIs by probabilistically efficient estimators tend to over-cover, and also have correspondingly larger widths than other methods; 3) Tight CIs by asymptotically efficient estimators appear to give more accurate coverage rates for a larger range of rarity levels than the tight CIs by probabilistically efficient estimators, even though they could still under-cover in some cases; 4) When tight CIs by probabilistically efficient and asymptotically efficient estimators both have accurate coverage rates, their widths appear to be comparable. Overall, it appears that tight CIs by asymptotically efficient estimators are more robust with respect to the rarity level of the problem, which is also in line with the comparison between Theorems~\ref{thm:general_interval2_GE} and \ref{thm:tight_interval_full_IS} (recall the discussion right after Theorem~\ref{thm:tight_interval_full_IS}). Nonetheless, recall that one motivation of us proposing the notion of probabilistic efficiency is that asymptotically efficient estimators, which require using more dominating points in their mixtures, could be computationally challenging to construct.

\section{Future Work}\label{sec:conclusion}
We conclude the paper with further discussions on the potential risk of the current framework and some future directions.

\subsection{Developing Diagnosis for Finite-Sample Under-Estimation}

Similar to the established notion of asymptotic efficiency in the rare-event simulation literature, our probabilistic efficiency framework is asymptotic. For a fixed rarity parameter $\gamma$ and given simulation size $n$ in practice, more work needs to be investigated to judge whether the obtained estimate is reliable or not. In particular, a risk of missing dominating points is finite-sample under-estimation. To be more specific, we look at the following example.

\begin{example} Suppose that our goal is to estimate $p=P(\frac{1}{\gamma}X_{\gamma}\in(-\infty,-k]\cup[1,\infty))$ where $X_{\gamma}\sim N(0,\gamma)$ and $k>1$. In this case, the most significant dominating point is 1, so we consider choosing $X_{\gamma}\sim N(\gamma,\gamma)$ as the IS distribution. That is, the IS estimator is $Z=I(\frac{1}{\gamma}X_{\gamma}\in(-\infty,-k]\cup[1,\infty))e^{- X_{\gamma}+\gamma/2}$ with $X_{\gamma}\sim N(\gamma,\gamma)$. We generate independent samples $Z^{(1)},\dots,Z^{(n)}$ and use $\hp=\frac{1}{n}\sum_{i=1}^n Z_i$ to estimate $p$. When $1<k<3$, the relative error of $Z$ grows exponentially in $\gamma$ (or $-\log p$), so $Z$ is not asymptotically efficient (see the proof of Proposition \ref{lem:example}). On the other hand, from our previous derivations, $Z$ is strongly probabilistically efficient, so we can still get a reliable estimate when $\gamma$ is sufficiently large. However, in the finite-sample case where $k$ is close to 1 and $\gamma$ is not large, if $n$ is chosen as a moderate size, then it could happen that $\{X\leq -k\gamma\}$ is not hit and hence we would get an estimate close to $p_1=P(\frac{1}{\gamma}X_{\gamma}\geq 1)$, but $p_2=P(\frac{1}{\gamma}X_{\gamma}\leq -k)$ is not negligible compared to $p_1$. In other words, we \emph{under-estimate} $p$. As a specific example, let $k=1.01$ and $\gamma=16$. Then $\tp_2=\tP(\frac{1}{16}X_{16}\leq -1.01)=\bar{\Phi}(8.04)\approx 4.44\times 10^{-16}$. In this case, for a moderate $n$ like 1000, with probability almost 1 the set $\{\frac{1}{16}X_{16}\leq -1.01\}$ is never hit. However, $p_2/p_1=\bar{\Phi}(4.04)/\bar{\Phi}(4)\approx 0.84$. This means with probability almost 1, we under-estimate $p$ by about $0.84/(1+0.84)\approx 46\%$.
\label{example:under_estimation}
\end{example}

Example \ref{example:under_estimation} shows that under finite parameter value and finite sample, it could be hard to tell whether we can safely drop less significant dominating points. A good aspect about the conclusion in this example, however, is that the under-estimation is arguably acceptable in relative term (i.e., the estimate is still in the same magnitude as the ground truth), pointing to an estimation resembling weak probabilistic efficiency. We should emphasize that a similar concern applies to asymptotically efficient estimators as well. That is, it is difficult to guarantee whether the sample size $n$ is large enough to give a reliable estimate for a given setup and rarity parameter value. Nevertheless, increasing $n$ would improve the performance of asymptotically efficient estimators, but for probabilistic efficiency, we do not have the luxury of increasing $n$ since our framework requires $n$ to be moderate in size. This points to more need of developing diagnostic methods to detect under-estimation due to finite-sample effects in the future.

\subsection{Further Developing Theory of IS with Missed Dominating Points}

Despite the presence of under-estimation risks described above, we maintain our motivation of probabilistic efficiency as  a theory to allow one to use few dominating points in problems where finding all of them is infeasible. In fact, what we have focused on in this paper is only one theory where dropping dominating points is valid, among other possibilities. 
To support this, we revisit Example \ref{example:under_estimation}:
\begin{example}[Example \ref{example:under_estimation} continued]
Consider the problem setting in Example \ref{example:under_estimation}. When $k\geq 3$, the relative error of $Z$ grows only polynomially in $\gamma$ (or $-\log p$), and hence $Z$ is asymptotically efficient (see the proof of Proposition \ref{lem:example}). The intuitive explanation is that $\tp_2=\tP(\frac{1}{\gamma}X_{\gamma}\leq -k)$ is extremely small which mitigates the blow-up of the likelihood ratio. More rigorously, when $\frac{1}{\gamma}X_{\gamma}\leq -k$ we have $Z\geq e^{(k+\frac12)\gamma}$, but $\tP(\frac{1}{\gamma}X_{\gamma}\leq -k)=\bar{\Phi}((k+1)\sqrt{\gamma})=\Theta(\frac{1}{\sqrt{\gamma}} e^{-\frac{(k+1)^2}{2}\gamma})$. Overall $\tE(Z^2)=\Theta(\frac{1}{\sqrt{\gamma}}e^{-\gamma})$, and hence $\tE(Z^2)/p^2=\Theta(\sqrt{\gamma})$ which grows polynomially in $-\log p=\Theta(\gamma)$.
\label{example:rare}
\end{example}

Example \ref{example:rare} shows that sometimes the missed dominating points are so rare that missing them does not even harm the asymptotic efficiency. Here, explaining the validity of IS with missed dominating points does not require probabilistic efficiency, but instead an alternate analysis of asymptotic efficiency that is tighter than the standard approach in the literature. Nonetheless, this phenomenon does not apply to our proposed estimators in Sections \ref{sec:prelim numerics} and \ref{sec:additional numerics}, as we have mathematically verified the asymptotic inefficiency in our considered estimator in each example.

\begin{figure}
    \centering
    \includegraphics[width=0.8\textwidth]{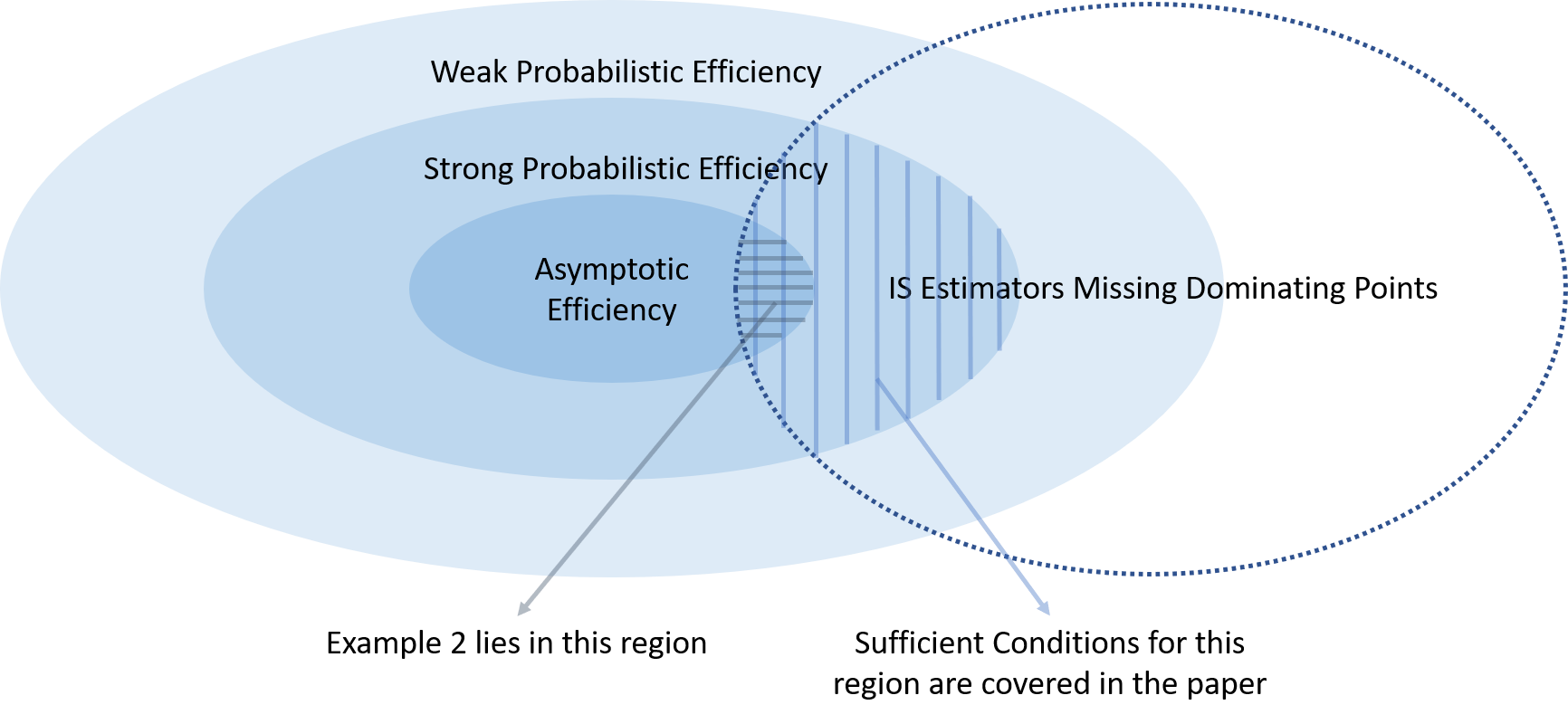}
    \caption{Relationships among different efficiency criteria and IS estimators missing dominating points.}
    \label{fig:diagram}
\end{figure}

More generally, we conclude our paper with Figure \ref{fig:diagram}, which shows the relations among different efficiency criteria and IS with missed dominating points. In this paper, we have built sufficient conditions to achieve strong probabilistic efficiency with missed dominating points. Example \ref{example:rare}, on the other hand, gives an example in achieving asymptotic efficiency with missed dominating points. Our immediate future endeavor is to fill in the regions in Figure \ref{fig:diagram} that are not covered by the current work, including the scenario depicted by Example \ref{example:rare}, relaxing the current sufficient conditions of probabilistic efficiency to allow $\tp_2$ to be less tiny, achieving weak instead of strong probabilistic efficiency with missed dominating points, and moreover, to understand the ``complementary" regions where missing dominating points would be guaranteed to violate the efficiency notions. In summary, this work serves as a first step in a new line of analysis aiming to relax existing variance-based efficiency criteria in rare-event simulation to be applicable to larger-scale and more complex problems.

\ACKNOWLEDGMENT{We gratefully acknowledge support from the National Science Foundation under grants CAREER CMMI-1834710, IIS-1849280 and IIS-1849304.}

\bibliographystyle{informs2014} 
\bibliography{citation} 
\newpage

\begin{APPENDICES}

\section{Algorithms}
\label{app:algorithm}
Algorithm \ref{algo:dominating_points} shows a procedure to obtain dominating sets. We briefly explain the idea here. First we minimize $I(x)$ over $x\in\cE$ and get the first dominating point $a_1$. If $\cE\subset\{x\in\R^d:s_{a_1}^\top(x-a_1)\geq 0\}$, then we know that $\{a_1\}$ is a dominating set and hence we could stop. Otherwise, we minimize $I(x)$ over $x\in \cE\setminus\{x\in\R^d:s_{a_1}^\top(x-a_1)\geq 0\}$ to get the second dominating point $a_2$. Then we check whether $\cE\in\bigcup_{i=1}^2\{x\in\R^d:s_{a_i}^\top(x-a_i)\geq 0\}$. By repeating this process, we would finally get a dominating set $A=\{a_1,\dots,a_r\}$ with $I(a_1)\leq\dots\leq I(a_r)$.
\begin{algorithm}[htbp]
\KwIn{Rarity parameter $\gamma$, rare-event set $\cE\subset\R^d$, rate function $I(y)$, function $a\mapsto s_a$.}
\KwOut{Dominating set $A$.}
\nl Start with $A = \emptyset,r=0$;\\

\nl {\bf While } $\{x\in\R^d: x\in\cE , s_{a_i}^\top(x-a_i) <0, \mbox{ $\forall i=1,\dots,r$} \} \neq  \emptyset$ {\bf do }\\ 
\nl \ \ \ \ \ \ Find a dominating point $a_{r+1}$ by solving the optimization problem \begin{align*} 
    a_{r+1} =\arg \min_{x} &\ \  I(x) \ \ \ \\
\text{s.t.}\ \ \  &x\in\cE,\ \  \\ &s_{a_i}^\top(x-a_i) <0, \ \mbox{$\forall i=1,\dots,r$}
\end{align*}

\ \ \ \ \ \ and update $A \leftarrow A \cup \{a_{r+1}\},r\leftarrow r+1$.\\
\nl {\bf End}\\

\caption{Sequentially find all the dominating points.} 
\label{algo:dominating_points}
\end{algorithm}

\section{A New Alternative Asymptotic Regime and Theoretical Guarantees}
\label{sec:guarantees}
In this section, we investigate a new regime that $\cA_{\gamma}=\{X\in\cE\}$ and $\cE=\{x\in\R^d:g(x)\geq\gamma\}$ where $X\in\R^d$ and $g:\R^d\to \R$ is a function. We propose this setting as it arises as a generic representation of recent problems in intelligent system safety testing \citep{bai2022rare}. There, $g$ could be highly complicated and leads to a gigantic number of dominating points, which in turn motivates the consideration of dropping most of them and our notion of probabilistic efficiency. We note that technically this setting is slightly different from the classical Gartner-Ellis regime in terms of the position of the scaling parameter $\gamma$ (see, e.g. \citealt{sadowsky1990large}), but conceptually similar.

In Section \ref{sec:guarantees_general}, we adapt the notions of rate function and dominating points to this new regime. Then we state the assumptions under which we could build on Theorem \ref{thm:strong_probabilistic_efficiency} to obtain reliable point estimates and CIs as in Section~\ref{sec:guarantees_GE}. We note that this new regime is harder to analyze as the rare-event set $\cE$ can change in a complicated way with $\gamma$, and thus inevitably our assumptions are relatively restrictive and need to be verified case by case. In Section \ref{sec:guarantees_gaussian}, we consider the special (but important) case where $X$ follows a Gaussian distribution and $g$ is piecewise linear, in particular propose a simple stopping strategy to determine whether it is safe to stop searching for the remaining dominating points. Throughout this section, we write $\alpha(\gamma)\sim\beta(\gamma)$ if $\alpha(\gamma)/\beta(\gamma)$ 
is subexponential in $-\log p(\gamma)$.

\subsection{Guarantees for General Input Distribution}
\label{sec:guarantees_general}
We define $\mu(x)=\log Ee^{x^\top X},x\in\R^d$ as the cumulant generating function of $X$ and $I(y)=\sup_{x\in\R^d}\{x^\top y-\mu(x)\},y\in\R^d$ as its Legendre transform. For any set $\cE\subset\R^d$, we denote $I(\cE)=\inf_{y\in\cE}I(y)$. Parallel to Assumptions~\ref{asm:mu_GE} and~\ref{asm:E_GE}, we make the following two assumptions.
\begin{assumption}
$\mu(x)$ satisfies the following conditions:
\begin{enumerate}
    \item $0\in \mathcal{D}(\mu)^{\circ}$;
    \item $\mu$ is essentially smooth, i.e., $\mathcal{D}(\mu)^{\circ}$ is non-empty, $\mu$ is differentiable everywhere in $\mathcal{D}(\mu)^{\circ}$ and $\mu$ is steep.
\end{enumerate}
\label{asm:mu}
\end{assumption}

\begin{assumption}
For any $\gamma$, $\cE=\cE(\gamma)\subset \R^d$ is a Borel set such that $\overline{\cE}=\overline{\cE^{\circ}}, \cE^\circ\cap\mathcal{D}(I)^\circ\ne\emptyset$ and $I(\cE)>0$.
\label{asm:E}
\end{assumption}

The concepts of dominating set and dominating points are also similar, but now they change with $\gamma$:
\begin{definition}[Dominating Set (New Regime)]
Suppose that Assumptions \ref{asm:mu} and \ref{asm:E} hold. We call $A=A(\gamma)\subset\overline{\cE}$ a dominating set for $\cE$ if 
\begin{enumerate}
    \item For each $a\in A$, there exists a unique $s_a\in\R^d$ such that $\nabla\mu(s_a)=a$;
    \item $\cE\subset\bigcup_{a\in A}\{x\in\R^d:s_a^\top (x-a)\geq 0\}$;
    \item $A\setminus\{a\}$ does not satisfy the first two conditions for any $a\in A$.
\end{enumerate}
We call any point in $A$ a dominating point. For two dominating points $a$ and $a'$, we say $a$ is more significant than $a'$ if $I(a)<I(a')$.
\label{def:dominating_set}
\end{definition}

Suppose that $A=\{a_1,\dots,a_r\}$ is a dominating set. Then the corresponding mixture IS distribution is given by $\frac{d\tP}{dP}(\omega)=\sum_{i=1}^r\alpha_i e^{s_{a_i}^\top X-\mu(s_{a_i})}$ with $\sum_{i=1}^r\alpha_i=1,\alpha_i>0,\forall i$. Like the discussion in Section \ref{sec:background}, we could split the rare-event set as $\cE=\bigcup_{i=1}^r\cE_i$ where $\cE_i$'s are disjoint and $a_i\in\overline{\cE_i}\subset\{x\in\R^d:s_{a_i}^\top (x-a_i)\geq 0\}$. Note that $s_{a_i}^\top (x-a_i)\geq 0$ still implies that $I(x)\geq I(a_i)$. Now, the likelihood ratio is given by
\begin{equation*}
    L(\omega):=\frac{dP}{d\tP}(\omega)=\frac{1}{\sum_{i=1}^r\alpha_i e^{s_{a_i}^\top X-\mu(s_{a_i})}}
\end{equation*}
and it satisfies that for any $i$,
\begin{equation*}
    L(\omega)\leq \frac{1}{\alpha_i e^{s_{a_i}^\top X-\mu(s_{a_i})}}=\frac{1}{\alpha_i}e^{-s_{a_i}^\top (X-a_i)-I(a_i)}.
\end{equation*}
Hence we have that $I_{\cE_i}(X)L(\omega)\leq \frac{1}{\alpha_i}e^{-I(a_i)}$. Then
\begin{equation}
    \tVar(I(X\in\cE)L(\omega))\leq \tE(I(X\in\cE)L^2(\omega))\leq \frac{1}{(\min_i\alpha_i)^2}e^{-2\min_iI(a_i)}.
    \label{eqn:vartilde}
\end{equation}
Thus, supposing we have a large deviations asymptotic given by $p=P(X\in\cE)\sim e^{-\min_iI(a_i)}$, then  combining with \eqref{eqn:vartilde} will give us that the IS estimator is asymptotically efficient.

Now we consider using a partial list of dominating points. In particular, we can still sequentially fill in the dominating set $A=\{a_1,\dots,a_r\}$ where $a_i=\arg\min\{I(y):y\in\cE,s_{a_j}^\top (y-a_j)<0,j=1,\dots,i-1\}$ and hence $I(a_1)\leq \dots\leq I(a_r)$, and suppose that we have a stopping strategy $k=k(\gamma)$ with $1\leq k\leq r$ before locating all the dominating points. We choose the mixture IS distribution given by
\begin{equation}
    \frac{d\tP}{dP}(\omega)=\sum_{i=1}^k\frac{1}{k} e^{s_{a_i}^\top X-\mu(s_{a_i})}.
    \label{eqn:IS_distribution}
\end{equation}
Unlike in \eqref{eqn:IS_distribution_GE}, we no longer assign an individual weight $\alpha_i$ to each dominating point $a_i$. Instead, we only consider the uniform mixture for simplicity, since now the number of dominating points $r$ and $k$ can both potentially change with $\gamma$. 

First of all, we summarize the basic assumptions on the problem setting to ensure that there exists a dominating set with moderate size:
\begin{assumption}
Consider the problem of estimating $p=P(X\in\cE)$ with $\cE=\{x\in\R^d:g(x)\geq\gamma\}$ where $X\in\R^d$ is a random vector and $g$ is a function. Assume that 
\begin{enumerate}
    \item $p\to 0$ as $\gamma\to\infty$;
    \item Assumption \ref{asm:mu} holds for the cumulant generating function of $X$ under $P$;
    \item Assumption~\ref{asm:E} holds for $\cE$;
    \item For any $\gamma$, there exists a dominating set for $\cE$, denoted as $A=\{a_1,\dots,a_r\}$ where $a_i=\arg\min\{I(y):y\in\cE,s_{a_j}^\top (y-a_j)<0,j=1,\dots,i-1\}$ and hence $I(a_1)\leq\dots\leq I(a_r)$. Besides, $r=r(\gamma)$ grows at most subexponentially in $-\log p$.
\end{enumerate}
\label{asm:problem}
\end{assumption}

Now consider a stopping strategy $k=k(\gamma)$ and the corresponding IS distribution \eqref{eqn:IS_distribution}. Under Assumption \ref{asm:problem}, for $k=k(\gamma)$ with $1\leq k\leq r$, denote $\cE_1=\cE\cap\bigcup_{i=1}^k\{x\in\R^d:s_{a_i}^\top (x-a_i)\geq 0\}$ and $\cE_2=\cE\setminus\cE_1$. Corresponding to the settings in Theorem \ref{thm:strong_probabilistic_efficiency}, we let $\cA_{\gamma}^j=\{X\in\cE_j\},j=1,2$. We introduce Assumptions \ref{asm:p2_negligible} and \ref{asm:p2tilde}.
\begin{assumption}
If $k=r$, then let $a_{k+1}=\infty\mathbf{1}_d$. Assume that $p_1:=P(\cA_{\gamma}^1)=P(X\in\cE_1)\sim e^{-I(a_1)}$ and that $p_2:=P(\cA_{\gamma}^2)=P(X\in\cE_2)$ is upper bounded by $e^{-I(a_{k+1})}$ up to subexponential factor in $-\log p$. Besides, $e^{I(a_1)-I(a_{k+1})}$ exponentially decays in $-\log p$.
\label{asm:p2_negligible}
\end{assumption}
\begin{assumption}
Assume that $\tp_2:=\tP(\cA_{\gamma}^2)=\tP(X\in\cE_2)$ exponentially decays in $-\log p$.
\label{asm:p2tilde}
\end{assumption}

Roughly, Assumption \ref{asm:p2_negligible} implies that $p_2$ is exponentially smaller than $p_1$ and Assumption \ref{asm:p2tilde} implies that $\cE_2$ is hardly hit even under $\tP$. Applying Theorem \ref{thm:strong_probabilistic_efficiency}, we have that 
\begin{theorem}[Attaining probabilistic efficiency with a partial list of dominating points (new regime)]
Under Assumptions \ref{asm:problem}, \ref{asm:p2_negligible} and \ref{asm:p2tilde}, the IS estimator $Z=I(X\in\cE)\frac{dP}{d\tP}(\omega)$ under $\tP$ given by \eqref{eqn:IS_distribution} is strongly probabilistically efficient.
\label{thm:general_point_estimate}
\end{theorem}

Similar to Section~\ref{sec:guarantees_GE}, we can also construct asymptotically valid CIs with the sample mean and sample variance.
\begin{theorem}[Constructing confidence intervals with probabilistically efficient estimators (new regime)]
Assume that Assumptions \ref{asm:problem}, \ref{asm:p2_negligible} and \ref{asm:p2tilde} hold. The IS estimator is  $Z=I(X\in\cE)\frac{dP}{d\tP}(\omega)$ under $\tP$ given by \eqref{eqn:IS_distribution}. We sample $X^{(1)},\dots,X^{(n)}$ i.i.d. from $\tP$ and let $Z^{(i)}=I(X^{(i)}\in\cE)\frac{dP}{d\tP},Z_1^{(i)}=I(X^{(i)}\in\cE_1)\frac{dP}{d\tP},i=1,\dots,n$. Use $\hat p$ and $\hat{V}$ to respectively denote the sample mean and sample variance of $Z^{(i)}$'s. In this case, If $n$ is subexponentially growing in $-\log p$ as $\gamma\to\infty$. Then, for any $0<\alpha<1$,
\begin{equation*}
\liminf_{\gamma\to\infty}\tilde{P}\left(|\hat p-p|\leq\sqrt{\frac{2\hat{V}\log(4/\alpha)}{n}}+\frac{7\log(4/\alpha)ke^{-I(a_1)}}{3(n-1)}\right)\geq1-\alpha.
\end{equation*}
That is, 
\begin{equation*}
\hat p\pm \left(\sqrt{\frac{2\hat{V}\log(4/\alpha)}{n}}+\frac{7\log(4/\alpha)ke^{-I(a_1)}}{3(n-1)}\right)
\end{equation*}
is an asymptotically valid $(1-\alpha)$-level CI for $p$.
\label{thm:general_interval1}
\end{theorem}

Similar to Section~\ref{sec:guarantees_GE}, the CI in Theorem~\ref{thm:general_interval1} is more conservative than the CLT-based interval $\hp\pm z_{1-\alpha/2}\sqrt{\frac{\hat{V}}{n}}$. If additionally the following assumption holds, then the CLT-based CI is also asymptotically valid:
\begin{assumption}
Denote $Z_1=I(X\in\cE_1)\frac{f(X)}{\tf(X)}$ under $\tP$. Assume that $\tVar(Z_1)\sim e^{-2I(a_1)}$.
\label{asm:BE_error}
\end{assumption}

Similar to Lemma \ref{lem:variance}, Assumption \ref{asm:BE_error} serves to control the Berry-Esseen error bound. To be more concrete, with this additional assumption, we have:
\begin{theorem}[Constructing tight confidence intervals (new regime)]
Assume that Assumptions \ref{asm:problem}, \ref{asm:p2_negligible}, \ref{asm:p2tilde} and \ref{asm:BE_error} hold. The IS estimator is  $Z=I(X\in\cE)\frac{dP}{d\tP}(\omega)$ under $\tP$ given by \eqref{eqn:IS_distribution}. We sample $X^{(1)},\dots,X^{(n)}$ i.i.d. from $\tP$ and let $Z^{(i)}=I(X^{(i)}\in\cE)\frac{dP}{d\tP},Z_1^{(i)}=I(X^{(i)}\in\cE_1)\frac{dP}{d\tP},i=1,\dots,n$. Use $\hat p$ and $\hat{V}$ to respectively denote the sample mean and sample variance of $Z^{(i)}$'s. In this case, we could choose $n$ subexponentially growing in $-\log p$ such that  $\frac{k^2e^{-2I(a_1)}}{n\tVar(Z_1^{(1)})}\to 0$ and $\frac{\tE^2|Z_1^{(1)}-p_1|^3}{n\tVar^3(Z_1^{(1)})}\to 0$ as $\gamma\to\infty$. Then, for any $0<\alpha<1$,
\begin{equation*}
\liminf_{\gamma\to\infty}\tilde{P}\left(|\hat p-p|\leq z_{1-\alpha/2}\sqrt{\frac{\hat{V}}{n}}\right)\geq1-\alpha.
\end{equation*}
That is, 
\begin{equation*}
\hat p\pm z_{1-\alpha/2}\sqrt{\frac{\hat{V}}{n}}
\end{equation*}
is an asymptotically valid $(1-\alpha)$-level CI for $p$.
\label{thm:general_interval2}
\end{theorem}

In the above, Assumptions~\ref{asm:problem}, \ref{asm:p2_negligible}, \ref{asm:p2tilde} and \ref{asm:BE_error} are technical conditions that need to be analyzed case by case. The next Section \ref{sec:guarantees_gaussian} focuses on a specific but important problem setting where we show how to verify all these conditions.

\subsection{Guarantees for Gaussian Input Distribution}
\label{sec:guarantees_gaussian}
Suppose $X\sim N(\lambda,\Sigma)$ under $P$ where $\lambda\in\R^d$ and $\Sigma\in\R^{d\times d}$ is positive definite. In this case $\mu(x)=\lambda^\top x+\frac{1}{2}x^\top \Sigma x$, $s_a=\Sigma^{-1}(a-\lambda)$ and $I(x)=\frac{1}{2}(x-\lambda)^\top \Sigma^{-1}(x-\lambda)$. Moreover, we suppose that $g$ is a piecewise linear function. Related theoretical guarantees for piecewise linear function of Gaussian input can be found in \citet{bai2022rare}. 

We use the following natural stopping strategy. First, fix a constant $C>1$. Then, we sequentially find the dominating points where we stop as long as 
$(a_{k+1}-\lambda)^\top \Sigma^{-1}(a_{k+1}-\lambda)>C(a_k-\lambda)^\top \Sigma^{-1}(a_k-\lambda)$. If there is no such $k$, then we let $k=r$, the number of all dominating points. The IS distribution is chosen as $X\sim\frac{1}{k}\sum_{i=1}^k\phi(x;a_i,\Sigma)$. We summarize these in Algorithm \ref{algo:gaussian}. 

\begin{algorithm}[htbp]
\KwIn{Piecewise linear function $g$, rarity parameter $\gamma$, input distribution $N(\lambda,\Sigma)$, threshold $C>1$, sample size $n$.}
\KwOut{IS estimate $\hp$.}
\nl Start with $k=0$;\\

\nl {\bf While } $\{x: g(x) \geq \gamma , (a_i-\lambda)^\top\Sigma^{-1}(x-a_i) <0, \mbox{ $\forall i=1,\dots,k$} \} \neq  \emptyset$ {\bf do }\\ 
\nl \ \ \ \ \ \ Find a dominating point $a_{k+1}$ by solving the optimization problem \begin{align*} 
    a_{k+1} =\arg \min_{x} &\ \  (x-\lambda)^\top\Sigma^{-1}(x-\lambda) \ \ \ \\
\text{s.t.}\ \ \  &g(x) \geq \gamma,\ \  \\ &(a_i-\lambda)^\top\Sigma^{-1}(x-a_i) <0, \mbox{ $\forall i=1,\dots,k$.}
\end{align*}\\
\nl \ \ \ \ \ \ {\bf If } $k>0$ and $(a_{k+1}-\lambda)^\top\Sigma^{-1}(a_{k+1}-\lambda)> C(a_k-\lambda)^\top\Sigma^{-1}(a_k-\lambda)$ {\bf do}\\
\nl \ \ \ \ \ \ \ \ \ \ \ \ {\bf Break}\\
\nl \ \ \ \ \ \ {\bf Else do}\\
\nl \ \ \ \ \ \ \ \ \ \ \ \ Update $k\leftarrow k+1$;\\
\nl {\bf End}\\
\nl Sample $X_1,\dots,X_n$ from the mixture distribution $\frac{1}{k}\sum_{i=1}^k\phi(x;a_i,\Sigma)$. \\
\nl Compute the IS estimate $\hp=\frac{1}{n}\sum_{i=1}^n I(g(X_i)\geq\gamma)L(X_i)$ where the likelihood ratio function is 
\begin{equation*}
    L(x)=\frac{e^{-\frac12(x-\lambda)^\top\Sigma^{-1}(x-\lambda)}}{\frac1k\sum_{i=1}^ke^{-\frac12(x-a_i)^\top\Sigma^{-1}(x-a_i)}}.
\end{equation*}

\caption{Estimate $P(g(X)\geq\gamma)$ with Gaussian input $X$ and piecewise linear function $g$.} 
\label{algo:gaussian}
\end{algorithm}

We have the following theorem suggesting that under the above setting, all the assumptions listed in Section \ref{sec:guarantees_general} are satisfied.
\begin{theorem}[Verification of assumptions for Gaussian inputs]
Suppose that $X\sim N(\lambda,\Sigma)$ under $P$ where $\lambda\in\R^d$ and $\Sigma\in\R^{d\times d}$ is positive definite. $\cE=\{x\in\R^d:g(x)\geq\gamma\}$ where $g$ is a piecewise linear function. $C>1$ is a constant. Assume that for any $\gamma$, $P(X\in\cE)>0$ and $\lambda\notin\overline{\cE}$. For any $\gamma$, we use the stopping strategy described above (more precisely Algorithm \ref{algo:gaussian}) to sequentially find dominating points $a_1,a_2,\dots,a_k$ in decreasing significance and set up the IS distribution $X\sim\frac{1}{k}\sum_{i=1}^k\phi(x;a_i,\Sigma)$. Denote $\cE_1=\cE\cap\bigcup_{i=1}^k\{x\in\R^d:(a_i-\lambda)^\top \Sigma^{-1}(x-a_i)\geq 0\}$ and $\cE_2=\cE\setminus\cE_1$. Then Assumptions \ref{asm:problem}, \ref{asm:p2_negligible}, \ref{asm:p2tilde} and \ref{asm:BE_error} hold.
\label{thm:gaussian_assumptions}
\end{theorem}

Note that in Theorem \ref{thm:gaussian_assumptions}, all assumptions including Gaussianity and piecewise linear $g$ are all straightforward to verify. With Theorem \ref{thm:gaussian_assumptions}, we thus get the following corollaries from Theorems \ref{thm:general_point_estimate}, \ref{thm:general_interval1} and \ref{thm:general_interval2}:
\begin{corollary}[Probabilistic efficiency for Gaussian inputs]
Under the settings in Theorem \ref{thm:gaussian_assumptions}, the IS estimator $Z=I(X\in\cE)\frac{dP}{d\tP}$ is strongly probabilistically efficient.
\label{cor:gaussian_point_estimate}
\end{corollary}
\begin{corollary}[Confidence intervals for Gaussian inputs]
Under the settings in Theorem \ref{thm:gaussian_assumptions}, we sample $X^{(1)},\dots,X^{(n)}$ i.i.d. from $\tP$ and let $Z^{(i)}=I(X^{(i)}\in\cE)\frac{dP}{d\tP}$. Use $\hp$ and $\hat{V}$ to denote the sample mean and sample variance of $Z^{(i)}$'s. We could choose $n$ subexponentially growing in $-\log p$ such that for any $0<\alpha<1$,
\begin{equation*}
    \liminf_{\gamma\to\infty}\tilde{P}\left(|\hat p-p|\leq\sqrt{\frac{2\hat{V}\log(4/\alpha)}{n}}+\frac{7\log(4/\alpha)ke^{-\frac12(a_1-\lambda)^\top \Sigma^{-1}(a_1-\lambda)}}{3(n-1)}\right)\geq1-\alpha.
\end{equation*}
\label{cor:gaussian_interval1}
\end{corollary}
\begin{corollary}[Tight confidence intervals for Gaussian inputs]
Under the settings in Theorem \ref{thm:gaussian_assumptions}, we sample $X^{(1)},\dots,X^{(n)}$ i.i.d. from $\tP$ and let $Z^{(i)}=I(X^{(i)}\in\cE)\frac{dP}{d\tP}$. Use $\hp$ and $\hat{V}$ to denote the sample mean and sample variance of $Z^{(i)}$'s. We could choose $n$ subexponentially growing in $-\log p$ such that for any $0<\alpha<1$,
\begin{equation*}
    \liminf_{\gamma\to\infty}\tilde{P}\left(|\hat p-p|\leq z_{1-\alpha/2}\sqrt{\frac{\hat{V}}{n}}\right)\geq1-\alpha.
\end{equation*}
\label{cor:gaussian_interval2}
\end{corollary}

While the choice of $C$ (as long as it is $>1$) does not affect the guarantee on strong probabilistic efficiency or the asymptotic validity of the CIs, it does affect the accuracy of the estimates for a given, finite $\gamma$. In particular, it is often the case that in practice we only need to solve one problem for a fixed $\gamma$ instead of solving a series of problems with varying $\gamma$. In this scenario, as long as $(a_2-\lambda)^\top \Sigma^{-1}(a_2-\lambda)>(a_1-\lambda)^\top \Sigma^{-1}(a_1-\lambda)$, we can pick a $C>1$ such that we would stop our dominating point search at $k=1$, i.e., use only the first point. Note that this does not imply that for this choice of $C$ we have $k=1$ for any $\gamma$. Thus, the finiteness of the rarity parameter comes into play in a subtle way. 





\section{Proofs}
\proof{Proof of Theorem \ref{thm:global_minimizer}.} First, we know that $I(a_i)=s_{a_i}^\top a_i-\mu(s_{a_i})$ for $i=1,\dots,r$. For any $x\in\R^d$ such that $s_{a_i}^\top(x-a_i)\geq 0$, we have $I(x)\geq s_{a_i}^\top x-\mu(s_{a_i})\geq s_{a_i}^\top a_i-\mu(s_{a_i})=I(a_i)$. By Definition \ref{def:dominating_set_GE}, we have $\cE\subset\bigcup_{i=1}^r\{x\in\R^d:s_{a_i}^\top(x-a_i)\geq 0\}$. Thus, for any $x\in\cE$, there exists $i$ such that $s_{a_i}^\top(x-a_i)\geq 0$, and hence $I(x)\geq I(a_i)$. Therefore, $I(\cE)=\inf_{x\in\cE}I(x)\geq \min_{i=1,\dots,r}I(a_i)$. On the other hand, $a_1,\dots,a_r\in\partial\cE\cap\mathcal{D}(I)^\circ$ and $I$ is differentiable in $\mathcal{D}(I)^\circ$, so $I(\cE)\leq \min_{i=1,\dots,r}I(a_i)$. Overall we have $I(\cE)=\min_{i=1,\dots,r}I(a_i)$.
\hfill\Halmos\endproof

\proof{Proof of Proposition \ref{classical AE}.}
We could split the rare-event set as $\cE=\bigcup_{i=1}^r\cE_i$ where $\cE_i$'s are disjoint and $a_i\in\overline{\cE_i}\subset\{x\in\R^d:s_{a_i}^\top (x-a_i)\geq 0\}$. Note that the hyperplane $\{x\in\R^d:s_{a_i}^\top (x-a_i)=0\}$ is tangent to the rate function level set $\{x\in\R^d:I(x)=I(a_i)\}$ at $a_i$, so $s_{a_i}^\top (x-a_i)\geq 0$ implies that $I(x)\geq I(a_i)$. Now, the likelihood ratio is given by

\begin{equation*}
    L(X_\gamma):=\frac{dP}{d\tP}(X_\gamma)=\frac{1}{\sum_{i=1}^r\alpha_i e^{s_{a_i}^\top X_{\gamma}-\gamma\mu_{\gamma}(s_{a_i})}}
\end{equation*}
and it satisfies that for any $i$,
\begin{equation*}
    L(X_\gamma)\leq \frac{1}{\alpha_i e^{s_{a_i}^\top X_{\gamma}-\gamma\mu_{\gamma}(s_{a_i})}}=\frac{1}{\alpha_i}e^{-s_{a_i}^\top (X_{\gamma}-\gamma a_i)-\gamma (s_{a_i}^\top a_i-\mu_{\gamma}(s_{a_i}))}.
\end{equation*}
Hence we have that 
\begin{equation}
    I\left(\frac{1}{\gamma}X_{\gamma}\in\cE_i\right)L(X_\gamma)\leq \frac{1}{\alpha_i}e^{-\gamma (s_{a_i}^\top a_i-\mu_{\gamma}(s_{a_i}))}
    \label{eqn:upper_bound_GE}
\end{equation}
and 
\begin{equation*}
    I\left(\frac{1}{\gamma}X_{\gamma}\in\cE\right)L(X_\gamma)\leq \max_{i=1,\dots,r}\left\{\frac{1}{\alpha_i}e^{-\gamma (s_{a_i}^\top a_i-\mu_{\gamma}(s_{a_i}))}\right\}.
\end{equation*}
Thus, using Theorem \ref{thm:GE}, we have that 
\begin{align*}
\liminf_{\gamma\to\infty}\frac{\log\tE\left(I\left(\frac{1}{\gamma}X_{\gamma}\in\cE\right)L^2(X_\gamma)\right)}{\log\tE\left(I\left(\frac{1}{\gamma}X_{\gamma}\in\cE\right)L(X_\gamma)\right)}&\geq \liminf_{\gamma\to\infty}\frac{2\log\left(\max_{i=1,\dots,r}\left\{\frac{1}{\alpha_i}e^{-\gamma (s_{a_i}^\top a_i-\mu_{\gamma}(s_{a_i}))}\right\}\right)}{\log P\left(\frac{1}{\gamma}X_{\gamma}\in\cE\right)}\\
&=\liminf_{\gamma\to\infty}\frac{2\max_{i=1,\dots,r}\left\{\frac{1}{\gamma}\log\left(\frac{1}{\alpha_i}e^{-\gamma (s_{a_i}^\top a_i-\mu_{\gamma}(s_{a_i}))}\right)\right\}}{\frac{1}{\gamma}\log P\left(\frac{1}{\gamma}X_{\gamma}\in\cE\right)}\\
&= \frac{-2\min_{i=1,\dots,r}I(a_i)}{-I(\cE)}=2,
\end{align*}
which verifies that the IS estimator is asymptotically efficient.\hfill\Halmos
\endproof

\proof{Proof of Proposition \ref{lem:example}.}
We derive the growth rate of the relative error for a more general case. Suppose that we estimate $p=P(\frac{1}{\gamma}X_{\gamma}\in(-\infty,-k]\cup[1,\infty))$ where $k>1$, $X_\gamma\sim N(0,\gamma)$ under $P$ and the IS distribution is chosen as $X_{\gamma}\sim N(\gamma,\gamma)$ under $\tP$. This is the setup of Examples \ref{example:under_estimation} and \ref{example:rare}. There are two dominating points, 1 and $-k$, and 1 is the more significant one. Under the IS distribution, the likelihood ratio is $L=e^{-X_{\gamma}+\gamma/2}$ and the IS estimator is $Z=I(\frac{1}{\gamma}X_{\gamma}\in (-\infty,-k]\cup[1,\infty))e^{-X_{\gamma}+\gamma/2}$. Then we have that 
\begin{align*}
\tE(Z^2)&=\tE\left(I\left(\frac{1}{\gamma}X_{\gamma}\in (-\infty,-k]\cup[1,\infty)\right)e^{-2X_{\gamma}+\gamma}\right)\\
&=\int_{x\leq -k\gamma\text{ or }x\geq\gamma}e^{-2x+\gamma}\frac{1}{\sqrt{2\pi\gamma}}e^{-\frac{(x-\gamma)^2}{2\gamma}}dx\\
&=\int_{x\leq -k\gamma\text{ or }x\geq\gamma}e^{\gamma}\frac{1}{\sqrt{2\pi\gamma}}e^{-\frac{(x+\gamma)^2}{2\gamma}}dx\\
&=\int_{y\leq (1-k)\sqrt{\gamma}\text{ or }y\geq 2\sqrt{\gamma}}e^{\gamma}\frac{1}{\sqrt{2\pi}}e^{-\frac{y^2}{2}}dy\ \ \ \ (y=(x+\gamma)/\sqrt{\gamma})\\
&=e^{\gamma}\left(\bar{\Phi}((k-1)\sqrt{\gamma})+\bar{\Phi}(2\sqrt{\gamma})\right)
\end{align*}
where $\bar{\Phi}$ denotes the tail distribution function of standard normal distribution. It is known that $\bar{\Phi}(x)=\Theta(\frac{1}{x}e^{-x^2/2})$ as $x\to\infty$, and hence $\tE(Z^2)=\Theta(\frac{1}{\sqrt{\gamma}}e^{(1-\frac{(k-1)^2}{2})\gamma})$ if $1<k<3$ and $\tE(Z^2)=\Theta(\frac{1}{\sqrt{\gamma}}e^{-\gamma})$ if $k\geq 3$. Besides, $p=\bar{\Phi}(\sqrt{\gamma})+\bar{\Phi}(k\sqrt{\gamma})=\Theta(\frac{1}{\sqrt{\gamma}}e^{-\gamma/2})$. Therefore, $\tE(Z^2)/p^2=\Theta(\sqrt{\gamma}e^{(2-\frac{(k-1)^2}{2})\gamma})$ which grows exponentially in $\gamma$ if $1<k<3$ and $\tE(Z^2)/p^2=\Theta(\sqrt{\gamma})$ which grows polynomially in $\gamma$ if $k\geq 3$.
\hfill\Halmos
\endproof

\proof{Proofs of Propositions \ref{prop:example2} and \ref{prop:example3}.} 
We first show a general result on IS estimator that exponentially tilts to the most significant dominating point. Consider a rare-event set $\cE$ with $a_1$ as the most significant dominating point. We estimate the target rare event $\{\frac{1}{\gamma}X_{\gamma}\in\cE\}$ using the IS estimator with likelihood ratio
$$
L = \frac{1}{e ^{s_{a_1}^\top X_{\gamma}-\gamma \mu_\gamma (s_{a_1})} }.
$$
That is, we use exponential tilting with respect to $a_1$ only. Specifically, we consider $X_{\gamma}\sim N(\gamma\lambda,\gamma\Sigma)$. In this case we have $\mu(x)=\mu_{\gamma}(x)=x^\top\lambda+\frac12x^\top\Sigma x$, $s_a=\Sigma^{-1}(a-\lambda)$, $I(y)=\frac12(y-\lambda)^\top\Sigma^{-1}(y-\lambda)$. 

The second moment of this IS estimator is
\begin{align*}
&\tE\left(I\left(\frac{1}{\gamma}X_{\gamma}\in\cE\right)L^2(\omega)\right)\\
=&E\left(I\left(\frac{1}{\gamma}X_{\gamma}\in\cE\right)L(\omega)\right)\\
=&\int_{\frac{x}{\gamma}\in\cE}e^{-s_{a_1}^\top x+\gamma \mu_\gamma (s_{a_1})}(2\pi)^{-\frac{d}{2}}|\gamma\Sigma|^{-1/2}e^{-\frac{1}{2\gamma}(x-\gamma\lambda)^\top\Sigma^{-1}(x-\gamma\lambda)}dx\\
=&\int_{\frac{x}{\gamma}\in\cE}e^{-(a_1-\lambda)^\top\Sigma^{-1}x+\gamma(a_1-\lambda)^\top\Sigma^{-1}\lambda+\frac{\gamma}{2}(a_1-\lambda)^\top\Sigma^{-1}(a_1-\lambda)}(2\pi)^{-\frac{d}{2}}|\gamma\Sigma|^{-1/2}e^{-\frac{1}{2\gamma}(x-\gamma\lambda)^\top\Sigma^{-1}(x-\gamma\lambda)}dx\\
=&\int_{\frac{x}{\gamma}\in\cE}(2\pi)^{-\frac{d}{2}}|\gamma\Sigma|^{-1/2}e^{\gamma[-(a_1-\lambda)^\top\Sigma^{-1}(\frac{x}{\gamma}-\lambda)+\frac12(a_1-\lambda)^\top \Sigma^{-1}(a_1-\lambda)-\frac12(\frac{x}{\gamma}-\lambda)^\top\Sigma^{-1}(\frac{x}{\gamma}-\lambda)]}dx\\
=&e^{\gamma(a_1-\lambda)^\top\Sigma^{-1}(a_1-\lambda)}\int_{\frac{x}{\gamma}\in\cE}(2\pi)^{-\frac{d}{2}}|\gamma\Sigma|^{-1/2}e^{-\frac{\gamma}{2}(\frac{x}{\gamma}+a_1-2\lambda)^\top\Sigma^{-1}(\frac{x}{\gamma}+a_1-2\lambda)}dx\\
=&e^{\gamma(a_1-\lambda)^\top\Sigma^{-1}(a_1-\lambda)}\bar{P}\left(\frac{1}{\gamma}X_{\gamma}\in\cE\right)
\end{align*}
where $\bar{P}$ is the probability measure given by the exponential tilting with respect to $2\lambda-a_1$ and hence $X_{\gamma}\sim N(\gamma(2\lambda-a_1),\gamma\Sigma)$ under $\bar{P}$. We correspondingly denote $\bar I$ as the rate function under $\bar P$. Now, since $X_\gamma$ is Gaussian Assumption \ref{asm:mu_GE} holds. Suppose also that Assumption \ref{asm:E_GE} holds. Then by Theorem~\ref{thm:GE} we know that 
$$
\lim_{\gamma\to\infty}\frac{1}{\gamma}\log\bar{P}\left(\frac{1}{\gamma}X_{\gamma}\in\cE\right)=-\bar{I}(\cE)=-\frac12\min_{y\in\cE}(y+a_1-2\lambda)^\top\Sigma^{-1}(y+a_1-2\lambda).
$$

We also know that $\lim_{\gamma\to\infty}\frac{1}{\gamma}\log P\left(\frac{1}{\gamma}X_{\gamma}\in\cE\right)=-I(a_1)=-\frac12(a_1-\lambda)^\top\Sigma^{-1}(a_1-\lambda)$. Therefore,
\begin{align*}
    &\lim_{\gamma\to\infty}\frac{\log\tE\left(I\left(\frac{1}{\gamma}X_{\gamma}\in\cE\right)L^2(\omega)\right)}{\log P\left(\frac{1}{\gamma}X_{\gamma}\in\cE\right)}\\
    =&\lim_{\gamma\to\infty}\frac{\frac{1}{\gamma}\log\tE\left(I\left(\frac{1}{\gamma}X_{\gamma}\in\cE\right)L^2(\omega)\right)}{\frac{1}{\gamma}\log P\left(\frac{1}{\gamma}X_{\gamma}\in\cE\right)}\\
    =&\lim_{\gamma\to\infty}\frac{(a_1-\lambda)^\top\Sigma^{-1}(a_1-\lambda)+\frac{1}{\gamma}\log\bar{P}\left(\frac{1}{\gamma}X_{\gamma}\in\cE\right)}{\frac{1}{\gamma}\log P\left(\frac{1}{\gamma}X_{\gamma}\in\cE\right)}\\
    =&\frac{(a_1-\lambda)^\top\Sigma^{-1}(a_1-\lambda)-\frac12\min_{y\in\cE}(y+a_1-2\lambda)^\top\Sigma^{-1}(y+a_1-2\lambda)}{-\frac12(a_1-\lambda)^\top\Sigma^{-1}(a_1-\lambda)}.
\end{align*}
By definition, the IS estimator is not asymptotically efficient if and only if $\min_{y\in\cE}(y+a_1-2\lambda)^\top\Sigma^{-1}(y+a_1-2\lambda)<4(a_1-\lambda)^\top\Sigma^{-1}(a_1-\lambda)$.

In order to check the existence of $\tilde{y}$ such that $(\tilde{y}+a_1-2\lambda)^\top\Sigma^{-1}(\tilde{y}+a_1-2\lambda)<4(a_1-\lambda)^\top\Sigma^{-1}(a_1-\lambda)$, we can formulate the following optimization\begin{equation}\label{eq:efficiency_optimization}
    \min_{y\in\cE}\quad (y+a_1-2\lambda)^\top\Sigma^{-1}(y+a_1-2\lambda)-4(a_1-\lambda)^\top\Sigma^{-1}(a_1-\lambda),
\end{equation}
and check whether the objective is negative for the optimal (or any feasible) solution. Since the objective function is quadratic on the decision vector $y$, the tractability of the optimization problem \eqref{eq:efficiency_optimization} is determined by the feasible region $\mathcal{E}$. In the example in Section~\ref{sec:numerical_sum}, we have $\bigcup_{m=1}^d\mathcal{H}_m=\bigcup_{m=1}^d\{x\in\R^d:\sum_{i=1}^mx_i\geq a\}$. Since each $\mathcal{H}_m$ is a half-space and hence convex, we can independently solve $$
    \min_{y\in \mathcal{H}_m}\quad (y+a_1-2\lambda)^\top\Sigma^{-1}(y+a_1-2\lambda)-4(a_1-\lambda)^\top\Sigma^{-1}(a_1-\lambda),
$$
for $m=1,...,d$. Since the feasible region of the above optimization associated with any $\mathcal{H}_m$ is a subset of the original problem \eqref{eq:efficiency_optimization}, any solution with negative objective indicates that the IS estimator with dominating points $a_1$ is not asymptotically efficient. In the example in Section~\ref{sec:numerical_mnist}, we formulate the problem \eqref{eq:efficiency_optimization} as a mixed integer programming problem following the treatment of the rare-event set formed by machine learning predictors in \cite{bai2022rare}. Using this approach, we can show that the IS estimators with the most dominating points are not asymptotically efficient for both examples considered in Sections~\ref{sec:numerical_sum} and~\ref{sec:numerical_mnist}.\hfill\Halmos
\endproof

To prove Theorem \ref{thm:general_interval1_GE}, we need the following lemma:
\begin{lemma}[Theorem 4 in \citealt{maurer2009empirical}]
Let $Y,Y_1,\dots,Y_n$ be i.i.d. random variables with values in $[0,1]$ and let $\delta>0$. Then with probability at least $1-\delta$
\begin{equation*}
EY-\frac{1}{n}\sum_{i=1}^nY_i\leq\sqrt{\frac{2V_n(\bY)\log(2/\delta)}{n}}+\frac{7\log(2/\delta)}{3(n-1)}
\end{equation*}
where $\bY=(Y_1,\dots,Y_n)$ and $V_n(\bY)$ is the sample variance of $Y_i$'s.
\label{lem:concentration}
\end{lemma}
Now we prove Theorem \ref{thm:general_interval1_GE}:
\proof{Proof of Theorem \ref{thm:general_interval1_GE}.}
Following the notation in Section \ref{sec:concept}, we split $\cE$ into $\cE_1=\cE\cap\bigcup_{i=1}^k\{x\in\R^d:s_{a_i}^\top (x-a_i)\geq 0\}$ and $\cE_2=\cE\setminus\cE_1$, and then we correspondingly define $\cA_{\gamma}^1=\{\frac{1}{\gamma}X_{\gamma}\in\cE_1\}$ and $\cA_{\gamma}^2=\{\frac{1}{\gamma}X_{\gamma}\in\cE_2\}$. 
We denote $Z_1^{(i)}=I(\frac{1}{\gamma}X^{(i)}\in\cE_1)\frac{dP}{d\tP}$. $\hat p_1$ and $\hat V_1$ respectively denote the sample mean and sample variance of $Z_1^{(i)}$'s. We also define $N=\sum_{i=1}^nI(\frac{1}{\gamma}X_{\gamma}^{(i)}\in\cE_2)$. Note that conditional on $N=0$, we have $Z^{(i)}=Z_1^{(i)}$. Then 
\begin{align*}
    &\tP\left(|\hat p-p|>\sqrt{\frac{2\hat{V}\log(4/\alpha)}{n}}+\frac{7\log(4/\alpha)M_{\gamma}}{3(n-1)}\right)\\
    =&\tP\left(|\hat p-p|>\sqrt{\frac{2\hat{V}\log(4/\alpha)}{n}}+\frac{7\log(4/\alpha)M_{\gamma}}{3(n-1)},N=0\right)\\
    &+\tP\left(|\hat p-p|>\sqrt{\frac{2\hat{V}\log(4/\alpha)}{n}}+\frac{7\log(4/\alpha)M_{\gamma}}{3(n-1)},N>0\right)\\
    \leq &\tP\left(|\hat p_1-p|>\sqrt{\frac{2\hat{V}_1\log(4/\alpha)}{n}}+\frac{7\log(4/\alpha)M_{\gamma}}{3(n-1)},N=0\right)+\tP(N>0)\\
    \leq &\tP\left(|\hat p_1-p|>\sqrt{\frac{2\hat{V}_1\log(4/\alpha)}{n}}+\frac{7\log(4/\alpha)M_{\gamma}}{3(n-1)}\right)+\tP(N>0)\\
    \leq &\tP\left(|\hat p_1-p_1|+p_2>\sqrt{\frac{2\hat{V}_1\log(4/\alpha)}{n}}+\frac{7\log(4/\alpha)M_{\gamma}}{3(n-1)}\right)+n\tilde{p}_2.
\end{align*}
Similar to the derivation of \eqref{eqn:upper_bound_GE}, we know that $0\leq Z_1^{(i)}\leq M_{\gamma},\forall i$. By applying Lemma \ref{lem:concentration} with $Y_i=Z_1^{(i)}/M_{\gamma}$ and $Y_i=1-Z_1^{(i)}/M_{\gamma}$ respectively, we get that 
\begin{equation*}
\tP\left(p_1>\hp_1+\sqrt{\frac{2\hat{V}_1\log(4/\delta)}{n}}+\frac{7\log(4/\delta)M_{\gamma}}{3(n-1)}\right)\leq \delta/2
\end{equation*}
and
\begin{equation*}
\tP\left(p_1<\hp_1-\sqrt{\frac{2\hat{V}_1\log(4/\delta)}{n}}-\frac{7\log(4/\delta)M_{\gamma}}{3(n-1)}\right)\leq \delta/2
\end{equation*}
for any $\delta>0$. Thus,
\begin{equation}
\tP\left(|\hp_1-p_1|>\sqrt{\frac{2\hat{V}_1\log(4/\delta)}{n}}+\frac{7\log(4/\delta)M_{\gamma}}{3(n-1)}\right)\leq \delta.
\label{eqn:concentration_GE}
\end{equation}
Find $\alpha'=\alpha'(\gamma)$ such that 
\begin{equation*}
\frac{7\log(4/\alpha')M_{\gamma}}{3(n-1)}=\frac{7\log(4/\alpha)M_{\gamma}}{3(n-1)}-p_2.
\end{equation*}
That is, 
\begin{equation*}
\alpha'=\alpha\exp\left(\frac{3(n-1)p_2}{7M_{\gamma}}\right).
\end{equation*}
Clearly $\alpha'>\alpha$ and $\log(4/\alpha')<\log(4/\alpha)$. From the proof of Theorem \ref{thm:general_point_estimate_GE}, we have that either (i) $k=r$ and $p_2=0$ or (ii) $k<r$ and $\limsup_{\gamma\to\infty}\frac{1}{\gamma}\log p_2\leq -I(a_{k+1})<-I(a_1)$. Moreover, we have 
\begin{equation*}
    \frac{1}{\gamma}\log M_{\gamma}=\max_{i=1,\dots,k}\left\{-\frac{1}{\gamma}\log\alpha_i-(s_{a_i}^\top a_i-\mu_{\gamma}(s_{a_i}))\right\}\to -I(a_1).
\end{equation*}
Thus, for subexponentially growing $n$, we have $\frac{3(n-1)p_2}{7M_{\gamma}}\to 0$ and hence $\alpha'\to\alpha$ as $\gamma\to\infty$. We replace $\delta$ with $\alpha'$ in \eqref{eqn:concentration_GE}, and then we get 
\begin{equation*}
\tP\left(|\hp_1-p_1|>\sqrt{\frac{2\hat{V}_1\log(4/\alpha')}{n}}+\frac{7\log(4/\alpha')M_{\gamma}}{3(n-1)}\right)\leq \alpha'.
\end{equation*}
Hence
\begin{equation*}
\tP\left(|\hp_1-p_1|>\sqrt{\frac{2\hat{V}_1\log(4/\alpha)}{n}}+\frac{7\log(4/\alpha)M_{\gamma}}{3(n-1)}-p_2\right)\leq \alpha'.
\end{equation*}
Therefore,
\begin{equation*}
\tP\left(|\hat p-p|>\sqrt{\frac{2\hat{V}\log(4/\alpha)}{n}}+\frac{7\log(4/\alpha)M_{\gamma}}{3(n-1)}\right)\leq \alpha'+n\tilde{p}_2\to\alpha\text{ as }\gamma\to\infty.
\end{equation*}\hfill\Halmos
\endproof

To prove Theorem \ref{thm:general_interval2_GE}, we first prove a lemma regarding the variance $\tVar(Z_1)$. In particular, we have assumed that $\tVar(Z_1)$ cannot be ``too large'', while this lemma implies that it cannot be ``too small'' as well. As we will see, this lemma is used to control the error of the normal approximation in analog to the Berry-Esseen theorem in order to argue the validity of $\mathcal I_2$.

\begin{lemma}
Under the same setting as Theorem \ref{thm:general_point_estimate_GE}, split $\cE$ into $\cE_1=\cE\cap\bigcup_{i=1}^k\{x\in\R^d:s_{a_i}^\top (x-a_i)\geq 0\}$ and $\cE_2=\cE\setminus\cE_1$. Denote $Z_1=I(\frac{1}{\gamma}X_{\gamma}\in\cE_1)\frac{dP}{d\tP}(\omega)$ with $X_{\gamma}\sim\tP$. Then $\lim_{\gamma\to\infty}\frac{1}{\gamma}\log\tVar(Z_1)=-2I(\cE)=-2I(a_1)$.
\label{lem:variance}
\end{lemma}
\proof{Proof of Lemma \ref{lem:variance}.}
From Proposition \ref{classical AE}, we know that $\limsup_{\gamma\to\infty}\frac{1}{\gamma}\log\tVar(Z_1)\leq -2I(\cE)$. Thus we only need to show that $\liminf_{\gamma\to\infty}\frac{1}{\gamma}\log\tVar(Z_1)\geq -2I(\cE)$. Indeed, we have that $\tVar(Z_1)=\tE((Z_1-p_1)^2)\geq \tE((Z_1-p_1)^2I(\frac{1}{\gamma}X_{\gamma}\notin\cE_1))=p_1^2\tP(\frac{1}{\gamma}X_{\gamma}\in\cE_1^c)\geq p_1^2\alpha_1\tP_1(\frac{1}{\gamma}X_{\gamma}\in\cE_1^c)$ where $\tP_1$ is as defined in the proof of Theorem \ref{thm:general_point_estimate_GE}, and hence $\liminf_{\gamma\to\infty}\frac{1}{\gamma}\log\tVar(Z_1)\geq\liminf_{\gamma\to\infty}\left(\frac{2}{\gamma}\log p_1+\frac{1}{\gamma}\log\alpha_1+\frac{1}{\gamma}\log\tP_1(\frac{1}{\gamma}X_{\gamma}\in\cE_1^c)\right)$. First, we know that $\lim_{\gamma\to\infty}\frac{1}{\gamma}\log p_1=-I(\cE)$. Second, we know that $\frac{1}{\gamma}\log\alpha_1\to 0$. Third, by Theorem \ref{thm:GE}, we know that $\liminf_{\gamma\to\infty}\frac{1}{\gamma}\log\tP_1(\frac{1}{\gamma}X_{\gamma}\in\cE_1^c)\geq -\tilde{I}_1((\cE_1^c)^\circ)$ where $\tilde{I}_1(y)=I(y)-s_{a_1}^\top y+\mu(s_{a_1})$ as in the proof of Theorem \ref{thm:general_point_estimate_GE}. As $a_1$ is a dominating point, we have required that $a_1\in\mathcal{D}(I)^\circ$, and thus $a_1\in\mathcal{D}(\tilde{I}_1)^\circ$. We note that $(\cE_1^c)^\circ\supset\bigcap_{i=1}^k\{x:s_{a_i}^\top(x-a_i)<0\}$. Then $\tilde{I}_1((\cE_1^c)^\circ)\leq\tilde{I}_1(\bigcap_{i=1}^k\{x:s_{a_i}^\top(x-a_i)<0\})$. We note that since $a_1$ is the most significant dominating point, we have $s_{a_i}(a_1-a_i)\leq 0$ (otherwise we get $I(a_i)<I(a_1)$, which is a contradiction). Thus $a_1$ is on the boundary of $\bigcap_{i=1}^k\{x:s_{a_i}^\top(x-a_i)<0\}$, and there exists a sequence of points $\{y_n\}_{n=1}^{\infty}\subset\mathcal{D}(\tilde{I}_1)^\circ\cap\bigcap_{i=1}^k\{x:s_{a_i}^\top(x-a_i)<0\}$ such that $y_n\to a_1$ as $n\to\infty$. Therefore, we get $\tilde{I}_1(\bigcap_{i=1}^k\{x:s_{a_i}^\top(x-a_i)<0\})\leq \tilde{I}_1(a_1)=0$. Overall, we have proved that $\liminf_{\gamma\to\infty}\frac{1}{\gamma}\log\tVar(Z_1)\geq -2I(\cE)$.
\hfill\halmos
\endproof

We also need a concentration result for the sample variance:
\begin{lemma}[Theorem 10 in \citealt{maurer2009empirical}]
Let $n\geq 2$ and $\bY=(Y_1,\dots,Y_n)$ be a vector of independent random variables with values in $[0,1]$. Then for $\delta>0$ we have
\begin{align*}
    P\left(\sqrt{EV_n(\bY)}>\sqrt{V_n(\bY)}+\sqrt{\frac{2\log(1/\delta)}{n-1}}\right)&\leq\delta;\\
    P\left(\sqrt{V_n(\bY)}>\sqrt{EV_n(\bY)}+\sqrt{\frac{2\log(1/\delta)}{n-1}}\right)&\leq\delta.
\end{align*}
\label{lem:concentration_var}
\end{lemma}


With Lemmas \ref{lem:variance} and \ref{lem:concentration_var}, now we prove Theorem \ref{thm:general_interval2_GE}:
\proof{Proof of Theorem \ref{thm:general_interval2_GE}.}
We denote $Z_1^{(i)}=I(\frac{1}{\gamma}X^{(i)}\in\cE_1)\frac{dP}{d\tP}$. $\hat p_1$ and $\hat V_1$ respectively denote the sample mean and sample variance of $Z_1^{(i)}$'s. We also define $N=\sum_{i=1}^nI(\frac{1}{\gamma}X_{\gamma}^{(i)}\in\cE_2)$. Similar to the proof of Theorem \ref{thm:general_interval1_GE}, we have 
\begin{align*}
    &\tP\left(|\hat p-p|>z_{1-\alpha/2}\sqrt{\frac{\hat{V}}{n}}\right)\\
    = &\tP\left(|\hat p-p|>z_{1-\alpha/2}\sqrt{\frac{\hat{V}}{n}},N=0\right)+\tP\left(|\hat p-p|>z_{1-\alpha/2}\sqrt{\frac{\hat{V}}{n}},N>0\right)\\
    \leq &\tP\left(|\hat p_1-p|>z_{1-\alpha/2}\sqrt{\frac{\hat{V}_1}{n}}\right)+\tP(N>0)\\
    \leq &\tP\left(|\hat p_1-p_1|>z_{1-\alpha/2}\sqrt{\frac{\hat{V}_1}{n}}-p_2\right)+n\tilde p_2\\
    \leq &\tP\left(|\hat p_1-p_1|>z_{1-\alpha/2}\left(\sqrt{\frac{\tE \hat{V}_1}{n}}-\sqrt{\frac{2\log(1/\delta)}{n(n-1)}}M_{\gamma}\right)-p_2\right)+\tP\left(\sqrt{\tE \hat{V}_1}>\sqrt{\hat{V}_1}+\sqrt{\frac{2\log(1/\delta)}{n-1}}M_{\gamma}\right)+n\tilde p_2
\end{align*}
for any $\delta=\delta(\gamma)>0$, where $M_{\gamma}$ is as defined in Theorem~\ref{thm:general_interval1_GE}. We know that $0\leq Z_1^{(i)}\leq M_{\gamma},\forall i$. By Lemma \ref{lem:concentration_var} with $Y_i=Z_1^{(i)}/M_{\gamma}$, we get that 
\begin{equation*}
\tP\left(\sqrt{\tE \hat{V}_1}>\sqrt{\hat{V}_1}+\sqrt{\frac{2\log(1/\delta)}{n-1}}M_{\gamma}\right)\leq\delta.
\end{equation*}
By Berry-Esseen theorem, we know that for any $x\in\R$
\begin{equation*}
\left|\tP\left(\frac{\sqrt{n}(\hp_1-p_1)}{\sqrt{\tVar(Z_1^{(1)})}}\leq x\right)-\Phi(x)\right|\leq\frac{C\tE|Z_1^{(1)}-p_1|^3}{\tVar^{3/2}(Z_1^{(1)})\sqrt{n}}
\end{equation*}
where $\Phi$ is the CDF of standard normal distribution and $C$ is a universal constant. Let 
\begin{equation*}
x=z_{1-\alpha/2}\left(\sqrt{\frac{\tE \hat{V}_1}{\tVar(Z_1^{(1)})}}-\sqrt{\frac{2\log(1/\delta)}{(n-1)\tVar(Z_1^{(1)})}}M_{\gamma}\right)-p_2\sqrt{\frac{n}{\tVar(Z_1^{(1)})}}.
\end{equation*}
Then we get that 
\begin{equation*}
\tP\left(\frac{\sqrt{n}|\hp_1-p_1|}{\sqrt{\tVar(Z_1^{(1)})}}> x\right)\leq 2\Phi(-x)+\frac{2C\tE|Z_1^{(1)}-p_1|^3}{\tVar^{3/2}(Z_1^{(1)})\sqrt{n}}.
\end{equation*}
Hence,
\begin{equation*}
\tP\left(|\hat p-p|>z_{1-\alpha/2}\sqrt{\frac{\hat{V}}{n}}\right)\leq 2\Phi(-x)+\frac{2C\tE|Z_1^{(1)}-p_1|^3}{\tVar^{3/2}(Z_1^{(1)})\sqrt{n}}+\delta+n\tilde{p}_2.
\end{equation*}
First, we have that $\tE|Z_1^{(1)}-p_1|^3\leq \max(p_1^3,(M_{\gamma}-p_1)^3)$. From the proof of Theorems~\ref{thm:general_point_estimate_GE} and \ref{thm:general_interval1_GE}, we know that $-\frac{1}{\gamma}\log p_1\to -I(a_1)$ and $-\frac{1}{\gamma}\log M_{\gamma}\to -I(a_1)$. Thus $\limsup_{\gamma\to\infty}\frac{1}{\gamma}\log\tE|Z_1^{(1)}-p_1|^3\leq -3I(a_1)$. By Lemma \ref{lem:variance}, we get that $\frac{\tE|Z_1^{(1)}-p_1|^3}{\tVar^{3/2}(Z_1^{(1)})}$ grows at most subexponentially in $\gamma$. Hence under the assumptions we could choose $n$ as required in the theorem. In particular, we have that 
\begin{equation*}
\frac{2C\tE|Z_1^{(1)}-p_1|^3}{\tVar^{3/2}(Z_1^{(1)})\sqrt{n}}\to 0.
\end{equation*}
Now we analyze $x$. We know that $\sqrt{\frac{E\hat{V}_1}{\tVar(Z_1^{(1)})}}=1$ and 
\begin{equation*}
p_2\sqrt{\frac{n}{\tVar(Z_1^{(1)})}}=\frac{p_2}{e^{-\gamma I(a_1)}}\sqrt{\frac{ne^{-2\gamma I(a_1)}}{\tVar(Z_1^{(1)})}}\to 0
\end{equation*}
since $\frac{p_2}{e^{-\gamma I(a_1)}}$ decays exponentially (proof of Theorem~\ref{thm:general_point_estimate_GE}) while $n$ and $\frac{e^{-2\gamma I(a_1)}}{\tVar(Z_1^{(1)})}$ grow subexponentially (Lemma \ref{lem:variance}) in $\gamma$. Now we consider
\begin{equation*}
\sqrt{\frac{2\log(1/\delta)}{(n-1)\tVar(Z_1^{(1)})}}M_{\gamma}=\sqrt{\frac{2\log(1/\delta)n}{n-1}}\sqrt{\frac{M_{\gamma}^2}{n\tVar(Z_1^{(1)})}}.
\end{equation*}
Since we assume that $\frac{M_{\gamma}^2}{n\tVar(Z_1^{(1)})}\to 0$, we could set $\delta$ such that $\delta\to 0$ and 
\begin{equation*}
\sqrt{\frac{2\log(1/\delta)}{(n-1)\tVar(Z_1^{(1)})}}M_{\gamma}\to 0.
\end{equation*}
In this case, $x\to z_{1-\alpha/2}$ and hence $\Phi(-x)\to \alpha/2$. Combining all the results, we get that 
\begin{equation*}
\limsup_{\gamma\to\infty}\tP\left(|\hat p-p|>z_{1-\alpha/2}\sqrt{\frac{\hat{V}}{n}}\right)\leq\alpha.
\end{equation*}

\hfill\Halmos
\endproof

\proof{Proof of Theorem \ref{thm:tight_interval_full_IS}.}
First, following the proof of Lemma \ref{lem:variance}, it is easy to get that $\lim_{\gamma\to\infty}\frac{1}{\gamma}\log\tVar(Z)=-2I(\cE)=-2I(a_1)$. Moreover, we have that 
\begin{align*}
    \tP\left(|\hat p-p|>z_{1-\alpha/2}\sqrt{\frac{\hat{V}}{n}}\right)
    \leq &\tP\left(|\hp-p|>z_{1-\alpha/2}\left(\sqrt{\frac{\tE \hat{V}}{n}}-\sqrt{\frac{2\log(1/\delta)}{n(n-1)}}M_{\gamma}\right)\right)\\
    &+\tP\left(\sqrt{\tE \hat{V}}>\sqrt{\hat{V}}+\sqrt{\frac{2\log(1/\delta)}{n-1}}M_{\gamma}\right)
\end{align*}
and 
\begin{align*}
    \tP\left(|\hat p-p|>z_{1-\alpha/2}\sqrt{\frac{\hat{V}}{n}}\right)
    \geq &\tP\left(|\hp-p|>z_{1-\alpha/2}\left(\sqrt{\frac{\tE \hat{V}}{n}}+\sqrt{\frac{2\log(1/\delta)}{n(n-1)}}M_{\gamma}\right)\right)\\
    &-\tP\left(\sqrt{\hat{V}}>\sqrt{\tE \hat{V}}+\sqrt{\frac{2\log(1/\delta)}{n-1}}M_{\gamma}\right)
\end{align*}
for any $\delta=\delta(\gamma)>0$. We know that $0\leq Z^{(i)}\leq M_{\gamma},\forall i$. By Lemma~\ref{lem:concentration_var}, we get that 
\begin{equation*}
    \tP\left(\sqrt{\tE \hat{V}}>\sqrt{\hat{V}}+\sqrt{\frac{2\log(1/\delta)}{n-1}}M_{\gamma}\right)\leq \delta
\end{equation*}
and
\begin{equation*}
    \tP\left(\sqrt{\hat{V}}>\sqrt{\tE \hat{V}}+\sqrt{\frac{2\log(1/\delta)}{n-1}}M_{\gamma}\right)\leq\delta.
\end{equation*}
By Berry-Esseen theorem, we know that 
\begin{equation*}
\tP\left(\frac{\sqrt{n}|\hp-p|}{\sqrt{\tVar(Z^{(1)})}}> x_1\right)\leq 2\Phi(-x_1)+\frac{2C\tE|Z^{(1)}-p|^3}{\tVar^{3/2}(Z^{(1)})\sqrt{n}}
\end{equation*}
and 
\begin{equation*}
\tP\left(\frac{\sqrt{n}|\hp-p|}{\sqrt{\tVar(Z^{(1)})}}> x_2\right)\geq 2\Phi(-x_2)-\frac{2C\tE|Z^{(1)}-p|^3}{\tVar^{3/2}(Z^{(1)})\sqrt{n}}
\end{equation*}
where $\Phi$ is the CDF of standard normal distribution, $C$ is a universal constant, and
\begin{align*}
    x_1&=z_{1-\alpha/2}\left(\sqrt{\frac{\tE \hat{V}}{\tVar(Z^{(1)})}}-\sqrt{\frac{2\log(1/\delta)}{(n-1)\tVar(Z^{(1)})}}M_{\gamma}\right),\\
    x_2&=z_{1-\alpha/2}\left(\sqrt{\frac{\tE \hat{V}}{\tVar(Z^{(1)})}}+\sqrt{\frac{2\log(1/\delta)}{(n-1)\tVar(Z^{(1)})}}M_{\gamma}\right).
\end{align*}
Combining the above derivations, we get that 
\begin{align*}
    \tP\left(|\hat p-p|>z_{1-\alpha/2}\sqrt{\frac{\hat{V}}{n}}\right)&
    \leq 2\Phi(-x_1)+\frac{2C\tE|Z^{(1)}-p|^3}{\tVar^{3/2}(Z^{(1)})\sqrt{n}}+\delta,\\
    \tP\left(|\hat p-p|>z_{1-\alpha/2}\sqrt{\frac{\hat{V}}{n}}\right)&
    \geq 2\Phi(-x_2)-\frac{2C\tE|Z^{(1)}-p|^3}{\tVar^{3/2}(Z^{(1)})\sqrt{n}}-\delta.
\end{align*}
First, we have that $\tE|Z^{(1)}-p|^3\leq \max(p^3,(M_{\gamma}-p)^3)$. From the proof of the previous theorems, we know that $-\frac{1}{\gamma}\log p\to -I(a_1)$ and $-\frac{1}{\gamma}\log M_{\gamma}\to -I(a_1)$. Thus, $\frac{\tE|Z^{(1)}-p|^3}{\tVar^{3/2}(Z^{(1)})}$ grows at most subexponentially in $\gamma$, and hence we could choose $n$ as required in the theorem. In particular, we have $\frac{2C\tE|Z^{(1)}-p|^3}{\tVar^{3/2}(Z^{(1)})\sqrt{n}}\to 0$. Next, since we assume that $\frac{M_{\gamma}^2}{n\tVar(Z^{(1)})}\to 0$, we could set $\delta$ such that $\delta\to 0$ and 
\begin{equation*}
    \sqrt{\frac{2\log(1/\delta)}{(n-1)\tVar(Z^{(1)})}}M_{\gamma}\to 0.
\end{equation*}
Hence, $x_1,x_2\to z_{1-\alpha/2}$ as $\gamma\to\infty$. Overall, we get that 
\begin{align*}
    \limsup_{\gamma\to\infty}\tP\left(|\hat p-p|>z_{1-\alpha/2}\sqrt{\frac{\hat{V}}{n}}\right)&\leq \alpha,\\
    \liminf_{\gamma\to\infty}\tP\left(|\hat p-p|>z_{1-\alpha/2}\sqrt{\frac{\hat{V}}{n}}\right)&\geq \alpha.
\end{align*}
Therefore, the theorem is proved.
\hfill\Halmos
\endproof

\proof{Proof of Theorem \ref{PE examples}.}
The conclusion follows from Theorem \ref{thm:unique}, by checking, for each setting in Sections \ref{sec:numerical_glasserman}, \ref{sec:numerical_sum} and \ref{sec:numerical_mnist}, Assumptions \ref{asm:mu_GE} and \ref{asm:E_GE} hold and the most significant dominating point is unique. Since we have $\mu(x)=1.5x+0.5 x^2- \log(1+x)$ for the case in Section~\ref{sec:numerical_glasserman} and $\mu(x)=\lambda^\top x + \frac{1}{2} x^\top \Sigma x$ for the cases in Sections \ref{sec:numerical_sum} and \ref{sec:numerical_mnist}, we verify the conditions in Assumption~\ref{asm:mu_GE}. On the other hand, since all rare-event sets in these examples are closed and contain unique optimal solutions for minimizing the corresponding rate function $I(x)$, we can verify Assumption \ref{asm:E_GE} and the uniqueness of most significant dominating point.  \hfill\Halmos
\endproof

\proof{Proof of Theorem \ref{thm:Weakly-PE-examples}.}
We check the three conditions in Theorem~\ref{thm:weak_probabilistic_efficiency} to show the probabilistic efficiency of the IS estimator using the dominating point $a_1$. Using the notation in Section \ref{sec:two-sided random walk}, we note that $\cA^1_\gamma$ only one most significant dominating point $a_{1}$ and $\cA^2_\gamma$ only one most significant dominating point $-a_{1}$, while $I(a_{1})=I(-a_{1})$. This indicates that $p_1/p_2 \to 1$ as $\gamma \to \infty$ and hence $p_1/p \to 1/2$ as $\gamma \to \infty$. Then, following the argument in Section~\ref{sec:experiment_validate_random_walk}, the IS estimator using the dominating point $a_1$ is asymptotically efficient for $\cA^1_\gamma$. Lastly we have $\tilde{p}_2 \to 0 $ exponentially fast and hence we have $n \tilde{p}_2  \to 0$ with subexponentially growing sample size $n$. We conclude that the IS estimator using the dominating point $a_1$ is weakly probabilistically efficient. 

On the other hand, based on Proposition~\ref{prop:example2}, the IS estimator using the dominating point $a_+$ is not asymptotically efficient for $\{ \bigcup_{m=1}^d \cH^+_m \}$ and hence is not asymptotically efficient for $\left\{  \left( \bigcup_{m=1}^d \cH^+_m \right) \bigcup \left( \bigcup_{m=1}^d \cH^-_m  \right)   \right\}$.\hfill\Halmos
\endproof

\proof{Proof of Theorem \ref{thm:general_point_estimate}.}
We only need to verify the conditions in Theorem \ref{thm:strong_probabilistic_efficiency}. First, from Assumption \ref{asm:p2_negligible}, we get that $p_2$ is exponentially smaller than $p_1$, and hence $\frac{p_1}{p}\to1$ as $\gamma\to\infty$. Second, Let $Z_1=I(X\in\cE_1)\frac{dP}{d\tP}$ under $\tP$. Clearly $\{a_1,\dots,a_k\}$ is a dominating set for $\cE_1$. Similar to \eqref{eqn:vartilde}, we get that $\tVar(Z_1)\leq k^2e^{-2I(a_1)}\leq r^2e^{-2I(a_1)}$. Assumption \ref{asm:p2_negligible} gives that $p_1\sim e^{-I(a_1)}$, and hence $Z_1$ is an asymptotically efficient estimator for $p_1$. Third, Assumption \ref{asm:p2tilde} implies that $n\tp_2\to 0$ for any $n$ subexponentially growing in $-\log p$. Therefore, all the conditions in Theorem \ref{thm:strong_probabilistic_efficiency} hold.\hfill\Halmos
\endproof

\proof{Proof of Theorem \ref{thm:general_interval1}.}
We denote $\hp_1=\frac{1}{n}\sum_{i=1}^nZ_1^{(i)}$ and $\bZ_1=(Z_1^{(1)},\dots,Z_1^{(n)})$. We also define $N=\sum_{i=1}^nI(X^{(i)}\in\cE_2)$. Note that conditional on $N=0$, we have $\bZ=\bZ_1$. Then 
\begin{align*}
    &\tilde{P}\left(|\hat p-p|>\sqrt{\frac{2\hat{V}\log(4/\alpha)}{n}}+\frac{7\log(4/\alpha)ke^{-I(a_1)}}{3(n-1)}\right)\\
    =&\tilde{P}\left(|\hat p-p|>\sqrt{\frac{2\hat{V}\log(4/\alpha)}{n}}+\frac{7\log(4/\alpha)ke^{-I(a_1)}}{3(n-1)},N=0\right)\\
    &+\tilde{P}\left(|\hat p-p|>\sqrt{\frac{2\hat{V}\log(4/\alpha)}{n}}+\frac{7\log(4/\alpha)ke^{-I(a_1)}}{3(n-1)},N>0\right)\\
    \leq &\tilde{P}\left(|\hat p_1-p|>\sqrt{\frac{2\hat{V}_1\log(4/\alpha)}{n}}+\frac{7\log(4/\alpha)ke^{-I(a_1)}}{3(n-1)},N=0\right)+P(N>0)\\
    \leq &\tilde{P}\left(|\hat p_1-p|>\sqrt{\frac{2\hat{V}_1\log(4/\alpha)}{n}}+\frac{7\log(4/\alpha)ke^{-I(a_1)}}{3(n-1)}\right)+P(N>0)\\
    \leq &\tilde{P}\left(|\hat p_1-p_1|+p_2>\sqrt{\frac{2\hat{V}_1\log(4/\alpha)}{n}}+\frac{7\log(4/\alpha)ke^{-I(a_1)}}{3(n-1)}\right)+n\tilde{p}_2.
\end{align*}
Similar to the derivation of \eqref{eqn:vartilde}, we know that $0\leq Z_1^{(i)}\leq ke^{-I(a_1)},\forall i$. By applying Lemma \ref{lem:concentration}  with $Y_i=Z_1^{(i)}/(ke^{-I(a_1)})$ and $Y_i=1-Z_1^{(i)}/(ke^{-I(a_1)})$ respectively, we get that 
\begin{equation*}
\tilde{P}\left(p_1>\hp_1+\sqrt{\frac{2\hat{V}_1\log(4/\delta)}{n}}+\frac{7\log(4/\delta)ke^{-I(a_1)}}{3(n-1)}\right)\leq \delta/2
\end{equation*}
and
\begin{equation*}
\tilde{P}\left(p_1<\hp_1-\sqrt{\frac{2\hat{V}_1\log(4/\delta)}{n}}-\frac{7\log(4/\delta)ke^{-I(a_1)}}{3(n-1)}\right)\leq \delta/2
\end{equation*}
for any $\delta>0$. Thus,
\begin{equation}
\tilde{P}\left(|\hp_1-p_1|>\sqrt{\frac{2\hat{V}_1\log(4/\delta)}{n}}+\frac{7\log(4/\delta)ke^{-I(a_1)}}{3(n-1)}\right)\leq \delta.
\label{eqn:concentration}
\end{equation}
Find $\alpha'=\alpha'(\gamma)$ such that 
\begin{equation*}
\frac{7\log(4/\alpha')ke^{-I(a_1)}}{3(n-1)}=\frac{7\log(4/\alpha)ke^{-I(a_1)}}{3(n-1)}-p_2.
\end{equation*}
That is, 
\begin{equation*}
\alpha'=\alpha\exp\left(\frac{3(n-1)p_2}{7ke^{-I(a_1)}}\right).
\end{equation*}
Clearly $\alpha'>\alpha$ and $\log(4/\alpha')<\log(4/\alpha)$. Moreover, we know that $\frac{p_2}{e^{-I(a_1)}}$ decays exponentially in $-\log p$ since Assumption \ref{asm:p2_negligible} holds and that $n$ grows subexponentially in $-\log p$, and thus $\alpha'\to\alpha$ as $\gamma\to\infty$. We replace $\delta$ with $\alpha'$ in \eqref{eqn:concentration}, and then we get 
\begin{equation*}
\tilde{P}\left(|\hp_1-p_1|>\sqrt{\frac{2\hat{V}_1\log(4/\alpha')}{n}}+\frac{7\log(4/\alpha')ke^{-I(a_1)}}{3(n-1)}\right)\leq \alpha'.
\end{equation*}
Hence
\begin{equation*}
\tilde{P}\left(|\hp_1-p_1|>\sqrt{\frac{2\hat{V}_1\log(4/\alpha)}{n}}+\frac{7\log(4/\alpha)ke^{-I(a_1)}}{3(n-1)}-p_2\right)\leq \alpha'.
\end{equation*}
Therefore,
\begin{equation*}
\tilde{P}\left(|\hat p-p|>\sqrt{\frac{2\hat{V}\log(4/\alpha)}{n}}+\frac{7\log(4/\alpha)ke^{-I(a_1)}}{3(n-1)}\right)\leq \alpha'+n\tilde{p}_2\to\alpha\text{ as }\gamma\to\infty.
\end{equation*}\hfill\Halmos
\endproof

\proof{Proof of Theorem \ref{thm:general_interval2}.}
We denote $\hp_1=\frac{1}{n}\sum_{i=1}^nZ_1^{(i)}$ and $\bZ_1=(Z_1^{(1)},\dots,Z_1^{(n)})$. We also define $N=\sum_{i=1}^nI(X^{(i)}\in\cE_2)$. We have 
\begin{align*}
    &\tilde{P}\left(|\hat p-p|>z_{1-\alpha/2}\sqrt{\frac{\hat{V}}{n}}\right)\\
    = &\tilde{P}\left(|\hat p-p|>z_{1-\alpha/2}\sqrt{\frac{\hat{V}}{n}},N=0\right)+\tilde{P}\left(|\hat p-p|>z_{1-\alpha/2}\sqrt{\frac{\hat{V}}{n}},N>0\right)\\
    \leq &\tilde{P}\left(|\hat p_1-p|>z_{1-\alpha/2}\sqrt{\frac{\hat{V}_1}{n}}\right)+P(N>0)\\
    \leq &\tilde{P}\left(|\hat p_1-p_1|>z_{1-\alpha/2}\sqrt{\frac{\hat{V}_1}{n}}-p_2\right)+n\tilde p_2\\
    \leq &\tilde{P}\left(|\hat p_1-p_1|>z_{1-\alpha/2}\left(\sqrt{\frac{\tE \hat{V}_1}{n}}-\sqrt{\frac{2\log(1/\delta)}{n(n-1)}}ke^{-I(a_1)}\right)-p_2\right)\\
    &+\tilde{P}\left(\sqrt{\tE \hat{V}_1}>\sqrt{\hat{V}_1}+\sqrt{\frac{2\log(1/\delta)}{n-1}}ke^{-I(a_1)}\right)+n\tilde p_2
\end{align*}
for any $\delta=\delta(\gamma)>0$. We know that $0\leq Z_1^{(i)}\leq ke^{-I(a_1)},\forall i$. By Lemma \ref{lem:concentration_var} with $Y_i=Z_1^{(i)}/(ke^{-I(a_1)})$, we get that 
\begin{equation*}
\tilde{P}\left(\sqrt{\tE \hat{V}_1}>\sqrt{\hat{V}_1}+\sqrt{\frac{2\log(1/\delta)}{n-1}}ke^{-I(a_1)}\right)\leq\delta.
\end{equation*}
By Berry-Esseen theorem, we know that for any $x\in\R$
\begin{equation*}
\left|\tilde{P}\left(\frac{\sqrt{n}(\hp_1-p_1)}{\sqrt{\widetilde{Var}(Z_1^{(1)})}}\leq x\right)-\Phi(x)\right|\leq\frac{C\tilde{E}|Z_1^{(1)}-p_1|^3}{\widetilde{Var}^{3/2}(Z_1^{(1)})\sqrt{n}}
\end{equation*}
where $\Phi$ is the CDF of standard normal distribution and $C$ is a universal constant. Let 
\begin{equation*}
x=z_{1-\alpha/2}\left(\sqrt{\frac{\tE \hat{V}_1}{\widetilde{Var}(Z_1^{(1)})}}-\sqrt{\frac{2\log(1/\delta)}{(n-1)\widetilde{Var}(Z_1^{(1)})}}ke^{-I(a_1)}\right)-p_2\sqrt{\frac{n}{\widetilde{Var}(Z_1^{(1)})}}.
\end{equation*}
Then we get that 
\begin{equation*}
\tilde{P}\left(\frac{\sqrt{n}|\hp_1-p_1|}{\sqrt{\widetilde{Var}(Z_1^{(1)})}}> x\right)\leq 2\Phi(-x)+\frac{2C\tilde{E}|Z_1^{(1)}-p_1|^3}{\widetilde{Var}^{3/2}(Z_1^{(1)})\sqrt{n}}.
\end{equation*}
Hence,
\begin{equation*}
\tilde{P}\left(|\hat p-p|>z_{1-\alpha/2}\sqrt{\frac{\hat{V}}{n}}\right)\leq 2\Phi(-x)+\frac{2C\tilde{E}|Z_1^{(1)}-p_1|^3}{\widetilde{Var}^{3/2}(Z_1^{(1)})\sqrt{n}}+\delta+n\tilde{p}_2.
\end{equation*}
First, we have that $\tilde{E}|Z_1^{(1)}-p_1|^3\leq \max(p_1^3,(ke^{-I(a_1)}-p_1)^3)\sim e^{-3I(a_1)}$. Since Assumption \ref{asm:BE_error} holds, we get that $\frac{\tilde{E}|Z_1^{(1)}-p_1|^3}{\widetilde{Var}^{3/2}(Z_1^{(1)})}$ grows at most subexponentially in $-\log p$. Hence under the assumptions we could choose $n$ as required in the theorem. In particular, we have that 
\begin{equation*}
\frac{2C\tilde{E}|Z_1^{(1)}-p_1|^3}{\widetilde{Var}^{3/2}(Z_1^{(1)})\sqrt{n}}\to 0.
\end{equation*}
Now we analyze $x$. We know that $\sqrt{\frac{E\hat{V}_1}{\widetilde{Var}(Z_1^{(1)})}}=1$ and 
\begin{equation*}
p_2\sqrt{\frac{n}{\widetilde{Var}(Z_1^{(1)})}}=\frac{p_2}{e^{-I(a_1)}}\sqrt{\frac{ne^{-2I(a_1)}}{\widetilde{Var}(Z_1^{(1)})}}\to 0
\end{equation*}
since $\frac{p_2}{e^{-I(a_1)}}$ decays exponentially (Assumption \ref{asm:p2_negligible}) while $n$ and $\frac{e^{-2I(a_1)}}{\widetilde{Var}(Z_1^{(1)})}$ grow subexponentially (Assumption \ref{asm:BE_error}) in $-\log p$. Now we consider
\begin{equation*}
\sqrt{\frac{2\log(1/\delta)}{(n-1)\widetilde{Var}(Z_1^{(1)})}}ke^{-I(a_1)}=\sqrt{\frac{2k^2\log(1/\delta)e^{-2I(a_1)}}{(n-1)\widetilde{Var}(Z_1^{(1)})}}.
\end{equation*}
Since we assume that $\frac{k^2e^{-2I(a_1)}}{n\tVar(Z_1^{(1)})}\to 0$, we could set $\delta$ such that $\delta\to 0$ and 
\begin{equation*}
\sqrt{\frac{2\log(1/\delta)}{(n-1)\widetilde{Var}(Z_1^{(1)})}}ke^{-I(a_1)}\to 0.
\end{equation*}
In this case, $x\to z_{1-\alpha/2}$ and hence $\Phi(-x)\to \alpha/2$. Combining all the results, we get that 
\begin{equation*}
\limsup_{\gamma\to\infty}\tilde{P}\left(|\hat p-p|>z_{1-\alpha/2}\sqrt{\frac{\hat{V}}{n}}\right)\leq\alpha.
\end{equation*}\hfill\Halmos
\endproof

\proof{Proof of Theorem \ref{thm:gaussian_assumptions}.}
First of all, it is easy to verify that the cumulant generating function $\mu(x)=\lambda^\top x+\frac12x^\top \Sigma x$ satisfies Assumption \ref{asm:mu}. Moreover, since $g$ is a piecewise linear function, we can express $\cE$ as the union of finite closed polyhedrons. If $P(X\in\cE)>0$ and $\lambda\notin\cE$, then $0<\inf_{x\in\cE}\frac12(x-\lambda)^\top \Sigma^{-1}(x-\lambda)<\infty$. From \citet{bai2022rare}, we know that
\begin{equation*}
    p=P(g(X)\geq \gamma)\sim e^{-\frac12(a_1-\lambda)^\top \Sigma^{-1}(a_1-\lambda)}
\end{equation*}
and for sufficient large $\gamma$, $(a_1-\lambda)^\top \Sigma^{-1}(a_1-\lambda)$ is a quadratic function of $\gamma$ which goes to $\infty$ as $\gamma\to\infty$. Finally, it is clear that the number of dominating points will not grow exponentially in $-\log p=\Theta((a_1-\lambda)^\top \Sigma^{-1}(a_1-\lambda))$. Therefore, Assumption \ref{asm:problem} is satisfied in this problem setting.

Next, we check Assumptions \ref{asm:p2_negligible} and \ref{asm:p2tilde}. Without loss of generality, we may assume that $k<r$, since otherwise $\cE_2=\emptyset$ and $p_2=\tp_2=0$. We know that $a_{k+1}=\arg\min_{x\in\cE_2}(x-\lambda)^\top \Sigma^{-1}(x-\lambda)$, and hence
\begin{equation*}
    \cE_2\subset\{x\in\R^d:(x-\lambda)^\top \Sigma^{-1}(x-\lambda)\geq (a_{k+1}-\lambda)^\top \Sigma^{-1}(a_{k+1}-\lambda)\}.
\end{equation*}
Denote $X'=\Sigma^{-1/2}(X-\lambda)$ and then $X'\sim N(0,I_d)$ under $P$. We have that
\begin{align*}
    p_2&\leq P((X-\lambda)^\top \Sigma^{-1}(X-\lambda)\geq (a_{k+1}-\lambda)^\top \Sigma^{-1}(a_{k+1}-\lambda))\\
    &=P(X^{'T}X'\geq (a_{k+1}-\lambda)^\top \Sigma^{-1}(a_{k+1}-\lambda))\\
    &\sim e^{-\frac12(a_{k+1}-\lambda)^\top \Sigma^{-1}(a_{k+1}-\lambda)}.
\end{align*}
We know that $(a_{k+1}-\lambda)^\top \Sigma^{-1}(a_{k+1}-\lambda)>C(a_k-\lambda)^\top \Sigma^{-1}(a_k-\lambda)\geq C(a_1-\lambda)^\top \Sigma^{-1}(a_1-\lambda)$. Thus, Assumption \ref{asm:p2_negligible} holds. Moreover, 
\begin{align*}
    \tp_2&=\frac{1}{k}\sum_{i=1}^kP_{X\sim N(a_i,\Sigma)}(X\in\cE_2)\\
    &\leq \frac{1}{k}\sum_{i=1}^kP_{X\sim N(a_i,\Sigma)}((X-\lambda)^\top \Sigma^{-1}(X-\lambda)\geq (a_{k+1}-\lambda)^\top \Sigma^{-1}(a_{k+1}-\lambda)).
\end{align*}
We have that 
\begin{align*}
    &(X-a_i)^\top \Sigma^{-1}(X-a_i)\\
    =&(X-\lambda)^\top \Sigma^{-1}(X-\lambda)+(a_i-\lambda)^\top \Sigma^{-1}(a_i-\lambda)-2(a_i-\lambda)^\top \Sigma^{-1}(X-\lambda)\\
    \geq &(X-\lambda)^\top \Sigma^{-1}(X-\lambda)+(a_i-\lambda)^\top \Sigma^{-1}(a_i-\lambda)-2\sqrt{(X-\lambda)^\top \Sigma^{-1}(X-\lambda)}\sqrt{(a_i-\lambda)^\top \Sigma^{-1}(a_i-\lambda)}\\
    =&\left(\sqrt{(X-\lambda)^\top \Sigma^{-1}(X-\lambda)}-\sqrt{(a_i-\lambda)^\top \Sigma^{-1}(a_i-\lambda)}\right)^2.
\end{align*}
Hence, if $(X-\lambda)^\top \Sigma^{-1}(X-\lambda)\geq (a_{k+1}-\lambda)^\top \Sigma^{-1}(a_{k+1}-\lambda)>C(a_k-\lambda)^\top \Sigma^{-1}(a_k-\lambda)\geq C(a_i-\lambda)^\top \Sigma^{-1}(a_i-\lambda)$, then $(X-a_i)^\top \Sigma^{-1}(X-a_i)\geq (\sqrt{C}-1)^2(a_i-\lambda)^\top \Sigma^{-1}(a_i-\lambda)$. Thus we get that 
\begin{align*}
    \tp_2&\leq\frac{1}{k}\sum_{i=1}^kP_{X\sim N(a_i,\Sigma)}((X-a_i)^\top \Sigma^{-1}(X-a_i)\geq (\sqrt{C}-1)^2(a_i-\lambda)^\top \Sigma^{-1}(a_i-\lambda))\\
    &\leq \frac{1}{k}\sum_{i=1}^kP_{X\sim N(a_i,\Sigma)}((X-a_i)^\top \Sigma^{-1}(X-a_i)\geq (\sqrt{C}-1)^2(a_1-\lambda)^\top \Sigma^{-1}(a_1-\lambda))\\
    &=P(X^{'T}X'\geq (\sqrt{C}-1)^2(a_1-\lambda)^\top \Sigma^{-1}(a_1-\lambda))\\
    &\sim e^{-\frac{1}{2}(\sqrt{C}-1)^2(a_1-\lambda)^\top \Sigma^{-1}(a_1-\lambda)}.
\end{align*}
Hence, $\tp_2$ exponentially decays in $-\log p$. That is, Assumption \ref{asm:p2tilde} holds.

Finally, we verify Assumption \ref{asm:BE_error}. We have $\tVar(Z_1)=\tE(Z_1-p_1)^2\geq\tE((Z_1-p_1)^2I_{\cE_1^c}(X))=p_1^2\tP(X\notin \cE_1)$. Since $p_1\sim e^{-\frac12(a_1-\lambda)^\top \Sigma^{-1}(a_1-\lambda)}$, it suffices to justify that $\tP(X\notin\cE_1)$ does not decay exponentially in $-\log p$. Indeed, we have that $\cE_1\subset\{x\in\R^d:(x-\lambda)^\top \Sigma^{-1}(x-\lambda)\geq (a_1-\lambda)^\top \Sigma^{-1}(a_1-\lambda)\}$ and thus
\begin{align*}
    \tP(X\in\cE_1)&=\frac{1}{k}\sum_{i=1}^kP_{X\sim N(a_i,\Sigma)}(X\in\cE_1)\\
    &\leq\frac{k-1}{k}+\frac{1}{k}P_{X\sim N(a_1,\Sigma)}(X\in\cE_1)\\
    &\leq\frac{k-1}{k}+\frac{1}{k}P_{X\sim N(a_1,\Sigma)}((X-\lambda)^\top \Sigma^{-1}(X-\lambda)\geq (a_1-\lambda)^\top \Sigma^{-1}(a_1-\lambda))\\
    &=\frac{k-1}{k}+\frac{1}{k}P((X'+\Sigma^{-1/2}(a_1-\lambda))^\top (X'+\Sigma^{-1/2}(a_1-\lambda))\geq (a_1-\lambda)^\top \Sigma^{-1}(a_1-\lambda)).
\end{align*}
Then we have that 
\begin{align*}
    \tP(X\notin\cE_1)&\geq \frac{1}{k}P((X'+\Sigma^{-1/2}(a_1-\lambda))^\top (X'+\Sigma^{-1/2}(a_1-\lambda))< (a_1-\lambda)^\top \Sigma^{-1}(a_1-\lambda))\\
    &\geq\frac{1}{r}P\left(X'\in B\left(-\Sigma^{-1/2}(a_1-\lambda),\sqrt{(a_1-\lambda)^\top \Sigma^{-1}(a_1-\lambda)}\right)\right)
\end{align*}
where $B(x,R):=\{y\in\R^d:\|y-x\|<R\}$. For sufficiently large $\gamma$, $(a_1-\lambda)^\top \Sigma^{-1}(a_1-\lambda)$ monotonely grows to $\infty$, and hence $\tP(X\notin\cE_1)\geq \frac{c}{r}$ for some constant $c>0$. As a result, Assumption \ref{asm:BE_error} holds.\hfill\Halmos
\endproof

\end{APPENDICES}

\end{document}